\documentclass[11pt]{article}   
\usepackage{amssymb,latexsym}
\usepackage{epsfig}
\usepackage{eufrak}
\usepackage{amsmath}
\usepackage{mathrsfs}
\usepackage{color}

\newtheorem{theorem}{Theorem}
\newtheorem{lemma}{Lemma}
\newtheorem{corollary}{Corollary}

\newcommand{\be}{\begin{equation}}
\newcommand{\ee}{\end{equation}}
\newcommand{\bee}{\begin{eqnarray*}}
\newcommand{\eee}{\end{eqnarray*}}
\newcommand{\bel}{\begin{eqnarray}}
\newcommand{\eel}{\end{eqnarray}}
\newcommand{\bec}{\begin{cases}}
\newcommand{\eec}{\end{cases}}
\newcommand{\bem}{\begin{bmatrix}}
\newcommand{\eem}{\end{bmatrix}}

\newcommand{\la}{\label}
\newcommand{\li}{\left}
\newcommand{\ri}{\right}

\newcommand{\DEF}{\stackrel{\mathrm{def}}{=}}

\newcommand{\ovl}{\overline}
\newcommand{\udl}{\underline}

\newcommand{\lc}{\lceil}
\newcommand{\rc}{\rceil}
\newcommand{\lf}{\lfloor}
\newcommand{\rf}{\rfloor}

\newcommand{\ep}{\epsilon}
\newcommand{\vep}{\varepsilon}
\newcommand{\lm}{\lambda}
\newcommand{\Lm}{\Lambda}
\newcommand{\Up}{\Upsilon}

\newcommand{\si}{\sigma}

\newcommand{\de}{\delta}

\newcommand{\vDe}{\varDelta}

\newcommand{\ga}{\gamma}
\newcommand{\Ga}{\Gamma}
\newcommand{\vse}{\vartheta}
\newcommand{\se}{\theta}
\newcommand{\Se}{\Theta}

\newcommand{\ze}{\zeta}
\newcommand{\al}{\alpha}
\newcommand{\ba}{\beta}

\newcommand{\vro}{\varrho}
\newcommand{\ro}{\rho}
\newcommand{\ka}{\kappa}
\newcommand{\om}{\omega}
\newcommand{\Om}{\Omega}

\newcommand{\f}{\frac}
\newcommand{\sq}{\sqrt}
\newcommand{\cd}{\cdots}

\newcommand{\qu}{\quad}
\newcommand{\qqu}{\qquad}
\newcommand{\fa}{\forall}

\newcommand{\mscr}{\mathscr}
\newcommand{\mcal}{\mathcal}
\newcommand{\mbf}{\mathbf}
\newcommand{\bb}{\mathbb}
\newcommand{\fra}{\mathfrak}

\newcommand{\wh}{\widehat}
\newcommand{\wt}{\widetilde}
\newcommand{\mrm}{\mathrm}
\newcommand{\bs}{\boldsymbol}

\newcommand{\ap}{\approx}

\newcommand{\sh}{\slash}

\newcommand{\tx}{\text}

\newcommand{\iy}{\infty}

\newcommand{\bed}{\begin{description}}
\newcommand{\eed}{\end{description}}
\newcommand{\bei}{\begin{itemize}}
\newcommand{\eei}{\end{itemize}}
\newcommand{\ben}{\begin{enumerate}}
\newcommand{\een}{\end{enumerate}}
\newcommand{\bib}{\bibitem}
\newcommand{\beL}{\begin{lemma}}
\newcommand{\eeL}{\end{lemma}}
\newcommand{\beT}{\begin{theorem}}
\newcommand{\eeT}{\end{theorem}}
\newcommand{\beC}{\begin{corollary}}
\newcommand{\eeC}{\end{corollary}}
\newcommand{\sect}{\section}

\newcommand{\bpf}{\begin{pf}}
\newcommand{\epf}{\end{pf}}
\newcommand{\bsk}{\bigskip}
\newcommand{\bi}{\binom}

\newcommand{\pfbox}{\hfill\mbox{$\Box$}}

\newenvironment{pf}{\paragraph*{Proof{\rm.}}}{\pfbox\bigskip}

\begin{document}

\title{{\bf A New Framework of Multistage Hypothesis Tests}
\thanks{The author had been previously working with
Louisiana State University at Baton Rouge, LA 70803, USA, and is now with Department of Electrical Engineering, Southern University and A\&M
College, Baton Rouge, LA 70813, USA; Email: chenxinjia@gmail.com. The main results of this paper have been presented in Proceedings of SPIE
Conferences, Orlando, April 5-9, 2010 and April 25-29, 2011.  The statistical methodology proposed in this paper has been applied to electrical
engineering and computer science, see recent literature \cite{Chen_SPIE8401, Chen_SPIE,  Chen_SPIE11a, Chen_SPIE11b, Chen_SPIE11c, Chen_J} and
the references therein.} }

\author{Xinjia Chen}

\date{First submitted in September 2008}

\maketitle

$\qqu \qqu  \qqu $ {\it In Memory of My Dear Father Hualong Chen
(1933--1990)}

\begin{abstract}

In this paper, we have established a general framework of multistage
hypothesis tests which applies to arbitrarily many mutually
exclusive and exhaustive composite hypotheses. Within the new
framework, we have constructed specific multistage tests which
rigorously control the risk of committing decision errors and are
more efficient than previous tests in terms of average sample number
and the number of sampling operations. Without truncation, the
sample numbers of our testing plans are absolutely bounded.

\end{abstract}

\tableofcontents

\sect{Introduction}

Statistical inference based on samples drawn from populations can be viewed as a rationale approach of making decisions based on observations of
stochastic processes.  Consider $\ka$ continuous-time stochastic processes $(X_t^\nu)_{t \in \bb{R}}, \; \nu = 1, \cd, \ka$ defined in a
probability space $(\Om, \mscr{F}, \Pr )$. Suppose that the probability measure $\Pr$ is determined by $q$ parameters $\vse_1, \cd, \vse_q$ in
the sense that for any $t \in \bb{R}$, the joint distribution of $X_t^\nu, \; \nu = 1, \cd, \ka$ is parameterized by $\vse_1, \cd, \vse_q$. Let
$\se$ be a function of $\vse_1, \cd, \vse_q$.  Let the set of all values of $\se$ be denoted by $\Se$. In many areas of engineering and
sciences, it is desirable to infer the true value of $\se$ based on the observation of such stochastic processes.  This topic can be formulated
as a general problem of testing $m$ mutually exclusive and exhaustive composite hypotheses: \be \la{mainpr} \mscr{H}_0: \se \in {\Se}_0, \qu
\mscr{H}_1: \se \in {\Se}_1, \qu \ldots, \qu \mscr{H}_{m - 1}: \se \in {\Se}_{m - 1}, \ee where $\Se_i = \{ \se \in \Se: \se_i < \se \leq
\se_{i+1} \}, \; i = 0, 1, \cd, m - 1$ with $- \iy = \se_0 < \se_1 < \cd < \se_{m - 1} < \se_m = \iy$. To control the probabilities of making
wrong decisions, for pre-specified numbers $\de_i \in (0, 1), \; i = 0, 1, \cd, m - 1$, it is typically required that \be \la{mainreq}
 \Pr \{ \tx{Reject} \; \mscr{H}_i \mid \se \} \leq \de_i, \qqu \fa \se
\in \varTheta_i, \qu i = 0, 1, \cd, m - 1 \ee  where $\varTheta_i = \{ \se \in \Se_i: \se_i^{\prime \prime} \leq \se \leq \se_{i+1}^{\prime} \},
\; i = 0, 1, \cd, m - 1$ with $\se_i^\prime, \se_i^{\prime \prime} \in \Se, \; i = 1, \cd, m - 1$ satisfying $- \iy = \se_0^{\prime \prime} <
\se_i^\prime < \se_i < \se_i^{\prime \prime} \leq \se_{i + 1}^\prime < \se_{i + 1} < \se_{i+1}^{\prime \prime} < \se_m^\prime = \iy$ for $i = 1,
\cd, m - 2$.  Here we consider continuous-time processes for the sake of generality, since discrete-time stochastic processes can be treated as
right-continuous processes in continuous time.   For $i = 0, 1, \cd, m - 1$, $\Pr \{ \tx{Accept} \; \mscr{H}_i \mid \se \}$ is referred to as an
Operating Characteristic (OC) function.  Since there is no requirement imposed for controlling the risk of making wrong decisions for $\se$ in
${\Se} \setminus \cup_{j = 0}^{m - 1} \varTheta_j = \cup_{i = 1}^{m-1} (\se_i^\prime, \se_i^{\prime \prime})$, such a remainder set, $\cup_{i =
1}^{m-1} (\se_i^\prime, \se_i^{\prime \prime})$, is referred to as an {\it indifference zone}. The concept of indifference zone was introduced
by Wald \cite{Wald} for two main reasons.  First, when the parameter $\se$ is close to $\se_i$, the margin between adjacent parameter subsets
$\Se_{i-1}$ and $\Se_i$, it is immaterial to decide whether $\mscr{H}_{i-1}$ or $\mscr{H}_i$ should be accepted. Second, the sample size
required to make a reliable decision between consecutive hypotheses $\mscr{H}_{i-1}$ and $\mscr{H}_i$ becomes increasingly intolerable as $\se$
tends to $\se_i$.  Undoubtedly, the indifference zone should be sufficiently ``narrow'' so that the consequence of making erroneous decision is
practically unimportant when $\se$ lies in it.

The general problem of hypothesis testing  described above has been
a fundamental issue of research for many decades. The well-known
sequential probability ratio test (SPRT) has been developed by Wald
\cite{Wald} to address such testing problem in the special case of
two hypotheses.  In addition to the limitation associated with the
number of hypotheses, the SPRT
 suffers from other drawbacks.  First, the sample number of SPRT
 is a random number which is not bounded. However, to be useful, the sample number of
 any testing plan should be bounded by a deterministic number.  Although this can be fixed
 by forced termination (see, e.g., \cite{Ghosh} and the references therein),
 the prescribed level of power may not be
 ensured as a result of truncation.  Second, the number of sampling
 operations of SPRT is as large as the number of samples.  In practice, it
 is usually much more economical to take a batch of samples at a
 time instead of one by one.  Third, the efficiency of SPRT is
 optimal only for the endpoints of the indifference zone. For other parametric values, the SPRT
  can be extremely inefficient.  Needless to say, a truncated version of SPRT may suffer from the same
 problem due to the partial use of the boundary of SPRT.

 In this paper, to overcome the limitations of SPRT and its
 variations, we have established a new framework of hypothesis testing which applies to arbitrary number of composite hypotheses.
 Our testing plans have the following features:
 i) The testing process has a finite number of stages and thus the cost of
 sampling operations is reduced as compared to SPRT; ii) The sample number is absolutely
 bounded without truncation; iii) The prescribed level of
 power is rigorously guaranteed;  iv) The testing is not only efficient for
 the endpoints of indifference zone, but also efficient for other parametric values.  The
 remainder of the paper is organized as follows.  In Section 2, we present a unified approach for multi-valued decision
 in the general framework of constructing sequential random intervals with pre-specified coverage probabilities. In Section 3, we
 present our general theory and computational mechanisms for the design and analysis of multistage testing
 plans. In Section 4, we first present more specific construction of testing procedures and then apply the general method to
 common important problems.  Specially, we demonstrate that the
 principle can be used for testing a binomial proportion,
  the proportion of a finite population, a Poisson parameter, the mean of a normal distribution with known variance,
  the variance of a normal distribution, the parameter of an exponential
 distribution, the scale parameter of a Gamma distribution and life testing.
 Section 5 is dedicated to tests of the mean of a normal distribution with unknown
 variance.  Section 6 addressed the problem of testing multiple hypotheses
 regarding the ratio of variances of two normal distributions.  In
 Section 7,  we have established an exact computational method of the OC function and average sample number
of the SPRT. Such computational method can be used to compare the performance of SPRTs with our tests.  In Section 8, we propose an exact and
efficient recursive method for computing boundary crossing probabilities, which can be applied to evaluate the risks of making incorrect
decisions in multistage hypothesis testing. Section 9 is the conclusion. All proofs of theorems are given in Appendices.

Throughout this paper, we shall use the following notations.  The
notation $\emptyset$ denotes an empty set. The set of real numbers
is denoted by $\bb{R}$. The set of integers is denoted by $\bb{Z}$.
The set of positive integers is denoted by $\bb{N}$.   The ceiling
function and floor function are denoted respectively by $\lc . \rc$
and $\lf . \rf$ (i.e., $\lc x \rc$ represents the smallest integer
no less than $x$; $\lf x \rf$ represents the largest integer no
greater than $x$). The gamma function is denoted by $\Ga(.)$. For
any integer $i$, the combinatoric function $\bi{i}{j}$ with respect
to integer $j$ takes value {\small $\f{ \Ga( i + 1) } { \Ga( j + 1)
\Ga (i - j + 1) }$} for $j \leq i$ and value $0$ otherwise. The
expectation of a random variable is denoted by $\bb{E}[.]$. We use
the notation $\Pr \{ . \mid \se \}$ to denote the probability of an
event which is defined in terms of random variables parameterized by
$\se$. The parameter $\se$ in $\Pr \{ . \mid \se \}$  may be dropped
whenever this can be done without introducing confusion.  If $Z$ is
parameterized by $\se$, we denote $\Pr \{  Z \leq z \mid \se \}$ by
$F_Z(z, \se)$ and $\Pr \{ Z \geq z \mid \se \}$ by $G_Z(z, \se)$
respectively.   The cumulative distribution function of a Gaussian
random variable is denoted by $\Phi(.)$. For $\al \in (0, 1)$,
$\mcal{Z}_\al$ denotes the critical value satisfying
$\Phi(\mcal{Z}_\al) = 1 - \al$. For $\al \in (0, 1)$, let $\chi_{n,
\al}^2$ denote the $100 \al \%$ percentile of a chi-square
distribution of $n$ degrees of freedom. For $\al \in (0, 1)$, let
$t_{n, \al}$ denote the $100(1 - \al) \%$ percentile of a Student
$t$-distribution of $n$ degrees of freedom. The support of a random
variable $Z$ is denoted by $I_Z$, i.e., $I_Z = \{ Z(\om): \om \in
\Om \}$.  We write $\bs{\de} = O(\ze)$ if $\bs{\de}$ is a function
of $\ze > 0$ such that there exist constants $A$ and $B$ such that $
A < \f{\bs{\de}}{\ze} < B$ provided that $\ze > 0$ is sufficiently
small.  The other notations will be made clear as we proceed.

\section{Sequential Random Intervals and Multi-hypotheses Testing} \la{secSRI}

As demonstrated in \cite{INCP}, the general hypothesis testing problem defined by (\ref{mainpr}) and (\ref{mainreq}) can be cast into the
framework of constructing a sequential random interval with pre-specified coverage probabilities. The methodology of \cite{INCP} is represented
in the sequel.

To reach a fast decision, it is desirable to solve the above problem by a multistage approach such that the sampling procedure is divided into
$s$ stages with observational times $t_{\nu, \ell}, \; \nu = 1, \cd, \ka; \; \ell = 1, \cd, s$, where $t_{\nu, \ell}$ is the observational time
for the $i$-th process at the $\ell$-th stage.  Starting from $\ell = 1$, at the $\ell$-th stage, based on the observation of $(X_t^\nu)_{0 \leq
t \leq t_{\nu, \ell} }, \; \nu = 1, \cd, \ka$,  pre-determined stopping and decision rules are applied to check whether the accumulated
observational data is sufficient to accept a hypothesis and terminate the sampling procedure. If the observational data is considered to be
insufficient for making a decision, then proceed to the next stage of observation.  The observation is continued stage by stage until a
hypothesis is accepted at some stage. Although the number of stages $s$ may be infinity, for practical considerations, the stopping and decision
rules are required to guarantee that the sampling procedure is surely terminated with a finite number of stages. Central to a multistage
procedure are the stopping and decision rules, which can be related to a sequential random interval described as follows.  Let $\se_0^\prime =
-\iy$ and $\se_m^{\prime \prime} = \iy$. For $i = 0, 1, \cd, m - 1$, let $\mscr{I}_i$ denote the open interval {\small $( \se_i^\prime,
\se_{i+1}^{\prime \prime} )$}. Let $\bs{l}$ be the index of stage at the termination of the sampling procedure. Let $\bs{\mcal{L}}$ and
$\bs{\mcal{U}}$ be random variables defined in terms of samples of the $\ka$ stochastic processes up to the $\bs{l}$-th stage such that the
sequential random interval {\small $( \bs{\mcal{L}}, \bs{\mcal{U}})$} has $m$ possible outcomes $\mscr{I}_i, \; i = 0, 1, \cd, m - 1$ and that
$\Pr \{ \bs{\mcal{L}}  < \se < \bs{\mcal{U}}  \mid \se \} > 1 - \de_i$ for any $\se \in \varTheta_i$ and $i = 0, 1, \cd, m - 1$.  Given that the
sequential random interval $( \bs{\mcal{L}}, \bs{\mcal{U}} )$ satisfying such requirements is constructed, the risk requirement (\ref{mainreq})
can be satisfied by using $( \bs{\mcal{L}}, \bs{\mcal{U}} )$ to define a decision rule such that, for $i = 0, 1, \cd, m - 1$, hypothesis
$\mscr{H}_i$ is accepted when the sequential random interval $( \bs{\mcal{L}}, \bs{\mcal{U}} )$ takes $\mscr{I}_i$ as its outcome at the
termination of the sampling process. It follows that $\{ \tx{Accept} \; \mscr{H}_i \} =\{ \bs{\mcal{L}} < \se < \bs{\mcal{U}} \}$ for any $\se
\in \varTheta_i$ and $i = 0, 1, \cd, m - 1$. Therefore, to solve the multi-valued decision problem defined by (\ref{mainpr}) and
(\ref{mainreq}), the objective is to ensure that $\se$ is included in the sequential random interval with pre-specified probabilities.  In the
sequel, we shall propose a general approach for defining stopping and decision rules for the construction of such sequential random interval.

\subsection{General Structure of Stopping and Decision Rules}

We shall first consider the structure of stopping and decision rules under the assumption that the number of stages, $s$, and the observational
times, $t_{\nu, \ell}, \; \nu = 1, \cd, \ka; \; \ell = 1, \cd, s$, are given.  The determination of stage number, observational times and the
issue of finite stopping will be addressed later in the parametrization of the stopping and decision rules. We propose to use one-sided
confidence sequences to control the coverage probability of the sequential random interval.  Assume that for $\ell = 1, \cd, s$ and $i = 1, \cd,
m-1$, random variables $L_{\ell, i}$ and $U_{\ell, i}$ can be defined in terms of positive numbers $\ze, \; a_i, \; b_i$ and the set of random
variables $(X_t^\nu)_{0 \leq t \leq t_{\nu, \ell} }, \; \nu = 1, \cd, \ka$ such that $\Pr \{ L_{\ell, i} \geq \se \mid \se \}$ and $\Pr \{
U_{\ell, i} \leq \se \mid \se \}$  can be made arbitrarily small by decreasing $\ze a_i$ and $\ze b_i$ respectively. Due to such assumption, we
call $(-\iy, L_{\ell, i}]$ and $[U_{\ell, i}, \iy)$ one-sided confidence intervals for $\se$. Accordingly, $(-\iy, L_{\ell, i}], \; \ell = 1,
\cd, s$ and $[U_{\ell, i}, \iy), \; \ell = 1, \cd, s$ are said to be one-sided confidence sequences for $\se$. In view of the controllability of
the coverage probabilities of the one-sided confidence intervals, the number $\ze$ is referred to as the {\it coverage tuning parameter}, and
$a_i, \; b_i, \; i = 1, \cd, m-1$ are called {\it weighting coefficients}.  Given that $\ze$ is sufficiently small, $\se
> \se_i^\prime$ will be credible if $L_{\ell, i} > \se_i^\prime$ is observed. Similarly, $\se < \se_i^{\prime \prime}$ will be credible if
$U_{\ell, i} < \se_i^{\prime \prime}$ is observed.  To figure out the general structure of stopping and decision rules, imagine that the
sampling procedure is stopped at the $\ell$-th stage and $\mscr{I}_i$ is to be designated as the outcome of the sequential random interval.
Since $\mscr{I}_i$ contains $[\se_i^{\prime \prime}, \se_{i + 1}^{\prime}]$, it follows that for $\se \in [\se_i^{\prime \prime}, \se_{i +
1}^{\prime}]$, it is true that $\se < \se_{i+1}^{\prime \prime}$ and $\se > \se_i^\prime$.  This implies that, if the coverage tuning parameter
$\ze$ is sufficiently small, then it is very likely to observe that $U_{\ell, i+1} < \se_{i+1}^{\prime \prime}$ and $L_{\ell, i} >
\se_i^\prime$.  Therefore, turning this thinking around leads to the following stopping and decision rules:

\vspace{0.05in}

\begin{tabular} {|l |}
\hline $ \tx{ {\it Continue observing the stochastic processes until for some $i \in \{0, 1, \cd, m - 1\}$,}}$\\
$\tx{ {\it the event $\{ U_{\ell, i+1} < \se_{i+1}^{\prime \prime} \; \tx{and} \; L_{\ell, i} > \se_i^\prime
\}$ occurs at some stage with index $\ell \in \{1, \cd, s\}$}}$.\\
$\tx{ {\it At the termination of the sampling process, make the following decision: If such index}}$\\
$ \tx{ {\it  $i$ is unique,  then designate $\mscr{I}_i$ as the outcome of the sequential random interval. If there}}$\\
$ \tx{ {\it  are multiple indexes satisfying the condition, then pick one of them and assign the}}$\\
$ \tx{ {\it  corresponding interval $\mscr{I}_i$ as the outcome of the sequential random interval based on}}$\\
$ \tx{ {\it  a predetermined policy.}}$
\\ \hline
\end{tabular}

\vspace{0.05in}

 The idea in the derivation of the above stopping and decision rules is to infer
 the location of $\se$ relative to the sequential random interval by comparing
 the confidence limits with the endpoints of the sequential random interval. Due to
the probabilistic nature of the comparison, such method of constructing stopping and decision rules is referred to as the {\it Principle of
Probabilistic Comparison}.  It should be noted that similar principles have been proposed in \cite{Chen_EST} for multistage estimation of
parameters.  Now consider the simplification of the above stopping and decision rules.  For $\ell = 1, \cd, s$, let $L_{\ell, 0} = - \iy$ and
$U_{\ell, m} = \iy$.  For $i = 0, 1, \cd, m - 1$, let $(L_{\ell, i}, U_{\ell, i+1}), \; \ell = 1, \cd, s$ be referred to as the $i$-th {\it
controlling confidence sequence}, where the word ``controlling'' is used to indicate that the confidence sequences are used to control the
coverage probability of the desired sequential random interval.  Assume that for $\ell = 1, \cd, s$, there exists a statistic $\wh{\se}_\ell$
dependent on random variables $(X_t^\nu)_{0 \leq t \leq t_{\nu, \ell} }, \; \nu = 1, \cd, \ka$ such that $L_{\ell, i} \leq \wh{\se}_\ell \leq
U_{\ell, i}$ for $i = 1, \cd, m - 1$.   In this setting, we propose the following stopping and decision rules: \vspace{0.05in}

\begin{tabular} {|l |}
\hline

$ \tx{ {\it Continue the sampling process if there exists no index $i\in \{0, 1, \cd, m - 1\}$ such that}}$\\
$ \tx{ {\it  the $i$-th controlling confidence sequence is included by interval $\mscr{I}_i$. At the termination}}$\\
$ \tx{ {\it   of the sampling process, make the following decision}}$:\\
$ \tx{ {\it (a):  If there exists a unique index $i \in \{0, 1, \cd, m - 1\}$ such that the $i$-th  controlling} }$\\
$ \qu \qu \; \tx{{\it  confidence sequence is included by interval $\mscr{I}_i$, then designate $\mscr{I}_i$ as the outcome}}$\\
$ \qu \qu \tx{ {\it  of the sequential random interval}}$.\\
$ \tx{ {\it (b):  If there exist two consecutive indexes $i - 1$ and $i$ in $\{1, \cd, m - 1\}$ such that the}}$\\
$ \qu \qu \tx{ {\it  $(i - 1)$-th and $i$-th controlling confidence sequences are included, respectively, by}}$\\
$ \qu \qu \tx{ {\it  intervals $\mscr{I}_{i-1}$ and $\mscr{I}_i$, then designate either $\mscr{I}_{i-1}$ or $\mscr{I}_i$
as the outcome of the}}$\\
$ \qu \qu \tx{ {\it  sequential random interval based on a pre-specified policy.}}$
\\ \hline
\end{tabular}

\vspace{0.05in}

In view of the inclusion relationship implemented in the termination conditions, this simplified method of constructing stopping and decision
rules is referred to as the {\it Inclusion Principle}. The properties of the above two types of stopping and decision rules are indicated by the
following probabilistic results.

\beT \la{Multi-Valued Inclusion Principle}

Let $A_0 = B_0 = - \iy, \; A_m = B_m = \iy$ and  $A_i < B_i \leq A_{i + 1} < B_{i + 1}$ for $i = 1, \cd, m - 2$.    Let $\varTheta_0 = (- \iy,
A_1], \; \varTheta_{m-1} = [B_{m - 1}, \iy )$ and $\varTheta_i = [B_i, A_{i + 1}]$ for $i = 1, \cd, m - 2$.  Let $(\Om, \mscr{F}, \{
\mscr{F}_\ell \}, \Pr )$ be a filtered space. Let $\bs{l}$ be a proper stopping time with support $I_{\bs{l}}$.  For $\ell \in I_{\bs{l}}$, let
$L_{\ell, m } = - \iy, \; U_{\ell, 0}  = \iy$ and let $L_{\ell, i}, \; U_{\ell, i}, \; i = 1,\cd, m - 1$ be random variables measurable in
$\mscr{F}_\ell$.  Let $\bs{\mcal{L}}$ and $\bs{\mcal{U}}$ be random variables such that $\cup_{i = 0}^{m-1} \{ \bs{\mcal{L}} = A_i, \;
\bs{\mcal{U}} = B_{i + 1} \} = \Om$ and that {\small $\{ \bs{l} = \ell, \; \bs{\mcal{L}} = A_j, \; \bs{\mcal{U}} = B_{j + 1} \} \subseteq \{
L_{\ell, j} \geq A_j \; \tx{and} \; U_{\ell, j+1} \leq B_{j+1} \}$ } for $\ell \in I_{\bs{l}}$ and $j = 0, 1, \cd, m - 1$.  Then,  $\Pr \{
\bs{\mcal{L}} \geq \se \} = \Pr \{ \bs{\mcal{L}} \geq A_{i+1} \} \leq \sum_{j > i} \Pr \{ L_{\ell, j} \geq A_j \; \tx{for some} \; \ell \in
I_{\bs{l}} \}$ and $\Pr \{ \bs{\mcal{U}} \leq \se \} = \Pr \{ \bs{\mcal{U}} \leq B_i \} \leq \sum_{j \leq i} \Pr \{ U_{\ell, j} \leq B_j \;
\tx{for some} \; \ell \in I_{\bs{l}} \}$ for $i = 0, 1, \cd, m - 1$ and $\se \in \varTheta_i$.

\eeT

\subsection{Parametrization of Stopping and Decision Rules}

In this section, we shall develop a general approach for parameterizing multistage testing plans by virtue of the inclusion principle proposed
in the preceding section.  The objective is to express the testing plans in terms of the {\it coverage tuning parameter} $\ze > 0$ and positive
{\it weighting coefficients} $a_i, b_i, \; i = 1, \cd, m - 1$.  For eliminating unnecessary waste of sampling and operational effort, we shall
emphasis the principle that {\it a multistage test should stop at or before the last stage with probability $1$ and the test should have a
positive probability to stop at the first stage}. To deal with the differences of the rates of taking samples from the stochastic processes, we
introduce functions $\mcal{T}_\nu (\tau), \; \nu = 1, \cd, \ka$ which are increasing with respect to $\tau
> 0$ and tend to $\iy$ as $\tau \to \iy$. Such functions are referred to as {\it time transformation functions}.  Let $\al_i = \ze a_i \in (0,
1)$ and $\ba_i = \ze b_i \in (0, 1)$ for $i = 1, \cd, m - 1$.  To ensure that the sampling process  associated with the desired testing plan
will eventually terminate with probability $1$, we seek methods to construct random variables $\mscr{L} (\tau, \al_i), \; \mscr{U} (\tau,
\ba_i), \; i = 1, \cd, m - 1$ from random variables $(X_t^\nu)_{0 \leq t \leq \mcal{T}_\nu(\tau) }, \; \nu = 1, \cd, \ka$ such that $\cap_{i =
1}^{m-1} \{ \se_i^\prime < \mscr{L} (\tau, \al_i), \; \mscr{U} (\tau, \ba_i) < \se_i^{\prime \prime} \} \neq \emptyset$ if $\tau$ is
sufficiently large. Once this can be accomplished, we define $\tau^*$ as the minimum number $\tau > 0$ such that $\cap_{i = 1}^{m-1} \{
\se_i^\prime < \mscr{L} (\tau, \al_i), \; \mscr{U} (\tau, \ba_i) < \se_i^{\prime \prime} \} \neq \emptyset$ and $\tau_\star$ as the minimum
number $\tau > 0$ such that $\cup_{i = 0}^{m - 1} \{ \se_i^\prime < \mscr{L} (\tau, \al_i), \; \mscr{U} (\tau, \ba_{i+1}) < \se_{i+1}^{\prime
\prime} \} \neq \emptyset$, where we assume that $\{ \se_0^\prime < \mscr{L} (\tau, \al_0) \}$ and $\{ \mscr{U} (\tau, \ba_m) < \se_m^{\prime
\prime} \}$ are sure events. It follows that we can define the number of stages as some number $s$ and choose the observational times as
$t_{\nu, \ell} = \mcal{T}_\nu(\tau_\ell), \; \ell = 1, \cd, s$ with $\tau_1 < \tau_2 < \cd < \tau_s$ such that $\tau_\star \leq \tau_1 < \tau^*
\leq \tau_s$. Accordingly, the confidence limits can be parameterized as $L_{\ell, i} = \mscr{L} (\tau_\ell, \al_i), \; U_{\ell, i} = \mscr{U}
(\tau_\ell, \ba_i)$ for $\ell = 1, \cd, s$ and $i = 1, \cd, m-1$.  In many situations, it is possible that for $\ell = 1, \cd, s$, there exists
a statistic $\wh{\se}_\ell$ dependent on random variables $(X_t^\nu)_{0 \leq t \leq t_{\nu, \ell} }, \; \nu = 1, \cd, \ka$ such that $L_{\ell,
i} \leq \wh{\se}_\ell \leq U_{\ell, i}$ for $\ell = 1, \cd, s$ and $i = 1, \cd, m - 1$.

\subsection{Coverage Tuning} \la{BRT}

Given that multistage tests can be parameterized as in the preceding section, we need to determine appropriate values for the coverage tuning
parameter $\ze$ and weighting coefficients $a_i, b_i, \; i = 1, \cd, m - 1$ so that the sampling cost is as low as possible, while guaranteeing
 that the sequential random interval $(\bs{\mcal{L}}, \bs{\mcal{U}})$ satisfying the coverage specification:
 $\Pr \{ \bs{\mcal{L}} < \se < \bs{\mcal{U}} \mid \se \} \geq 1 - \de_i$ for any $\se \in \varTheta_i$ and $i = 0, 1, \cd, m - 1$.
The computational process for accomplishing this task is called {\it coverage tuning}. As a consequence of the connection between sequential
random intervals and hypothesis testing we previously established,  the specification for the coverage probability of $(\bs{\mcal{L}},
\bs{\mcal{U}})$ is equivalent to the risk requirement (\ref{mainreq}).  Clearly, if the weighting coefficients are given, one can determine the
coverage tuning parameter $\ze$ to meet the risk requirement by the following two steps:  First, find the maximum number, $\udl{\ze}$, in the
set $\{ 10 \times 2^{-i}: i \in \bb{N} \}$ such that the risk requirement  is satisfied when the coverage tuning parameter $\ze$ assumes value
$\udl{\ze}$. Second, apply a bisection search method to obtain a number $\ze^\star$ as large as possible from interval $[\udl{\ze}, 2
\udl{\ze})$  such that the risk requirement is satisfied when the coverage tuning parameter $\ze$ assumes value $\ze^\star$.  However, these two
steps are not sufficient to produce testing plans of satisfactory efficiency if the weighting coefficients are not properly chosen.   To
overcome this limitation, we observe that to make a testing plan efficient, it is an effective approach to make the testing plan efficient for
parametric values corresponding to the endpoints of the indifference zone. Hence, for determining the appropriate values of the weighting
coefficients, we formulate a minimax optimization problem as follows.  For the family of testing plans parameterized by $a_i, \; b_i, \; i = 1,
\cd, m - 1$ and $\ze$, define {\small $Q = \max_{i \in \{ 1, \cd, m - 1 \} } \max \li \{ A_i, \; B_i  \ri \}$} and {\small $R = \min_{i \in \{
1, \cd, m - 1 \} } \min \li \{ A_i, \; B_i \ri \}$}, where $A_i = \f{ \de_{i-1} } { \Pr \{ \tx{Reject} \; \mscr{H}_{i-1}  \mid \se_i^\prime  \}
}$ and $B_i = \f{ \de_{i} } { \Pr \{ \tx{Reject} \; \mscr{H}_{i} \mid \se_i^{\prime \prime} \} }$ for $i = 1, \cd, m-1$.  It can be seen that if
$R \geq 1$, then the risk requirement is satisfied for parametric values at the endpoints of the indifference zone.  Moreover, under the
restriction that $R \geq 1$, if $Q$ is smaller, then the associated testing plan is more efficient for parametric values at the endpoints of the
indifference zone. Therefore, we propose the following minimization problem: {\it Determine coverage tuning parameter $\ze$ and weighting
coefficients $a_i, \; b_i, \; i = 1, \cd, m-1$ such that $Q$ is minimized under the constraint that $R$ is no less than $1$}. To accomplish such
a task of minimax optimization, we propose the following iterative algorithm.

\vspace{0.05in}

\begin{tabular} {|l |}
\hline $ \nabla  \;  \tx{ Set maximum number of iterations as $k_{\mrm max}$.    Choose initial values of weighting}$\\
$\qu \; \; \tx{coefficients as $a_i
=\de_{i-1}, \; b_i = \de_i$ for $i = 1, \cd, m - 1$. Let $\wh{Q} \leftarrow \iy$ and $k \leftarrow 0$}$.\\
$ \nabla  \;  \tx{While $k \leq k_{\mrm max}$, do the following}$:\\
$ \indent  \; \diamond  \; \tx{Use a bisection search method to determine a number $\ze^*
> 0$ as large as possible}$\\
$ \indent \qu \; \; \tx{for $\ze$ such that the value of $R$ associated with $\ze^*$ and $a_i, \; b_i, \; i = 1, \cd, m-1$}$\\
$\indent \qu \; \;  \tx{is no less than $1$.  Let $Q^*$ and $A_i^*, \; B_i^* , \; i = 1, \cd, m-1$
respectively denote the}$\\
$ \indent \qu \; \; \tx{corresponding values of $Q$ and $A_i, \; B_i, \; i = 1, \cd, m-1$, which are associated}$\\
$ \indent \qu \; \; \tx{with $\ze^*$ and $a_i, \; b_i, \; i = 1, \cd,
m-1$}$.\\
$ \indent \; \diamond  \;  \tx{If $Q^* < \wh{Q}$, then let $\wh{Q} \leftarrow Q^*$ and $\wh{a}_i \leftarrow \ze^* a_i, \; \; \wh{b}_i \leftarrow
\ze^* b_i$ for $i = 1, \cd, m - 1$.}$\\
$\indent \qu \; \; \tx{If there exists an index $j \in \{1, \cd, m-1\}$ such that $A_j^* = Q^*$, then let}$\\
$ \indent \qu \; \; \tx{{\small $a_j \leftarrow \ze^* a_j ( 1 + \f{ Q^* - 1 }{5} )$}. If there exists an index $j \in \{1, \cd, m-1\}$ such that}$\\
$\indent \qu \; \; \tx{$B_j^* = Q^*$, then let {\small $b_j
\leftarrow \ze^* b_j ( 1 + \f{ Q^* - 1 }{5} )$}.  Let $k \leftarrow k + 1$}$.\\
$ \nabla  \;  \tx{Return $\wh{a}_i, \wh{b}_i, \; i = 1, \cd, m - 1$ as the desired weighting coefficients}$.
\\ \hline
\end{tabular}

\vspace{0.05in}

Clearly, the above algorithm returns weighting coefficients such that $Q$ is approximately minimized with $\ze = 1$ and $a_i = \wh{a}_i, \; b_i
= \wh{b}_i, \; i = 1, \cd, m - 1$ subject to the constraint that $R \geq 1$.  With weighting coefficients obtained from the above minimax
optimization procedure, we can use the two steps mentioned at the beginning of this section to obtain $\ze$  as large as possible such that the
risk requirement (\ref{mainreq}) is guaranteed.  Our computational experiences indicate that in many situations, the resultant value of $\ze$ is
equal or very close to $1$.  The intuition behind the above algorithm is that for fixed $\ze$ and $i = 1, \cd, m-1$, $\Pr \{ \tx{Reject} \;
\mscr{H}_{i-1} \mid \se_i^\prime \}$ and $\Pr \{ \tx{Reject} \; \mscr{H}_{i}  \mid \se_i^{\prime \prime} \}$ are ``roughly'' increasing with
respect to $a_i$ and $b_i$, respectively, which can be explained by the following heuristic arguments.

From the parametrization of the stopping and decision rules and their connection with the multi-hypothesis problem defined by (\ref{mainpr}) and
(\ref{mainreq}), it can be seen that $\Pr \{ \tx{Reject} \; \mscr{H}_{i-1} \mid \se_i^\prime \} = \Pr \{  \se_i^\prime \notin  (\bs{\mcal{L}},
\bs{\mcal{U}}) \mid \se_i^\prime  \} = \Pr \{ \bs{\mcal{L}} \geq \se_i^\prime  \mid \se_i^\prime  \} + \Pr \{ \bs{\mcal{U}} \leq \se_i^\prime
\mid \se_i^\prime  \}$.  It follows from statement (II) of Theorem \ref{Multi-Valued Inclusion Principle} that $\Pr \{ \bs{\mcal{U}} \leq
\se_i^\prime \mid \se_i^\prime \} = \Pr \{ \bs{\mcal{U}} \leq \se_{i-1}^{\prime \prime} \mid \se_i^\prime  \} \leq \sum_{j < i} \Pr \{ U_{\ell,
j} \leq \se_{j}^{\prime \prime} \; \tx{for some} \; \ell \mid \se_i^\prime \} \ap \Pr \{ U_{\ell, i-1} \leq  \se_{i-1}^{\prime \prime} \;
\tx{for some} \; \ell \mid \se_i^\prime \}$ and that $\Pr \{ \bs{\mcal{L}} \geq \se_i^\prime \mid \se_i^\prime  \} \leq \sum_{j \geq i} \Pr \{
L_{\ell, j} \geq \se_{j}^{\prime} \; \tx{for some} \; \ell \mid \se_i^\prime \} \ap \Pr \{ L_{\ell, i} \geq  \se_{i}^{\prime} \; \tx{for some}
\; \ell \mid \se_i^\prime \}$.  If the gap between $\se_i^\prime$ and $\se_{i-1}^{\prime \prime}$ is sufficiently large, then $\Pr \{ U_{\ell,
i-1} \leq \se_{i-1}^{\prime \prime} \; \tx{for some} \; \ell \mid \se_i^\prime \}$ will be much smaller than $\Pr \{ U_{\ell, i-1} \leq
\se_{i-1}^{\prime \prime} \; \tx{for some} \; \ell \mid \se_{i-1}^{\prime \prime} \}$ and consequently, it is reasonable to believe that $\Pr \{
L_{\ell, i} \geq  \se_{i}^{\prime} \; \tx{for some} \; \ell \mid \se_i^\prime \}$ is much greater than $\Pr \{ U_{\ell, i-1} \leq
\se_{i-1}^{\prime \prime} \; \tx{for some} \; \ell \mid \se_i^\prime \}$.  This implies that $\Pr \{ \tx{Reject} \; \mscr{H}_{i-1} \mid
\se_i^\prime \}$ is ``dominated'' by $\Pr \{ L_{\ell, i} \geq \se_{i}^{\prime} \; \tx{for some} \; \ell \mid \se_i^\prime \}$, which can be
increased by increasing $a_i$. By a similar argument, it can be seen that $\Pr \{ \tx{Reject} \; \mscr{H}_{i}  \mid \se_i^{\prime \prime} \}$ is
``dominated'' by $\Pr \{ U_{\ell, i} \leq \se_i^{\prime \prime} \; \tx{for some} \; \ell \mid \se_i^{\prime \prime} \}$, which can be increased
by increasing $b_i$.

\section{Stopping and Decision Rules in Terms of Point Estimators}

In this section, we shall apply the general approach presented in Section \ref{secSRI} to the special case of a single discrete-time stochastic
process $(X_n)_{n \in \bb{N}}$ such that $X_1, X_2, \cd$ are identical samples of $X$ which is parameterized by $\se \in \Se$.   Our objective
is to design multistage procedures for the multi-hypotheses testing problem defined by (\ref{mainpr}) and (\ref{mainreq}).  We shall apply the
inclusion principle to construct stopping and decision rules which can be expressed in terms of point estimators.

In general, a testing plan in our proposed framework consists of $s$
stages. For $\ell = 1, \cd, s$, the number of available samples
(i.e., sample size) of the $\ell$-th stage is denoted by $n_\ell$.
For the $\ell$-th stage, a decision variable $\bs{D}_\ell =
\mscr{D}_\ell (X_1, \cd, X_{n_\ell})$ is defined in terms of samples
$X_1, \cd, X_{n_\ell}$ such that $\bs{D}_\ell$ assumes $m+1$
possible values $0, 1, \cd, m$  with the following notion:

(i) Sampling is continued until $\bs{D}_\ell \neq 0$ for some $\ell
\in \{1, \cd, s\}$.

(ii) The hypothesis $\mscr{H}_j$ is accepted at the $\ell$-th stage
if $\bs{D}_\ell = j + 1$ and $\bs{D}_i = 0$ for $1 \leq i < \ell$.

For practical considerations, we shall only focus on sampling
schemes which are closed in the sense that $\Pr \{ \bs{D}_s = 0 \} =
0$.  For efficiency, a sampling scheme should satisfy the condition
that both  $\Pr \{ \bs{D}_1 \neq 0 \}$ and $\Pr \{ \bs{D}_{s-1} = 0
\}$ are greater than zero.

Let $\bs{l}$ denote the index of stage when the sampling is
terminated.  Then, the sample number when the sampling is
terminated, denoted by $\mbf{n}$,  is  equal to $n_{\bs{l}}$. For
the $\ell$-th stage, an estimator $\wh{\bs{\se}}_\ell$ for $\se$ can
be defined based on samples $X_1, \cd, X_{n_\ell}$. Consequently,
the overall estimator for $\se$, denoted by $\wh{\bs{\se}}$, is
equal to $\wh{\bs{\se}}_{\bs{l}}$.  In many cases, decision
variables $\bs{D}_\ell$ can be defined in terms of
$\wh{\bs{\se}}_\ell$. Specially, if $\wh{\bs{\se}}_\ell$ is a
Unimodal-Likelihood Estimator (ULE) of $\se$ for $\ell = 1, \cd, s$,
the design and analysis of multistage sampling schemes can be
significantly simplified.  For a random tuple $X_1, \cd,
X_{\mbf{r}}$ (of deterministic or random length $\mbf{r}$)
parameterized by $\se$, we say that the estimator $\varphi (X_1,
\cd, X_{\mbf{r}})$ is a ULE of $\se$ if $\varphi$ is a multivariate
function such that, for any observation $(x_1, \cd, x_{r})$ of
$(X_1, \cd, X_{\mbf{r}})$, the likelihood function is non-decreasing
with respect to $\se$ no greater than $\varphi (x_1, \cd, x_{r})$
and is non-increasing with respect to $\se$ no less than $\varphi
(x_1, \cd, x_{r})$. For discrete random variables $X_1, \cd, X_{r}$,
the associated likelihood function is the joint probability mass
function $\Pr \{ X_i = x_i, \; i = 1, \cd, r \mid \se \}$. For
continuous random variables $X_1, \cd, X_{r}$, the corresponding
likelihood function is, $f_{X_1, \cd, X_r} (x_1, \cd, x_r, \se)$,
the joint probability density function of random variable $X_1, \cd,
X_{r}$. It should be noted that a ULE may not be a
maximum-likelihood estimator (MLE). On the other side, a MLE may not
be a ULE.

In the sequel, we shall focus on multistage sampling schemes which
can be defined in terms of estimator $\bs{\varphi}_n = \varphi(X_1,
\cd, X_n)$ such that $\bs{\varphi}_n$ is a ULE of $\se$ for every
$n$ and that $\bs{\varphi}_n$ converges in probability to $\se$ in
the sense that, for any $\vep > 0$ and $\de \in (0, 1)$, $\Pr \{ |
\bs{\varphi}_n - \se | \geq \vep \} < \de$ provided that $n$ is
sufficiently large. Such estimator $\bs{\varphi}_n$ is referred to
as a {\it Unimodal-likelihood Consistent Estimator} (ULCE) of $\se$.
For the $\ell$-th stage, the estimator $\wh{\bs{\se}}_\ell$ is
defined as $\bs{\varphi}_{n_\ell} = \varphi (X_1, \cd, X_{n_\ell})$.
Accordingly, the decision variables $\bs{D}_\ell$ can be defined in
terms of estimator $\wh{\bs{\se}}_\ell = \bs{\varphi}_{n_\ell}$.

\subsection{Stopping and Decision Rules from Inclusion Principle}

In this section, we shall propose our general stopping and decision rules
 derived from the inclusion principle proposed in \cite{INCP} and represented in Section \ref{secSRI} of this paper.

To avoid prohibitive burden of computational complexity in the
design process, our global strategy is to construct multistage
sampling schemes of certain structure such that the risks of
erroneously accepting or rejecting a hypothesis can be adjusted by
some parameter $\ze > 0$. Such a parameter $\ze$ is referred to as a
{\it risk tuning parameter} in this paper to convey the idea that
$\ze$ is used to ``tune'' the risk of making a wrong decision to be
acceptable.  As will be seen in the sequel, by virtue of the concept
of ULE, we are able to construct a class of multistage testing plans
such that the risks can be ``tuned'' to be no greater than
prescribed levels by making the risk tuning parameter $\ze$
sufficiently small.  Moreover, the risk tuning can be accomplished
by a bisection search method.  Furthermore, the OC functions of
these multistage testing plans possess some monotonicity which makes
it possible to control the probabilities of committing decision
errors by checking the endpoints of indifference zone.

For the ease of presentation of our sampling schemes, we need to
introduce some multivariate functions regarding estimator
$\bs{\varphi}_n = \varphi( X_1, \cd, X_n)$ of $\se$.  For $n \in
\bb{N}, \se \in \Se, \de \in (0, 1)$, define {\small \bee & & f (n,
\se, \de) = \bec \max \{ z \in I_{ \bs{\varphi}_n } : F_{
\bs{\varphi}_n } (z, \se) \leq \de, \; z \leq \se \} & \tx{if} \; \{
F_{ \bs{\varphi}_n } (\bs{\varphi}_n, \se)
\leq \de, \; \bs{\varphi}_n \leq \se \} \neq \emptyset,\\
- \iy & \tx{otherwise}  \eec\\
&  & g (n, \se, \de) = \bec \min \{ z \in I_{ \bs{\varphi}_n } : G_{
\bs{\varphi}_n } (z, \se) \leq \de, \; z \geq \se \} & \tx{if} \; \{
G_{ \bs{\varphi}_n } (\bs{\varphi}_n, \se) \leq
\de, \; \bs{\varphi}_n \geq \se \} \neq \emptyset,\\
\iy & \tx{otherwise}  \eec
 \eee}
For $\se^\prime < \se^{\prime \prime}$ contained in ${\Se}$ and
$\de^\prime, \de^{\prime \prime} \in (0, 1)$, define {\small \bee &
& \udl{f} (n, \se^\prime, \se^{\prime \prime}, \de^\prime,
\de^{\prime \prime} ) = \min \li \{ f (n, \se^{\prime \prime},
\de^{\prime \prime}), \qu \f{1}{2} [ f (n, \se^{\prime \prime},
\de^{\prime \prime}) + g (n, \se^{\prime}, \de^{\prime}) ]  \ri \}, \\
 &  & \ovl{g} (n, \se^\prime, \se^{\prime \prime},
\de^\prime, \de^{\prime \prime} ) = \max \li \{ g (n, \se^{\prime}, \de^{\prime}), \qu \f{1}{2} [ f (n, \se^{\prime \prime}, \de^{\prime
\prime}) + g (n, \se^{\prime}, \de^{\prime}) ] \ri \}. \eee} By virtue of the inclusion principle, we have derived a general method for
constructing multistage test plans in terms point estimators and their properties described by Theorem \ref{Multi_Comp_Exact} as follows.

 \beT
\la{Multi_Comp_Exact}

Let $\al_i = O(\ze) \in (0, 1), \ba_i = O(\ze) \in (0, 1)$ for $i =
1, \cd, m - 1$ and $\al_m = \ba_0 = 0$. Define $\ovl{\al}_i = \max
\{ \al_j: i < j \leq m \}$ and $\ovl{\ba}_i = \max \{ \ba_j: 0 \leq
j \leq i \}$ for $i = 0, 1, \cd, m - 1$. Suppose that
$\bs{\varphi}_n$ is a ULCE of $\se$. Suppose that the maximum sample
size $n_s$ is no less than the minimum integer $n$ such that $f (n,
\se_i^{\prime \prime}, \ba_i) \geq g (n, \se_i^{\prime}, \al_i)$ for
$i = 1, \cd, m-1$. Define $f_{\ell, i} = \udl{f} (n_\ell,
\se_i^{\prime}, \se_i^{\prime \prime}, \al_i, \ba_i)$ and $g_{\ell,
i} = \ovl{g} (n_\ell, \se_i^{\prime}, \se_i^{\prime \prime}, \al_i,
\ba_i)$ for $i = 1, \cd, m-1$. Define \be \la{defgoog}
 \bs{D}_\ell = \bec 1 & \tx{if} \; \;
\wh{\bs{\se}}_\ell \leq
f_{\ell ,1}, \\
i & \tx{if} \; \;  g_{\ell, i-1} < \wh{\bs{\se}}_\ell \leq
f_{\ell,i}
\; \tx{where} \; 2 \leq i \leq m - 1,\\
m  & \tx{if} \; \;  \wh{\bs{\se}}_\ell > g_{\ell ,m-1},\\
0 & \tx{else} \eec \ee for $\ell = 1, \cd, s$.  The following
statements (I)-(VI) hold true for $m \geq 2$.

(I) $\Pr \{ \tx{Reject} \; \mscr{H}_0 \mid \se \}$ is non-decreasing
with respect to $\se \in \varTheta_0$.

(II) $\Pr \{ \tx{Reject} \; \mscr{H}_{m-1} \mid \se \} \; \tx{is
non-increasing with respect to} \; \se \in \varTheta_{m-1}$.

(III)  $\Pr \{ \tx{Reject} \; \mscr{H}_i  \mid \se \} \leq s (
\ovl{\al}_i  + \ovl{\ba}_i  )$ for any $\se \in \varTheta_i$ and $i
= 0, 1, \cd, m - 1$.

(IV) For $0 < i \leq m- 1$, $\Pr \{ \tx{Accept} \; \mscr{H}_i \mid
\se \}$ is no greater than $s \al_i$ and is non-decreasing with
respect to $\se \in {\Se}$ no greater than $\se_i^\prime$.

(V) For $0 \leq i \leq m-2$, $\Pr \{ \tx{Accept} \; \mscr{H}_i \mid
\se \}$ is no greater than $s \ba_{i + 1}$ and is non-increasing
with respect to $\se \in {\Se}$ no less than $\se_{i+1}^{\prime
\prime}$.

(VI)  Assume that $\bb{E} [e^{\ro X}]$ exists for any $\ro \in
\bb{R}$ and  that {\small $\bs{\varphi}_n = \f{ \sum_{i=1}^n X_i
}{n}$} is an unbiased and unimodal-likelihood estimator of $\se$,
where $X_1, X_2,  \cd$ are i.i.d. samples of $X$.   Then, for $i =
0, 1, \cd, m - 1$, $\lim_{\ze \to 0} \Pr \{ \tx{Reject} \;
\mscr{H}_i \mid \se \} = 0$ for any $\se \in \varTheta_i$.

 Moreover, the following statements
(VII), (VIII) and (IX) hold true for $m \geq 3$.

(VII) \bee & & \Pr \{ \tx{Reject} \; \mscr{H}_i \mid \se \} \leq \Pr
\{ \tx{Reject} \; \mscr{H}_i, \; \wh{\bs{\se}} \leq a \mid a \} +
\Pr \{ \tx{Reject} \; \mscr{H}_i, \; \wh{\bs{\se}} \geq b \mid b \},\\
&  & \Pr \{ \tx{Reject} \; \mscr{H}_i \mid \se \} \geq  \Pr \{
\tx{Reject} \; \mscr{H}_i, \; \wh{\bs{\se}} \leq a \mid b \} + \Pr
\{ \tx{Reject} \; \mscr{H}_i, \; \wh{\bs{\se}} \geq b \mid a \} \eee
for any $\se \in [a, b ] \subseteq \varTheta_i$ and $1 \leq i \leq m
- 2$.

(VIII) $\Pr \{ \tx{Reject} \; \mscr{H}_0 \; \tx{and} \;
\mscr{H}_{m-1} \mid \se \}$ is non-decreasing with respect to $\se
\in \varTheta_0$ and is non-increasing with respect to $\se \in
\varTheta_{m-1}$.

(IX) $\Pr \{ \tx{Reject} \; \mscr{H}_0 \; \tx{and} \; \mscr{H}_{m-1}
\mid \se \}$ is no greater than $s \times \max \{ \al_i: 1 \leq i
\leq m - 2 \}$ for $\se \in \varTheta_0$ and is no greater than $s
\times \max \{ \ba_i: 2 \leq i \leq m - 1 \}$ for $\se \in
\varTheta_{m-1}$.

\eeT

See Appendix \ref{Multi_Comp_Exact_Ap} for a proof.

In situations that the parameter $\se$ to be tested is the expectation of $X$, we can apply normal approximation to simplify the stopping and
decision rules.  Assume that $X_1, X_2, \cd$ are identical samples of $X$ and that the variance of $\bs{\varphi}_n = \f{\sum_{i=1}^n X_i}{n}$ is
a bivariate function, denoted by $\mscr{V} (\se, n)$, of $\se$ and $n$. If all sample sizes are large, then the central limit theorem may be
applied to establish the normal approximation \bee &  & F_{\bs{\varphi}_n} \li ( z, \se \ri ) \DEF \Pr \{  \bs{\varphi}_n \leq z \mid \se  \}
\ap \Phi \li ( \f{ z - \se } { \sqrt{ \mscr{V} ( \se, \; n ) } } \ri ), \\
&  &  G_{\bs{\varphi}_n} \li ( z, \se \ri ) \DEF \Pr \{  \bs{\varphi}_n \geq z \mid \se  \} \ap \Phi \li ( \f{ \se - z } { \sqrt{ \mscr{V} (
\se, \; n ) } } \ri ) \eee and consequently, the stopping and decision rule described by Theorem \ref{Multi_Comp_Exact} can be simplified by
applying the approximation to redefine $f(n, \se, \de)$ and $g(n, \se, \de)$ as follows: {\small \bee & & f (n, \se, \de) = \bec \max \li \{ z
\in I_{ \bs{\varphi}_n } : \Phi \li ( \f{ z - \se } { \sqrt{ \mscr{V} ( \se, \; n ) } } \ri ) \leq \de, \; z \leq \se \ri \} & \tx{if} \; \li \{
\Phi \li ( \f{ \bs{\varphi}_n - \se }
{ \sqrt{ \mscr{V} ( \se, \; n ) } } \ri ) \leq \de, \; \bs{\varphi}_n \leq \se \ri \} \neq \emptyset,\\
- \iy & \tx{otherwise}  \eec\\
&  & g (n, \se, \de) = \bec \min \li \{ z \in I_{ \bs{\varphi}_n } : \Phi \li ( \f{ \se - z } { \sqrt{ \mscr{V} ( \se, \; n ) } } \ri ) \leq
\de, \; z \geq \se \ri \} & \tx{if} \; \li \{ \Phi \li ( \f{ \se - \bs{\varphi}_n } { \sqrt{ \mscr{V} ( \se, \; n ) } } \ri ) \leq
\de, \; \bs{\varphi}_n \geq \se \ri \} \neq \emptyset,\\
\iy & \tx{otherwise}  \eec
 \eee}
for $n \in \bb{N}, \se \in \Se, \de \in (0, 1)$. Except this modification, the definition of the stopping and decision rules remain unchanged.
It should be noted that this is not the best approximation method for simplifying the stopping and decision rules.  Our computational
experiences indicate that the accuracy of normal approximation can be improved by replacing $\se$ in $\mscr{V} ( \se, \; n )$ as $z + w (\se -
z)$, where $w \in [0, 1]$. In other words, we propose a new normal approximation as follows: \bee &  & F_{\bs{\varphi}_n} \li ( z, \se \ri )
\ap \Phi \li ( \f{ z - \se } { \sqrt{ \mscr{V} ( z + w (\se - z), \; n ) } } \ri ), \\
&  &  G_{\bs{\varphi}_n} \li ( z, \se \ri ) \ap \Phi \li ( \f{ \se - z } { \sqrt{ \mscr{V} ( z + w (\se - z), \; n ) } } \ri ). \eee
Accordingly, the stopping and decision rule described by Theorem \ref{Multi_Comp_Exact} can be simplified  by redefining  $f(n, \se, \de)$ and
$g(n, \se, \de)$ as follows: {\small \bee & & f (n, \se, \de) = \bec \max \li \{ z \in I_{ \bs{\varphi}_n } : \Phi \li ( \f{ z - \se } { \sqrt{
\mscr{V} ( z + w (\se - z), \; n ) } } \ri ) \leq \de, \; z
\leq \se \ri \} & \tx{if} \;  \mscr{A} \neq \emptyset,\\
- \iy & \tx{otherwise}  \eec\\
&  & g (n, \se, \de) = \bec \min \li \{ z \in I_{ \bs{\varphi}_n } : \Phi \li ( \f{ \se - z } { \sqrt{ \mscr{V} ( z + w (\se - z), \; n ) } }
\ri ) \leq \de, \; z \geq \se \ri \} & \tx{if} \;  \mscr{B} \neq \emptyset,\\
\iy & \tx{otherwise}  \eec
 \eee}
for $n \in \bb{N}, \se \in \Se, \de \in (0, 1)$, where $\mscr{A} \DEF \li \{ \Phi \li ( \f{ \bs{\varphi}_n - \se } { \sqrt{ \mscr{V} (
\bs{\varphi}_n + w (\se - \bs{\varphi}_n), \; n ) } } \ri ) \leq \de, \; \bs{\varphi}_n \leq \se \ri \}$ and $\mscr{B} \DEF \li \{ \Phi \li (
\f{ \se - \bs{\varphi}_n } { \sqrt{ \mscr{V} ( \bs{\varphi}_n + w (\se - \bs{\varphi}_n), \; n ) } } \ri ) \leq \de, \; \bs{\varphi}_n \geq \se
\ri \}$.  As before, except this modification, the definition of the stopping and decision rules remain unchanged.

Although approximation methods are used, for many problems, the risk requirements can be guaranteed by choosing $\ze$ to be a sufficiently small
number.  Moreover, the performance of the testing plans can be optimized with respect to $w \in [0, 1]$.  Clearly, this approach of constructing
simple stopping and decision rules applies to the problems of testing binomial proportion, Poisson parameter, and finite population proportion.

In addition to the normal approximation, bounds of tail
probabilities of $\bs{\varphi}_n = \f{\sum_{i=1}^n X_i}{n}$, where
$X_1, X_2, \cd$ are identical samples of $X$ as before,  can be used
to simplify stopping and decision rules.  To proceed in this
direction, define multivariate functions {\small \bee &  & f_c (n,
\se, \de) = \bec \max \{ z \in I_{ \bs{\varphi}_n } : [ \mscr{C} (z,
\se) ]^n \leq \de, \; z \leq \se \} & \tx{if} \; \{ [ \mscr{C}
(\bs{\varphi}_n, \se) ]^n \leq \de, \; \bs{\varphi}_n \leq \se \} \neq \emptyset,\\
- \iy & \tx{otherwise}  \eec\\
&  & g_c (n, \se, \de) = \bec \min \{ z \in I_{ \bs{\varphi}_n } : [
\mscr{C} (z, \se) ]^n \leq \de, \; z \geq \se \} & \tx{if} \; \{ [
\mscr{C} (\bs{\varphi}_n, \se) ]^n \leq
\de, \; \bs{\varphi}_n \geq \se \} \neq \emptyset,\\
\iy & \tx{otherwise}  \eec
 \eee}
for $n \in \bb{N}, \se \in \Se, \de \in (0, 1)$, where $\mscr{C} (z,
\se) = \inf_{\ro \in \bb{R} } \bb{E} [  e^{\ro ( X - z)} ]$.
Moreover, define {\small \bee & & \udl{f}_c (n, \se^\prime,
\se^{\prime \prime}, \de^\prime, \de^{\prime \prime} ) = \min\li \{
f_c (n, \se^{\prime \prime}, \de^{\prime \prime}), \qu \f{1}{2} [
f_c (n, \se^{\prime \prime},
\de^{\prime \prime}) + g_c (n, \se^{\prime}, \de^{\prime}) ] \ri \}, \\
&  & \ovl{g}_c (n, \se^\prime, \se^{\prime \prime}, \de^\prime,
\de^{\prime \prime} ) = \max \li \{ g_c (n, \se^{\prime},
\de^{\prime}), \qu \f{1}{2} [ f_c (n, \se^{\prime \prime},
\de^{\prime \prime}) + g_c (n, \se^{\prime}, \de^{\prime}) ] \ri \}
\eee} for $\se^\prime < \se^{\prime \prime}$ in $\Se, \; \de^\prime,
\de^{\prime \prime} \in (0, 1)$ and $n \in \bb{N}$.

Our sampling schemes and their properties can be described by
Theorem \ref{Multi_Comp_Chernoff} as follows.

\beT \la{Multi_Comp_Chernoff}  Let $\al_i = O(\ze) \in (0, 1), \ba_i
= O(\ze) \in (0, 1)$ for $i = 1, \cd, m - 1$ and $\al_m = \ba_0 =
0$. Define $\ovl{\al}_i = \max \{ \al_j: i < j \leq m \}$ and
$\ovl{\ba}_i = \max \{ \ba_j: 0 \leq j \leq i \}$ for $i = 0, 1,
\cd, m - 1$. Suppose that $\bb{E} [e^{\ro X}]$ exists for any $\ro
\in \bb{R}$ and {\small $\bs{\varphi}_n = \f{ \sum_{i=1}^n X_i
}{n}$} is an unbiased and unimodal-likelihood estimator of $\se$.
Suppose that the maximum sample size $n_s$ is no less than the
minimum integer $n$ such that $f_c (n, \se_i^{\prime \prime}, \ba_i)
\geq g_c (n, \se_i^{\prime}, \al_i)$ for $i = 1, \cd, m-1$ . Define
decision variable $\bs{D}_\ell$ by (\ref{defgoog}) for $\ell = 1,
\cd, s$ with $f_{\ell, i} = \udl{f}_c (n_\ell, \se_i^{\prime},
\se_i^{\prime \prime}, \al_i, \ba_i)$ and $g_{\ell, i} = \ovl{g}_c
(n_\ell, \se_i^{\prime}, \se_i^{\prime \prime}, \al_i, \ba_i)$ for
$i = 1, \cd, m-1$.   Then, the same conclusion as described by
statements (I)--(IX) of Theorem \ref{Multi_Comp_Exact} holds true.

\eeT

Theorem \ref{Multi_Comp_Chernoff} can be established by making use of Lemmas \ref{ProbTrans}, \ref{ULE_Basic}, and \ref{Unify_CH} and an
argument similar to the proof of Theorem \ref{Multi_Comp_Exact}.

In the preceding discussion, Chernoff bounds and normal approximations are used to simplifying stopping and decision rules.  In addition to
these techniques, we can also use the bounds on the distribution of the likelihood ratios to construct simple stopping and decision rules.  For
this purpose, we have derived the following results.

\beT \la{boundCDFLR} Let $\al$ be a positive number and $n$ be a positive integer.  Let $f_n( x_1, \cd, x_n; \se)$ denote the joint probability
density (or mass) function of random variables $X_1, \cd, X_n$ parameterized by $\se \in \Se$. Assume that $\f{ f_n( X_1, \cd, X_n; \se_1)
}{f_n( X_1, \cd, X_n; \se_0)}$ can be expressed as a function, $\Lm (\varphi_n , \se_0, \se_1)$,  of $\se_0, \se_1$ and $\varphi_n =
\varphi(X_1, \cd, X_n)$ such that $\Lm (\varphi_n , \se_0, \se_1)$ is increasing with respect to $\varphi_n$. Let $\wh{\se}_n$ be a function of
$\varphi_n$ such that $\wh{\se}_n$ takes values in $\Se$. Then, \bel &  & \Pr \li \{ \f{ f_n( X_1, \cd, X_n; \se) }{f_n( X_1, \cd, X_n;
\wh{\se}_n)} \leq \f{\al}{2}, \; \wh{\se}_n \leq \se \mid \se \ri \}
\leq \f{\al}{2}, \la{RB1}\\
&  & \Pr \li \{ \f{ f_n( X_1, \cd, X_n; \se) }{f_n( X_1, \cd, X_n; \wh{\se}_n)} \leq \f{\al}{2}, \; \wh{\se}_n \geq \se \mid \se \ri \} \leq \f{\al}{2}, \la{RB2}\\
&  & \Pr \li \{ \f{ f_n( X_1, \cd, X_n; \se) }{f_n( X_1, \cd, X_n; \wh{\se}_n)} \leq \f{\al}{2} \mid \se \ri \} \leq \al \la{RB3} \eel for $\se
\in \Se$.  Moreover, under additional assumption that $\wh{\se}_n$ is a ULE for $\se$, the following inequalities \bel &  & \Pr \li \{  \f{
\sup_{\vse \in \mscr{S}} f_n( X_1, \cd, X_n; \vse) }{\sup_{\vse \in \Se} f_n( X_1, \cd, X_n; \vse)} \leq \f{\al}{2}, \;
\wh{\se}_n \leq \inf \mscr{S} \mid \se \ri \} \leq \f{\al}{2} \la{LBA},\\
 &  & \Pr \li \{  \f{ \sup_{\vse \in \mscr{S}} f_n( X_1, \cd, X_n;
\vse) }{\sup_{\vse \in \Se} f_n( X_1, \cd,
X_n; \vse)} \leq \f{\al}{2}, \; \wh{\se}_n \geq \sup \mscr{S} \mid \se \ri \} \leq \f{\al}{2},  \la{LBB}\\
&  & \Pr \li \{  \f{ \sup_{\vse \in \mscr{S}} f_n( X_1, \cd, X_n; \vse) }{\sup_{\vse \in \Se} f_n( X_1, \cd, X_n; \vse)} \leq \f{\al}{2} \mid
\se \ri \} \leq \al \la{LBC}
 \eel hold true for arbitrary nonempty subset $\mscr{S}$ of $\Se$ and all
$\se \in \mscr{S}$.

\eeT

See Appendix \ref{boundCDFLR_app} for a proof.  By virtue of Theorem \ref{boundCDFLR},  we can construct confidence intervals for $\se$, which
can be used to construct stopping and decision rules based on the inclusion principle.

\subsection{Bisection Risk Tuning} \la{BRT}

In this section, we shall propose bisection risk tuning method based on the groundwork established in Theorems \ref{Multi_Comp_Exact} and
\ref{Multi_Comp_Chernoff}.    In applications, the number of stages $s$ and the sample sizes $n_1, \cd, n_s$ can be defined as certain functions
of $\ze$ and $\al_i, \ba_i, \; i = 1, \cd, m - 1$.  It can be seen from Theorems \ref{Multi_Comp_Exact} and \ref{Multi_Comp_Chernoff} that if
$\al_i, \ba_i, \; i = 1, \cd, m - 1$ are given functions of the risk tuning parameter $\ze$, then the corresponding stopping and decision rules
are actually parameterized by the risking tuning parameter $\ze$.   Assuming that $\al_i,  \ba_i, \; i = 1, \cd, m - 1$ are given functions of
the risk tuning parameter $\ze$, the objective of bisection risk tuning is to determine $\ze$ as large as possible such that the risk
requirement (\ref{mainreq}) is satisfied. The procedure of bisection risk tuning is illustrated as follows.

According to statement (VI) of Theorem \ref{Multi_Comp_Exact}, $\Pr \{ \tx{Reject} \; \mscr{H}_i \mid \se \}$ tends to $0$ as $\ze$ tends to
$0$.  This implies that we can ensure (\ref{mainreq}) by choosing a sufficiently small risk tuning parameter $\ze$. Clearly, every value of
$\ze$ determines a test plan and consequently its performance specifications such as average sample number (ASN) and risks of making wrong
decisions. Intuitively, under the constraint of risk requirements, the risk tuning parameter $\ze$ should be chosen as large as possible in
order to reduce the sample number.  To achieve such an objective, it is a critical subroutine to determine whether a given $\ze$ is sufficient
to ensure the risk requirement (\ref{mainreq}).  Since there may be an extremely large number or infinite parametric values in $\cup_{i=0}^{m-1}
\varTheta_i$, it is essential to develop an efficient method to check the risk requirement (\ref{mainreq}) without exhaustive computation.  For
this purpose, statements (I), (II) and (VI) of Theorem \ref{Multi_Comp_Exact} can be very useful.
 As a consequence of statement (I),
 to check if $\Pr \{ \tx{Reject} \; \mscr{H}_0 \mid \se \} \leq
\de_0$ for any $\se \in \varTheta_0$, it suffices to check whether
$\Pr \{ \tx{Reject} \; \mscr{H}_0 \mid \se_1^\prime \} \leq \de_0$
is true.  By virtue of statement (II), for purpose of determining
whether $\Pr \{ \tx{Reject} \; \mscr{H}_{m-1} \mid \se \} \leq
\de_{m-1}$ for any $\se \in \varTheta_{m-1}$, it is sufficient to
check if $\Pr \{ \tx{Reject} \; \mscr{H}_{m-1} \mid
\se_{m-1}^{\prime \prime} \} \leq \de_{m-1}$ is true.   For $i \in
\{1, \cd, m-2 \}$, to determine whether $\Pr \{ \tx{Reject} \;
\mscr{H}_i \mid \se \} \leq \de_i$ for any $\se \in \varTheta_i$, we
can apply the bounding results in statement (VI) of Theorem
\ref{Multi_Comp_Exact} and the Adaptive Maximum Checking Algorithm
(AMCA) established in \cite{Chen_EST}.  Therefore, it is clear that
we can develop an efficient subroutine to determine whether a given
$\ze$ guarantees the risk requirement (\ref{mainreq}).  Now, let
$\udl{\ze}$ be the maximum number in the set $\{ 10 \times 2^{-i}: i
\in \bb{N} \}$ such that the risk requirement (\ref{mainreq}) is
satisfied when the risk tuning parameter $\ze$ assumes value
$\udl{\ze}$.  Such number $\udl{\ze}$ can be obtained by using the
subroutine to check the risk requirement (\ref{mainreq}). Once
$\udl{\ze}$ is found, we can apply a bisection search to obtain a
number $\ze^\star$ as large as possible from interval $[\udl{\ze}, 2
\udl{\ze})$  such that the risk requirement (\ref{mainreq}) is
satisfied when the risk tuning parameter $\ze$ assumes value
$\ze^\star$.

The above bisection risk tuning technique can be straightforwardly
extended to control the following error probabilities: \bee &  & \Pr
\{ \tx{Accept} \; \mscr{H}_i \mid \se \in \varTheta_j \}, \qqu 0 \leq i < j \leq m - 1\\
&  & \Pr \{  \tx{Accept} \; \mscr{H}_i \mid \se \in \varTheta_j \},
\qqu 0 \leq j < i \leq m - 1\\
&  & \Pr \{ \tx{Accept} \; \mscr{H}_i \mid \se \in \Se_j \},
\qqu 0 \leq i \leq j - 2 < j \leq m - 1\\
&  & \Pr \{
 \tx{Accept} \; \mscr{H}_i \mid \se \in \Se_j \}, \qqu 0 \leq j \leq i - 2 < i \leq m - 1
 \eee
For this purpose, statements (IV) and (V) of Theorem
\ref{Multi_Comp_Exact} can be used to develop efficient method of
checking the above risk requirements.  In a similar spirit, by
virtue of statements (VII) and (VIII) of Theorem
\ref{Multi_Comp_Exact}, the control of $\Pr \{ \tx{Reject} \;
\mscr{H}_0 \; \tx{and} \; \mscr{H}_{m - 1} \mid \se \in \varTheta_0
\cup \varTheta_{m-1} \}$ can be incorporated in the bisection risk
tuning technique.  As can be seen from above discussion, a critical
idea in the tuning technique is to avoid exhaustive computation by
making use of monotonicity of error probabilities with respect to
$\se$.

\subsection{Minimax Optimization for Determining Weighting Coefficients}  \la{minmaxcom}

As can be seen from Section \ref{BRT}, to construct testing plans based on the bisection risk tuning technique, it is necessary to choose
appropriate forms for $\al_i, \ba_i, \; i = 1, \cd, m - 1$ as functions of $\ze$. It is indicated in Theorems \ref{Multi_Comp_Exact} and
\ref{Multi_Comp_Chernoff} that these functions should be taken within the class $O (\ze)$.  Specifically, in order to apply the bisection risk
tuning technique, we propose to choose $\al_i, \ba_i, \; i = 1, \cd, m - 1$ as \be \la{weight}
 \al_i = \ze a_i, \qqu \ba_i = \ze b_i, \qqu i = 1, \cd, m-1, \ee where the constants $a_i, b_i, \; i = 1, \cd, m-1$ are referred to as {\it weighting coefficients} in
this paper.  The notion of this terminology is derived from the following intuition:

(i) For $i = 1, \cd, m-1$, $\Pr \{ \tx{Reject} \; \mscr{H}_{i-1}  \mid \se_i^\prime  \}$ is ``roughly'' increasing with respect to $a_i$;

(ii) For $i = 1, \cd, m-1$,  $\Pr \{ \tx{Reject} \; \mscr{H}_{i}  \mid \se_i^{\prime \prime} \}$ is ``roughly'' increasing with respect to
$b_i$.

The probabilities of making wrong decisions are affected by the functions $\al_i, \ba_i, \; i = 1, \cd, m - 1$ through the risk tuning parameter
$\ze$ and the weighting coefficients.

Using  $\al_i, \; \ba_i$ defined by (\ref{weight}), we can apply Theorems \ref{Multi_Comp_Exact} and \ref{Multi_Comp_Chernoff} to define testing
plans.  Clearly, the weighting coefficients significantly impact the efficiency of the resultant testing plans.  We observe that to make a
testing plan efficient, it is an effective approach to make the testing plan efficient for parametric values corresponding to the endpoints of
the indifference zones. Hence, for determining the appropriate values of the weighting coefficients, we formulate a minimax optimization problem
as follows.

Consider a family of testing plans associated with $\al_i, \; \ba_i$ defined by (\ref{weight}).   Define
\[
A_i = \f{ \de_{i-1} } { \Pr \{ \tx{Reject} \; \mscr{H}_{i-1}  \mid \se_i^\prime  \} }, \qqu B_i = \f{ \de_{i} } { \Pr \{ \tx{Reject} \;
\mscr{H}_{i}  \mid \se_i^{\prime \prime} \} } \qu \tx{for $i = 1, \cd, m-1$}
\]
and
\[
Q = \max_{i \in \{ 1, \cd, m - 1 \} } \max \li \{ A_i, \; B_i  \ri \}, \qqu R = \min_{i \in \{ 1, \cd, m - 1 \} } \min \li \{ A_i, \; B_i \ri
\}.
\]
Clearly, $Q, \; R$ and $A_i, B_i, \; i = 1, \cd, m-1$ are functions of $\ze$ and $a_i, \; b_i, \; i = 1, \cd, m-1$.  It can be seen that if $R
\geq 1$, then the risk requirement is satisfied for parametric values at the endpoints of the indifference zones.  Moreover, under the
restriction that $R \geq 1$, if $Q$ is smaller, then the associated testing plan is more efficient for parametric values at the endpoints of the
indifference zones.  Therefore, we propose the following minimization problem:

{\it Determine risk tuning parameter $\ze$ and weighting coefficients $a_i, \; b_i, \; i = 1, \cd, m-1$ such that $Q$ is minimized under the
constraint that $R$ is no less than $1$}.

Actually, this is a minimax optimization problem, since the quantity to be minimized is also a maximum over a discrete set. To accomplish such a
task of minimax optimization, we propose the following iterative algorithm based on the intuition as stated in the above statements (i) and
(ii).

\bed

\item [\underline{Step 1}]: Set the maximum number of iterations as
$k_{\mrm max}$.    Choose the initial values of weighting coefficients as $a_i = \de_{i-1}, \; b_i = \de_i$ for $i = 1, \cd, m - 1$. Let $\wh{Q}
\leftarrow \iy$ and $k \leftarrow 0$.

\item  [\underline{Step 2}]:  While $k \leq k_{\mrm max}$, do the following:

\bed

\item [\underline{Step 2-1}]: Based on $\al_i = \ze a_i, \;  \ba_i = \ze b_i$ for $i = 1, \cd, m-1$,  use a bisection search method to determine $\ze > 0$ as large as possible such that the value of $R$
associated with $\ze$ and $a_i, \; b_i, \; i = 1, \cd, m-1$ is no less than $1$. Let the value of $\ze$ obtained at this step be denoted by
$\ze^*$.  Let $Q^*, \; R^*$ and $A_i^*, \; B_i^* , \; i = 1, \cd, m-1$ respectively denote the corresponding values of $Q, \; R$ and $A_i, \;
B_i, \; i = 1, \cd, m-1$, which are associated with $\ze^*$ and $a_i, \; b_i, \; i = 1, \cd, m-1$.

\item [\underline{Step 2-2}]: If $Q^* < \wh{Q}$, then let $\wh{Q} \leftarrow Q^*$
and $\wh{a}_i \leftarrow \ze^* a_i, \; \; \wh{b}_i  \leftarrow \ze^* b_i$ for $i = 1, \cd, m - 1$.

\item [\underline{Step 2-3}]: If there exists an index $j \in \{1, \cd, m-1\}$ such that
$A_j^* = Q^*$, then let $a_j \leftarrow \ze^* a_j \li ( 1 + \f{ Q^* - 1 }{5} \ri )$.

\item [\underline{Step 2-4}]: If there exists an index $j \in \{1, \cd, m-1\}$ such that $B_j^* = Q^*$,
then let $b_j \leftarrow \ze^* b_j \li ( 1 + \f{ Q^* - 1 }{5} \ri )$.

\item [\underline{Step 2-5}]:  Let $k \leftarrow k + 1$.

\eed

\item [\underline{Step 3}]: Return $\wh{a}_i, \wh{b}_i, \; i = 1, \cd, m - 1$ as the desired weighting coefficients.

 \eed

We would like to point out that the number ``$5$'' in the denominator of $\f{ Q^* - 1 }{5}$ appeared in Steps 2-3 and 2-4 may not be the optimal
value for the efficiency of minimization. Our computational experiences indicate that the numbers chosen from interval $[4, 8]$ work reasonably
well.

 Clearly, the above algorithm returns weighting coefficients such that $Q$ is approximately minimized with $\ze = 1$ and
$a_i = \wh{a}_i, \; b_i = \wh{b}_i, \; i = 1, \cd, m - 1$ subject to the constraint that $R \geq 1$.  With weighting coefficients obtained from
the above minimax optimization procedure, we can use the bisection risk tuning technique proposed in Section \ref{BRT} to obtain $\ze$  as large
as possible such that the risk requirement (\ref{mainreq}) is guaranteed.  Our computational experiences indicate that in many situations, the
resultant value of $\ze$ is equal or very close to $1$.  This is consistent with our observation that for many problem cases, the maximum of
probabilities of incorrectly rejecting hypotheses is attained at some endpoint of the indifference zones.

\subsection{Recursive Computation}

As will be seen in the sequel, for most multistage sampling schemes for testing parameters of discrete variables, the computation of the OC
functions involve probabilistic terms like $\Pr \{ K_i \in \mscr{K}_i, \; i = 1, \cd, \ell  \}, \; \ell = 1, 2, \cd$,  where $K_\ell = \sum_{i =
1}^{n_\ell} X_i$ and $\mscr{K}_i$ is a subset of integers. The calculation of such terms can be performed by virtue of the following recursive
relationship: \bel & & \Pr \{ K_{\ell + 1} = k_{\ell + 1}; \; K_i \in \mscr{K}_i, \; i
= 1, \cd, \ell \} \nonumber\\
&   & = \sum_{k_\ell \in \mscr{K}_\ell} \li [ \Pr \{ K_\ell = k_\ell; \; K_i \in \mscr{K}_i, \; i = 1, \cd, \ell - 1 \} \ri. \\
&   & \qqu \qu \li. \times \Pr \{ K_{\ell + 1} - K_\ell = k_{\ell + 1} - k_\ell \mid K_\ell = k_\ell; \;  K_i \in \mscr{K}_i, \; i = 1, \cd,
\ell - 1 \} \ri ], \la{recur1} \eel where the computation of the conditional probability $\Pr \{ K_{\ell + 1} - K_\ell = k_{\ell + 1} - k_\ell
\mid K_\ell = k_\ell; \; K_i \in \mscr{K}_i, \; i = 1, \cd, \ell - 1 \}$ depends on specific problems. In the context of testing a binomial
parameter $p$, we have {\small \[ \Pr \{ K_{\ell + 1} - K_\ell = k_{\ell + 1} - k_\ell \mid K_\ell = k_\ell; \;  K_i \in \mscr{K}_i, \; i = 1,
\cd, \ell - 1 \} = \bi{ n_{\ell + 1} - n_\ell }{ k_{\ell + 1} - k_\ell } p^{ k_{\ell + 1} - k_\ell } ( 1 - p )^{ n_{\ell + 1} - n_\ell - k_{\ell
+ 1} + k_\ell }.
\]}
In the context of testing a Poisson parameter $\lm$, we have {\small \[ \Pr \{ K_{\ell + 1} - K_\ell = k_{\ell + 1} - k_\ell \mid K_\ell =
k_\ell; \; K_i \in \mscr{K}_i, \; i = 1, \cd, \ell - 1 \} = \f{ [ (n_{\ell + 1} - n_\ell ) \lm ]^{k_{\ell + 1} - k_\ell} \exp ( - (n_{\ell + 1}
- n_\ell) \lm  ) } { (k_{\ell + 1} - k_\ell)!  }.
\]}
In the context of testing the proportion, $p$, of a finite population of size $N$ using multistage sampling schemes to be described in Section
\ref{Finite_Plan}, we have {\small \be \la{confinite}
 \Pr \{ K_{\ell + 1} - K_\ell = k_{\ell + 1} - k_\ell \mid K_\ell = k_\ell; \;  K_i \in \mscr{K}_i,
\; i = 1, \cd, \ell - 1 \} = \f{ \bi{p N - k_\ell}{k_{\ell + 1} - k_\ell} \bi{N - n_\ell - p N + k_\ell} { n_{\ell + 1} - n_\ell - k_{\ell + 1}
+ k_\ell } } { \bi{N - n_\ell}{n_{\ell + 1} - n_\ell} }. \ee} The conditional probability in (\ref{confinite}) can be viewed as the probability
of seeing $k_{\ell + 1} - k_\ell$ units having a certain attribute in the course of drawing $n_{\ell + 1} - n_\ell$ units, based on a simple
sampling without replacement, from a population of $N - n_\ell$ units, among which  $p N - k_\ell$ units having the attribute.  Actually, as can
be seen from Appendix \ref{multsamwithout_app}, the recursive formulae (\ref{recur1}) and (\ref{confinite}) for multistage sampling without
replacement can be proved by virtue of the notion of probability space.

It should be noted that the domain truncation technique to be described in subsection \ref{secDT} can be used to significantly reduce
computation.

\subsection{Domain Truncation} \la{secDT}

In the design and analysis of multistage sampling schemes, a frequent problem is to compute probabilistic terms like $\Pr \{ \mscr{W}
(\wh{\bs{\se}}, \mbf{n} ) \in \mscr{R} \}$, where $\mscr{W}$ is a bivariate function and $\mscr{R}$ is a subset of real numbers.  The
computational complexity associated with this type of problems can be high because the domain of summation or integration is large.  The
truncation techniques recently established in \cite{Chen1} have the power to considerably simplify the computation by reducing the domain of
summation or integration to a much smaller set. The following result derived from a similar method as that of \cite{Chen1}, shows that the
truncation can be done with controllable error.

\beT \la{Trun_THM5} Let $\eta \in (0, 1)$. Let $\udl{\se}_\ell, \; \ovl{\se}_\ell, \; \ell = 1, \cd, s$ be real numbers such that $\Pr \{
\udl{\se}_\ell \leq \wh{\bs{\se}}_\ell \leq \ovl{\se}_\ell \; \tx{for} \; \ell = 1, \cd, s \} \geq 1 - \eta$.  Assume that there exist subsets
of real numbers $\mscr{A}_\ell, \; \ell = 1, \cd, s$ such that $\{ \bs{l} = \ell \} = \{ \wh{\bs{\se}}_i \in \mscr{A}_i \; \tx{for} \; 1 \le i
\leq \ell \}$ for $\ell = 1, \cd, s$.  Then, {\small \be \la{good8} \Pr \{ \mscr{W} (\wh{\bs{\se}}, \mbf{n} ) \in \mscr{R} \} - \eta \leq
\sum_{\ell = 1}^s \Pr \{ \mscr{W} (\wh{\bs{\se}}_\ell, \mbf{n}_\ell ) \in \mscr{R} \; \tx{and} \; \wh{\bs{\se}}_i \in \mscr{B}_i \; \tx{for} \;
1 \leq i \leq \ell \} \leq \Pr \{ \mscr{W} (\wh{\bs{\se}}, \mbf{n} ) \in \mscr{R} \}, \ee} where $\mscr{B}_\ell = \{ \vse \in \mscr{A}_\ell: \;
\udl{\se}_\ell \leq \vse \leq \ovl{\se}_\ell  \}$ for $\ell = 1, \cd, s$.  \eeT

To determine numbers $\udl{\se}_\ell, \; \ovl{\se}_\ell, \; \ell = 1, \cd, s$ such that  $\Pr \{ \udl{\se}_\ell \leq \wh{\bs{\se}}_\ell \leq
\ovl{\se}_\ell \; \tx{for} \; \ell = 1, \cd, s \} \geq 1 - \eta$, we can follow a similar method as that of \cite{Chen1}.

\section{Construction of Sampling Schemes}

In this section, we shall discuss the applications of the
fundamental principle described in the previous section to the
design and analysis of multistage testing plans.

\subsection{Tests of Simple Hypotheses}

In some situations, it may be interesting to test multiple simple hypotheses $\mscr{H}_i : \se = \se_i$ for $i = 0, 1, \cd, m - 1$. For risk
control purpose, it is typically required that, for prescribed numbers $\de_i \in (0, 1)$, \be \la{simpreq}
 \Pr \li \{ \tx{Accept} \; \mscr{H}_i
\mid \se_i \ri \} \geq 1 - \de_i, \qqu i = 0, 1, \cd, m - 1. \ee Applying Theorem \ref{Multi_Comp_Exact} to the following hypotheses
\[
\mcal{H}_0: \se \leq \vse_1, \qu \mcal{H}_1: \vse_1 < \se \leq \vse_2, \qu \ldots, \qu \mcal{H}_{m-2}: \vse_{m-2} < \se \leq \vse_{m-1}, \qu
\mcal{H}_{m-1}: \se > \vse_{m-1}
\]
with $\vse_i = \f{ \se_{i - 1} + \se_i} {2}, \; i = 1, \cd, m-1$ and indifference zone $\cup_{i=1}^{m-1} (\se_{i-1}, \se_i)$, we have the
following results.

\beC \la{Multi_Simple_Exact} Let $\al_i, \ba_i \in (0, 1)$ for $i = 1, \cd, m - 1$ and $\al_m = \ba_0 = 0$. Define $\ovl{\al}_i = \max \{ \al_j:
i < j \leq m \}$ and $\ovl{\ba}_i = \max \{ \ba_j: 0 \leq j \leq i \}$ for $i = 0, 1, \cd, m - 1$. Suppose that $\bs{\varphi}_n$ is a ULCE of
$\se$ and that the maximum sample size $n_s$ is no less than  the minimum integer $n$ such that $f (n, \se_i, \ba_i) \geq g (n, \se_{i-1},
\al_i)$ for $i = 1, \cd, m-1$. Define $f_{\ell, i} = \udl{f} (n_\ell, \se_{i-1}, \se_i, \al_i, \ba_i)$ and $g_{\ell, i} = \ovl{g} (n_\ell,
\se_{i-1}, \se_i, \al_i, \ba_i)$ for $i = 1, \cd, m-1$. Define decision variable $\bs{D}_\ell$ by (\ref{defgoog}) for $\ell = 1, \cd, s$. Then,
$\Pr \{ \tx{Reject} \; \mscr{H}_i \mid \se_i \} \leq s ( \ovl{\al}_i + \ovl{\ba}_i )$ for $i = 0, 1, \cd, m - 1$.
 \eeC

Applying Theorem \ref{Multi_Comp_Chernoff} to hypotheses $\mcal{H}_i, \; i = 0, 1, \cd, m-1$ with indifference zone $\cup_{i=1}^{m-1}
(\se_{i-1}, \se_i)$, we have the following results.

\beC \la{Multi_Simple_Chernoff} Let $\al_i, \ba_i \in (0, 1)$ for $i = 1, \cd, m - 1$ and $\al_m = \ba_0 = 0$. Define $\ovl{\al}_i = \max \{
\al_j: i < j \leq m \}$ and $\ovl{\ba}_i = \max \{ \ba_j: 0 \leq j \leq i \}$ for $i = 0, 1, \cd, m - 1$. Suppose that \be \la{vg8} f_c (n,
\se_i, \ba_i) \geq g_c (n, \se_{i-1}, \al_i), \qqu i = 1, \cd, m-1 \ee if $n$ is sufficiently large. Suppose that the maximum sample size $n_s$
is no less than  the minimum integer $n$ such that (\ref{vg8}) is satisfied.  Define $f_{\ell, i} = \udl{f}_c (n_\ell, \se_{i-1}, \se_i, \al_i,
\ba_i)$ and $g_{\ell, i} = \ovl{g}_c (n_\ell, \se_{i-1}, \se_i, \al_i, \ba_i)$ for $i = 1, \cd, m-1$. Define decision variable $\bs{D}_\ell$ by
(\ref{defgoog}) for $\ell = 1, \cd, s$. Suppose that $\wh{\bs{\se}}_\ell$ is an unbiased and unimodal-likelihood estimator of $\se$ for $\ell =
1, \cd, s$. Then, $\Pr \{ \tx{Reject} \; \mscr{H}_i \mid \se_i \} \leq s ( \ovl{\al}_i + \ovl{\ba}_i )$ for $i = 0, 1, \cd, m - 1$.

\eeC

Corollaries  \ref{Multi_Simple_Exact} and \ref{Multi_Simple_Chernoff} provide methods to define testing plans and reveal that the risk
requirement \ref{simpreq} can be satisfied by choosing sufficiently small $\al_i, \; \ba_i$ for $i = 1, \cd, m - 1$.  The concrete determination
of such parameters is addressed in the sequel.

\subsubsection{Risk Tuning and Minimax Optimization} \la{minmaxcomsimple}

As can be seen from Corollaries \ref{Multi_Simple_Exact} and \ref{Multi_Simple_Chernoff}, to construct efficient testing plans satisfying the
risk requirement,  it is necessary to choose appropriate forms for $\al_i, \ba_i, \; i = 1, \cd, m - 1$ as functions of $\ze$. Specifically, we
propose to choose  \be \la{weight2}
 \al_i = \ze a_i, \qqu \ba_i = \ze b_i \qqu \tx{for} \; \; i = 1, \cd, m-1, \ee where the constants $a_i, b_i, \; i = 1, \cd, m-1$ are referred to as the weighting
coefficients.  Using  $\al_i, \; \ba_i$ defined by (\ref{weight2}), we can apply Corollaries \ref{Multi_Simple_Exact} and
\ref{Multi_Simple_Chernoff} to define testing plans.  For purpose of efficiency, it is crucial to determine appropriate values for the risk
tuning parameter $\ze$ and the weighting coefficients $a_i, b_i, \; i = 1, \cd, m-1$.  This task can be formulated as a minimax optimization
problem as follows.

Consider a family of testing plans associated with $\al_i, \; \ba_i$ defined by (\ref{weight2}).   Define
\[
A_i = \f{ \de_{i} } { \Pr \{ \tx{Reject} \; \mscr{H}_i  \mid \se_i \} }  \qqu \tx{for $i = 0, 1, \cd, m-1$}
\]
and
\[
Q = \max_{i \in \{ 0, 1, \cd, m - 1 \} } A_i, \qqu \qqu R = \min_{i \in \{ 0, 1, \cd, m - 1 \} } A_i.
\]
Clearly, $A_i, \; i = 0, 1, \cd, m-1$ are functions of $\ze$ and $a_i, \; b_i, \; i = 1, \cd, m-1$.  It can be seen that if $R \geq 1$, then the
risk requirement (\ref{simpreq}) is satisfied.  Moreover, under the restriction that $R \geq 1$, if $Q$ is smaller, then the associated testing
plan is more efficient.  Hence, we propose the following minimization problem:

{\it Determine risk tuning parameter $\ze$ and weighting coefficients $a_i, \; b_i, \; i = 1, \cd, m-1$ such that $Q$ is minimized under the
constraint that $R$ is no less than $1$}.

Clearly, this is a minimax optimization problem, since the quantity to be minimized is also a maximum over a discrete set.  For notational
convenience, we define $a_m = b_0 = 0$ and use it for describing our algorithm throughout the remainder of this section. To resolve the minimax
optimization problem, we propose the following iterative algorithm, which is motivated by the intuition that for $i = 0, 1, \cd, m-1$, $\Pr \{
\tx{Reject} \; \mscr{H}_i  \mid \se_i \}$ is ``roughly'' increasing with respect to $a_{i+1} + b_i$.

\bed

\item [\underline{Step 1}]: Set the maximum number of iterations as
$k_{\mrm max}$.   Choose the initial values of weighting coefficients as $a_1 = \de_0, \; b_{m-1} = \de_{m-1}$ and $a_{i+1} = b_i =
\f{\de_i}{2}$ for $i = 1, \cd, m - 2$.  Let $\wh{Q} \leftarrow \iy$ and $k \leftarrow 0$.

\item  [\underline{Step 2}]:  While $k \leq k_{\mrm max}$, do the following:

\bed

\item [\underline{Step 2-1}]: Based on $\al_i = \ze a_i, \;  \ba_i = \ze b_i$ for $i = 1, \cd, m-1$,  use a bisection search method to determine $\ze > 0$ as large as possible such that the value of $R$
associated with $\ze$ and $a_i, \; b_i, \; i = 1, \cd, m-1$ is no less than $1$. Let the value of $\ze$ obtained at this step be denoted by
$\ze^*$.  Let $Q^*, \; R^*$ and $A_i^*, \; i = 0, 1, \cd, m-1$ respectively denote the corresponding values of $Q, \; R$ and $A_i, \; i = 0, 1,
\cd, m-1$, which are determined by $\ze^*$ and $a_i, \; b_i, \; i = 0, 1, \cd, m-1$.

\item [\underline{Step 2-2}]: If $Q^* < \wh{Q}$, then let $\wh{Q} \leftarrow Q^*$ and
$\wh{a}_i \leftarrow \ze^* a_i, \; \; \wh{b}_i  \leftarrow \ze^* b_i$ for $i = 1, \cd, m - 1$.

\item [\underline{Step 2-3}]: For index $j \in \{0, 1, \cd, m-1\}$ such that
$A_j^* = Q^*$, let $a_{j+1} \leftarrow \ze^* a_{j+1} \li ( 1 + \f{ Q^* - 1 }{5} \ri )$ and $b_j \leftarrow \ze^* b_j \li ( 1 + \f{ Q^* - 1 }{5}
\ri )$.

\item [\underline{Step 2-4}]:  Let $k \leftarrow k + 1$.

\eed

\item [\underline{Step 3}]: Return $\ze = 1$ as the desired value of risk tuning parameter and $\wh{a}_i, \wh{b}_i, \; i = 1, \cd, m - 1$ as the desired weighting coefficients.

 \eed

Clearly, the above iterative algorithm returns weighting coefficients such that $Q$ is approximately minimized with $\ze = 1$ and $a_i =
\wh{a}_i, \; b_i = \wh{b}_i, \; i = 1, \cd, m - 1$ subject to the constraint that $R \geq 1$.  Hence, the iterative algorithm accomplishes the
task of coverage tuning.

\subsection{One-sided Tests}  \la{secones}

In order to infer from random samples $X_1, X_2, \cd$ of $X$ whether
the true value of $\se$ is greater or less than a certain number
$\vse \in {\Se}$, a classical problem is to test one-sided
hypothesis $\mscr{H}_0 : \se \leq \vse$ versus $\mscr{H}_1 : \se >
\vse$.  This problem can be cast in the general formulation
(\ref{mainpr}) with $m = 2, \; {\Se}_0 = \{ \se \in {\Se}: \se \leq
\vse \}$ and ${\Se}_1 = \{ \se \in {\Se}: \se
> \vse \}$. To control the probabilities of making wrong decisions, it
is typically required that, for {\it a priori} numbers $\al, \ba \in
(0, 1)$, \bel & & \Pr \li \{ \tx{Reject} \; \mscr{H}_0 \mid \se \ri
\} \leq \al \qu \tx{for any $\se \in
\varTheta_0$}, \la{re1gen} \\
&  & \Pr \li \{ \tx{Accept} \; \mscr{H}_0 \mid \se \ri \} \leq \ba
\qu \tx{for any $\se \in \varTheta_1$} \la{re2gen} \eel with
$\varTheta_0 = \{ \se \in {\Se}_0: \se \leq \se_0 \}$ and
$\varTheta_1 = \{ \se \in {\Se}_1: \se \geq \se_1 \}$, where $\se_0$
and $\se_1$ are numbers in ${\Se}$ such that $\se_0 < \vse < \se_1$.
The inequalities in (\ref{re1gen}) and (\ref{re2gen}) specify,
respectively,  the upper bounds for the probabilities of committing
a Type I error and a Type II error.  Clearly, the interval $(\se_0,
\se_1)$ is an indifference zone, since there is no requirement
imposed on probabilities of committing decision errors for $\se \in
(\se_0, \se_1)$.

Applying Theorem \ref{Multi_Comp_Exact} to the special case of $m =
2$, we have the following results.

\beC \la{One_sided_Exact_Lmit}    Let $\al_0, \ba_1 \in (0, 1)$.
 Suppose that $\bs{\varphi}_n$ is a ULCE of $\se$ and
that the maximum sample size $n_s$ is no less than the minimum
integer $n$ such that $f (n, \se_1, \ba_1) \geq g (n, \se_0,
\al_0)$. Define
\[ \bs{D}_\ell = \bec 1 & \tx{if} \;
\wh{\bs{\se}}_\ell \leq \udl{f} (n_\ell, \se_0,
\se_1, \al_0, \ba_1),\\
2 & \tx{if} \; \wh{\bs{\se}}_\ell > \ovl{g} (n_\ell, \se_0,
\se_1, \al_0, \ba_1),\\
0 & \tx{else} \eec \] for $\ell = 1, \cd, s$. Then, $\Pr \{
\tx{Accept} \; \mscr{H}_0 \mid \se \} \leq s \ba_1$ for $\se \in
{\Se}$ no less than $\se_1$, and $\Pr \{ \tx{Reject} \; \mscr{H}_0
\mid \se \} \leq s \al_0$ for $\se \in {\Se}$ no greater than
$\se_0$. Moreover, $\Pr \{ \tx{Accept} \; \mscr{H}_0 \mid \se \}$ is
non-increasing with respect to $\se \in {\Se}$ such that $\se \notin
(\se_0, \se_1)$.

\eeC

Applying Theorem \ref{Multi_Comp_Chernoff} to the special case of $m
= 2$, we have the following results.

 \beC \la{One_Sided_Chernoff} Let $\al_0, \ba_1 \in (0, 1)$.
 Suppose that $f_c (n, \se_1, \ba_1) \geq g_c (n, \se_0, \al_0)$ if $n$
is sufficiently large. Suppose that the maximum sample size $n_s$ is
no less than  the minimum integer $n$ such that $f_c (n, \se_1,
\ba_1) \geq g_c (n, \se_0, \al_0)$.  Suppose that
$\wh{\bs{\se}}_\ell$ is an unbiased and unimodal-likelihood
estimator of $\se$ for $\ell = 1, \cd, s$. Define
\[ \bs{D}_\ell = \bec 1 & \tx{if} \;
\wh{\bs{\se}}_\ell \leq \udl{f}_c (n_\ell, \se_0,
\se_1, \al_0, \ba_1),\\
2 & \tx{if} \; \wh{\bs{\se}}_\ell > \ovl{g}_c (n_\ell, \se_0,
\se_1, \al_0, \ba_1),\\
0 & \tx{else} \eec \] for $\ell = 1, \cd, s$.  Then, the same
conclusion as that of Corollary \ref{One_sided_Exact_Lmit} holds
true.

\eeC

In order to develop a class of test plans with OC functions being
monotone in the overall parameter space $\Se$, we shall introduce
multivariate functions {\small \bee & & \wh{F} (n, \se, \de) = \bec
\max \{ z \in I_{ \bs{\varphi}_n } : F_{ \bs{\varphi}_n } (z, \se)
\leq \de \} & \tx{if} \; \{ F_{ \bs{\varphi}_n }
(\bs{\varphi}_n, \se) \leq \de \} \neq \emptyset,\\
- \iy & \tx{otherwise}  \eec\\
&  & \wh{G} (n, \se, \de) = \bec \min \{ z \in I_{ \bs{\varphi}_n }
: G_{ \bs{\varphi}_n } (z, \se) \leq \de \}
& \tx{if} \; \{ G_{ \bs{\varphi}_n } (\bs{\varphi}_n, \se) \leq \de \} \neq \emptyset,\\
\iy & \tx{otherwise}  \eec \eee} for $n \in \bb{N}, \; \se \in \Se,
\; \de \in (0, 1)$ and {\small \bee &  & \udl{F} (n, \se^\prime,
\se^{\prime \prime}, \de^\prime, \de^{\prime \prime} ) = \min \li \{
\wh{F} (n, \se^{\prime \prime}, \de^{\prime \prime}), \qu \f{1}{2} [
\wh{F} (n, \se^{\prime \prime}, \de^{\prime \prime})
+ \wh{G} (n, \se^{\prime}, \de^{\prime}) ] \ri \},\\
&  & \ovl{G} (n, \se^\prime, \se^{\prime \prime}, \de^\prime,
\de^{\prime \prime} ) = \max \li \{ \wh{G} (n, \se^{\prime},
\de^{\prime}), \qu \f{1}{2} [ \wh{F} (n, \se^{\prime \prime},
\de^{\prime \prime}) + \wh{G} (n, \se^{\prime}, \de^{\prime}) ] \ri
\} \eee} for $\se^\prime < \se^{\prime \prime}$ in $\Se$ and
$\de^\prime, \de^{\prime \prime} \in (0, 1)$.  Moreover, we need to
make use of the concept of monotone likelihood ratio.  The
likelihood ratio is said to be monotonically increasing with respect
to $\bs{\varphi}_n$ if, for arbitrary $\se^{\prime} < \se^{\prime
\prime}$ in ${\Se}$, the likelihood ratio $\f{ \Pr \{ X_i = x_i, \;
i = 1, \cd, n \mid \se^{\prime \prime} \} }{ \Pr \{ X_i = x_i, \; i
= 1, \cd, n \mid \se^{\prime} \} }$ (or $\f{ f_{X_1, \cd, X_n} (x_1,
\cd, x_n \mid \se^{\prime \prime} )  } { f_{X_1, \cd, X_n} (x_1,
\cd, x_n \mid \se^{\prime} ) }$ for the continuous case) is
monotonically increasing with respect to $\bs{\varphi}_n$.

Now we are ready to describe a new class of test plans by Theorem
\ref{One_sided_Exact} as follows.

\beT \la{One_sided_Exact}    Let $\al_0, \ba_1 \in (0, 1)$. Suppose
that $\bs{\varphi}_n$ is a ULCE of $\se$ and that the likelihood
ratio is monotonically increasing with respect to $\bs{\varphi}_n$.
Suppose that the maximum sample size $n_s$ is no less than the
minimum integer $n$ such that $\wh{F} (n, \se_1, \ba_1) \geq \wh{G}
(n, \se_0, \al_0)$.  Define
\[ \bs{D}_\ell = \bec 1 & \tx{if} \; \wh{\bs{\se}}_\ell \leq \udl{F}
(n_\ell, \se_0,
\se_1, \al_0, \ba_1),\\
2 & \tx{if} \; \wh{\bs{\se}}_\ell > \ovl{G} (n_\ell, \se_0,
\se_1, \al_0, \ba_1),\\
0 & \tx{else} \eec \] for $\ell = 1, \cd, s$. Then, $\Pr \{
\tx{Accept} \; \mscr{H}_0 \mid \se \} \leq s \ba_1$ for $\se \in
{\Se}$ no less than $\se_1$, and $\Pr \{ \tx{Reject} \; \mscr{H}_0
\mid \se \} \leq s \al_0$ for $\se \in {\Se}$ no greater than
$\se_0$. Moreover, $\Pr \{ \tx{Accept} \; \mscr{H}_0 \mid \se \}$ is
non-increasing with respect to $\se \in {\Se}$.

\eeT

See Appendix \ref{One_sided_Exact_Ap} for a proof.

\bsk

In Theorem \ref{One_sided_Exact},  Corollaries \ref{One_sided_Exact_Lmit} and \ref{One_Sided_Chernoff},  we have proposed different testing
plans for testing the one-sided hypotheses.  To satisfy the risk requirements (\ref{re1gen}) and (\ref{re2gen}), we need to perform risk tuning
procedure.  Specifically, we propose to choose $\al_0 = \ze a_0, \; \ba_1 = \ze b_1$, where $\ze$ is the risk tuning parameter and $a_0, \; b_1$
are weighting coefficients.  Due to the monotonicity of the OC function, it suffices to ensure the risk requirements (\ref{re1gen}) and
(\ref{re2gen}) for parametric values $\se_0$ and $\se_1$.  Consider the family of testing plans associated with $\al_0 = \ze a_0$ and $\ba_1 =
\ze b_1$.   Define
\[
A = \f{ \al } { \Pr \{ \tx{Reject} \; \mscr{H}_0  \mid \se_0  \} }, \qu B = \f{ \ba } { \Pr \{ \tx{Reject} \; \mscr{H}_1 \mid \se_1 \} }, \qu Q
=  \max \li \{ A, \; B  \ri \}, \qu R = \min \li \{ A, \; B \ri \}
\]
as functions of $a_0, \; b_1$ and $\ze$.  For purpose of developing an efficient testing plan satisfying the risk requirement, we propose to
determine risk tuning parameter $\ze$ and weighting coefficients $a_0, \; b_1$ such that $Q$ is minimized under the constraint that $R$ is no
less than $1$. This task can be accomplished by applying the iterative minimax optimization algorithm developed in Section \ref{minmaxcom} to
 the present special problem.  The adapted algorithm is as follows.

\bed

\item [\underline{Step 1}]: Set the maximum number of iterations as
$k_{\mrm max}$.    Choose the initial values of weighting coefficients as $a_0 = \al$ and $b_1 = \ba$. Let $\wh{Q} \leftarrow \iy$ and $k
\leftarrow 0$.

\item  [\underline{Step 2}]:  While $k \leq k_{\mrm max}$, do the following:

\bed

\item [\underline{Step 2-1}]: Based on $\al_0 = \ze a_0$ and $\ba_1 = \ze b_1$,  use a bisection search method to determine $\ze > 0$ as large as possible such that the value of $R$ associated with $a_0, \;
b_1$ and $\ze$ is no less than $1$. Let the value of $\ze$ obtained at this step be denoted by $\ze^*$.  Let $Q^*, \; R^*$ and $A^*, \; B^*$
respectively denote the corresponding values of $Q, \; R$ and $A, \; B$.

\item [\underline{Step 2-2}]: If $Q^* < \wh{Q}$, then let
$\wh{a}_0 \leftarrow \ze^* a_0, \; \; \wh{b}_1  \leftarrow \ze^* b_1$ and $\wh{Q} \leftarrow Q^*$.

\item [\underline{Step 2-3}]: If $A^* = Q^*$, then let $a_0 \leftarrow \ze^* a_0 \li ( 1 + \f{ Q^* - 1 }{5} \ri )$. Otherwise,
let $b_1 \leftarrow \ze^* b_1 \li ( 1 + \f{ Q^* - 1 }{5} \ri )$.

\item [\underline{Step 2-4}]:  Let $k \leftarrow k + 1$.

\eed

\item [\underline{Step 3}]: Return $\ze = 1$ as the desired risk tuning parameter and $\wh{a}_0, \wh{b}_1$ as the desired weighting coefficients.

\eed

\subsection{Two-sided Tests} \la{sec Twos}

In order to infer from random samples $X_1, X_2, \cd$ of $X$ whether
the true value of $\se$ is equal to a certain number $\se_1 \in
{\Se}$, it is a frequent problem to test two-sided hypothesis
$\mscr{H}_0: \se = \se_1$ versus $\mscr{H}_1: \se \neq \se_1$.  To
control the probabilities of making wrong decisions, it is typically
required that, for {\it a priori} numbers $\al, \ba \in (0, 1)$,
\bel &  & \Pr \li \{ \tx{Reject} \; \mscr{H}_0 \mid \se_1 \ri \}
\leq \al,
\la{re1genTwo}\\
&  & \Pr \li \{ \tx{Accept} \; \mscr{H}_0 \mid \se \ri \} \leq \ba
\qu \tx{for $\se \in \Se$ such that $\se \notin (\se_0, \se_2 )$ },
\la{re2genTwo} \eel  where $\se_0$ and $\se_2$ are two numbers in
${\Se}$ such that $\se_0 < \se_1 < \se_2$. The inequalities in
(\ref{re1genTwo}) and (\ref{re2genTwo}) specify, respectively, the
upper bounds for the probabilities of committing a Type I error and
a Type II error. Since there is no requirement imposed on
probabilities of committing errors for $\se \in (\se_0, \se_1) \cup
(\se_1, \se_2)$, the union of intervals $(\se_0, \se_1) \cup (\se_1,
\se_2)$ is referred to as an indifference zone.

Applying Theorem \ref{Multi_Comp_Exact} to test hypotheses
\[
\mcal{H}_0: \se \leq \f{ \se_0 + \se_1} {2}, \qu \mcal{H}_1: \f{
\se_0 + \se_1} {2} < \se \leq \f{ \se_1 + \se_2} {2}, \qu
\mcal{H}_2: \se > \f{ \se_1 + \se_2} {2}
\]
with indifference zone $(\se_0, \se_1) \cup (\se_1, \se_2)$, we have
$\Pr \{ \tx{Reject} \; \mcal{H}_0 \; \tx{and} \; \mcal{H}_2 \mid \se
\} = \Pr \{ \tx{Accept} \; \mscr{H}_0  \mid \se \}$ and the
following results follow immediately.

\beC \la{Two_Sided_Exact}   Let $\al_0, \al_1, \ba_1, \ba_2 \in (0,
1)$. Suppose that $\bs{\varphi}_n$ is a ULCE of $\se$ and that the
maximum sample size $n_s$ is no less than  the minimum integer $n$
such that $f (n, \se_1, \ba_1) \geq g (n, \se_0, \al_0)$ and $f (n,
\se_2, \ba_2) \geq g (n, \se_1, \al_1)$. Define
\[ \bs{D}_\ell = \bec
1 & \tx{if} \; \ovl{g} (n_\ell, \se_0, \se_1, \al_0, \ba_1)
< \wh{\bs{\se}}_\ell \leq \udl{f} (n_\ell, \se_1, \se_2, \al_1, \ba_2),\\
2 & \tx{if} \; \wh{\bs{\se}}_\ell \leq \udl{f} (n_\ell, \se_0,
\se_1, \al_0, \ba_1) \; \tx{or} \; \wh{\bs{\se}}_\ell > \ovl{g}
(n_\ell, \se_1, \se_2, \al_1, \ba_2),\\
0 & \tx{else} \eec \] for $\ell = 1, \cd, s$. Then, $\Pr \{
\tx{Accept} \; \mscr{H}_0  \mid \se \} \leq s \times \max \{ \al_0,
\ba_2 \}$ for $\se \in \Se$ such that $\se \notin (\se_0, \se_2)$,
and $\Pr \{ \tx{Reject} \; \mscr{H}_0 \mid \se_1 \} \leq s ( \al_1 +
\ba_1)$. Moreover, $\Pr \{ \tx{Accept} \; \mscr{H}_0 \mid \se \}$ is
non-decreasing with respect to $\se \in {\Se}$ no greater than
$\se_0$ and is non-increasing with respect to $\se \in {\Se}$ no
less than $\se_2$.

\eeC

Applying Theorem \ref{Multi_Comp_Chernoff} to test hypotheses
$\mcal{H}_0, \; \mcal{H}_1$ and $\mcal{H}_2$ with indifference zone
$(\se_0, \se_1) \cup (\se_1, \se_2)$, we have the following results.

\beC \la{Two_Sided_Chernoff}

Let $\al_0, \al_1, \ba_1, \ba_2 \in (0, 1)$. Suppose that \be
\la{find}
 f_c (n, \se_1, \ba_1) \geq g_c (n, \se_0, \al_0), \qqu f_c (n, \se_2, \ba_2)
\geq g_c (n, \se_1, \al_1) \ee if $n$ is sufficiently large. Suppose
that the maximum sample size $n_s$ is no less than the minimum
integer $n$ such that (\ref{find}) is satisfied. Suppose that
$\wh{\bs{\se}}_\ell$ is an unbiased and unimodal-likelihood
estimator of $\se$ for $\ell = 1, \cd, s$. Define
\[ \bs{D}_\ell = \bec
1 & \tx{if} \; \ovl{g}_c (n_\ell, \se_0, \se_1, \al_0,
\ba_1)  < \wh{\bs{\se}}_\ell \leq \udl{f}_c (n_\ell, \se_1, \se_2, \al_1, \ba_2),\\
2 & \tx{if} \; \wh{\bs{\se}}_\ell \leq \udl{f}_c (n_\ell, \se_0,
\se_1, \al_0, \ba_1) \; \tx{or} \; \wh{\bs{\se}}_\ell > \ovl{g}_c
(n_\ell, \se_1, \se_2, \al_1, \ba_2),\\
0 & \tx{else} \eec \]  for $\ell = 1, \cd, s$.  Then, the same
conclusion as that of Corollary \ref{Two_Sided_Exact} holds true.

\eeC

In Corollaries \ref{Two_Sided_Exact} and \ref{Two_Sided_Chernoff},  we have proposed different testing plans for testing the two-sided
hypotheses. To satisfy the risk requirements (\ref{re1genTwo}) and (\ref{re2genTwo}), we need to develop a concrete risk tuning procedure.
Specifically, we propose to choose $\al_0 = \ze a_0, \; \al_1 = \ze a_1, \; \ba_1 = \ze b_1, \; \ba_2 = \ze b_2$, where $\ze$ is the risk tuning
parameter and $a_0, a_1, b_1, b_2$ are weighting coefficients.  Due to the monotonicity of the OC function, it suffices to ensure the risk
requirements (\ref{re1genTwo}) and (\ref{re2genTwo}) for parametric values $\se_0, \; \se_1$ and $\se_2$. Consider the family of testing plans
associated with $\al_0 = \ze a_0, \; \al_1 = \ze a_1, \; \ba_1 = \ze b_1, \; \ba_2 = \ze b_2$. Define
\[
Q_0 = \f{ \ba } { \Pr \{ \tx{Reject} \; \mscr{H}_1  \mid \se_0  \} }, \qqu Q_1 = \f{ \al } { \Pr \{ \tx{Reject} \; \mscr{H}_0 \mid \se_1 \} },
\qqu Q_2 = \f{ \ba } { \Pr \{ \tx{Reject} \; \mscr{H}_1 \mid \se_2 \} }
\]
\[ Q = \max \li \{ Q_0, \; Q_1, \; Q_2  \ri \}, \qqu R = \min \li \{ Q_0, \; Q_1, \; Q_2 \ri \}
\]
as functions of $a_0, a_1, b_1,  b_2$ and $\ze$.  In order to develop an efficient testing plan satisfying the risk requirements, we propose to
determine risk tuning parameter $\ze$ and weighting coefficients $a_0, a_1, b_1,  b_2$ such that $Q$ is minimized under the constraint that $R$
is no less than $1$. Our computational experiences indicate the truth of the following intuition:

(i) $\Pr \{ \tx{Reject} \; \mscr{H}_1  \mid \se_0  \}$ is ``roughly'' increasing with respect to $a_0$.

(ii) $\Pr \{ \tx{Reject} \; \mscr{H}_0  \mid \se_1  \}$ is ``roughly'' increasing with respect to $a_1 + b_1$.

(iii) $\Pr \{ \tx{Reject} \; \mscr{H}_1  \mid \se_2  \}$ is ``roughly'' increasing with respect to $b_2$.

Making use of the above intuition, we can adapt the iterative minimax optimization algorithm developed in Section \ref{minmaxcom} to solve our
present special problem.  The adapted algorithm is as follows.

\bed

\item [\underline{Step 1}]: Set the maximum number of iterations as
$k_{\mrm max}$.    Choose the initial values of weighting coefficients as $a_0 = b_2 = \ba$ and $a_1 = b_1 = \f{\al}{2}$. Let $\wh{Q} \leftarrow
\iy$ and $k \leftarrow 0$.

\item  [\underline{Step 2}]:  While $k \leq k_{\mrm max}$, do the following:

\bed

\item [\underline{Step 2-1}]: Based on $\al_0 = \ze a_0, \; \al_1 = \ze a_1, \; \ba_1 = \ze b_1, \; \ba_2 = \ze b_2$,
use a bisection search method to determine $\ze > 0$ as large as possible such that the value of $R$ associated with
$a_0, \; a_1, \; b_1, \; b_2$ and $\ze$ is no less than $1$. Let the value of $\ze$ obtained at this step be denoted by $\ze^*$. Let $Q^*, \;
R^*$ and $Q_0^*, \; Q_1^*, \; Q_2^*$ respectively denote the corresponding values of $Q, \; R$ and $Q_0, \; Q_1, \; Q_2$.

\item [\underline{Step 2-2}]: If $Q^* < \wh{Q}$, then let
$\wh{a}_0 \leftarrow \ze^* a_0, \; \wh{a}_1 \leftarrow \ze^* a_1, \;  \wh{b}_1  \leftarrow \ze^* b_1, \; \wh{b}_2 \leftarrow \ze^* b_2$ and
$\wh{Q} \leftarrow Q^*$.

\item [\underline{Step 2-3}]: If $Q_0^* = Q^*$, then let $a_0 \leftarrow \ze^* a_0 \li ( 1 + \f{ Q^* - 1 }{5} \ri )$.

\item [\underline{Step 2-4}]: If $Q_1^* = Q^*$, then let $a_1 \leftarrow \ze^* a_1 \li ( 1 + \f{ Q^* - 1 }{5} \ri )$
and $b_1 \leftarrow \ze^* b_1 \li ( 1 + \f{ Q^* - 1 }{5} \ri )$.

\item [\underline{Step 2-5}]: If $Q_2^* = Q^*$, then let $b_2 \leftarrow \ze^* b_2 \li ( 1 + \f{ Q^* - 1 }{5} \ri )$.

\item [\underline{Step 2-6}]:  Let $k \leftarrow k + 1$.

\eed

\item [\underline{Step 3}]: Return $\ze = 1$ as the desired risk tuning parameter
and $\wh{a}_0, \wh{a}_1, \wh{b}_1, \wh{b}_2$ as the desired weighting coefficients.

\eed

\subsection{Tests of Triple Hypotheses}  \la{secTriple}

As compared to two-sided tests, a more realistic formulation is to test three hypotheses $\mscr{H}_0: \se < \se_1, \; \mscr{H}_1: \se = \se_1$
and $\mscr{H}_2: \se > \se_1$, where $\se_1 \in \Se$. To control the risks of committing  decision errors, it is typically required that, for
prescribed numbers $\de_0, \de_1, \de_2 \in (0, 1)$, \bel & & \Pr \li \{ \tx{Accept} \; \mscr{H}_0 \mid \se \ri
\} \geq 1 - \de_0 \qu \tx{for  $\se \in \Se$ such that $\se \leq \se_0$,} \la{triska}\\
&  & \Pr \li \{ \tx{Accept} \;
\mscr{H}_1 \mid \se_1 \ri \} \geq 1 - \de_1, \la{triskb}\\
&  & \Pr \li \{ \tx{Accept} \; \mscr{H}_2 \mid \se \ri \} \geq 1 - \de_2 \qu \tx{for $\se \in \Se$ such that $\se \geq \se_2$,} \la{triskc} \eel
where $\se_0$ and $\se_2$ are numbers in ${\Se}$ such that $\se_0 < \se_1 < \se_2$.  Clearly, $(\se_0, \se_1) \cup (\se_1, \se_2)$ is an
indifference zone. Applying Theorem \ref{Multi_Comp_Exact} to test hypotheses {\small $\mcal{H}_0: \se \leq \f{ \se_0 + \se_1} {2}, \;
\mcal{H}_1: \f{ \se_0 + \se_1} {2} < \se \leq \f{ \se_1 + \se_2} {2}$} and {\small $\mcal{H}_2: \se > \f{ \se_1 + \se_2} {2}$} with indifference
zone $(\se_0, \se_1) \cup (\se_1, \se_2)$, we have the following results.

\beC \la{Exact_Triple}   Let $\al_0, \al_1, \ba_1, \ba_2 \in (0,
1)$.  Suppose that $\bs{\varphi}_n$ is a ULCE of $\se$. Suppose that
the maximum sample size $n_s$ is no less than the minimum integer
$n$ such that $f (n, \se_1, \ba_1) \geq g (n, \se_0, \al_0)$ and $f
(n, \se_2, \ba_2) \geq g (n, \se_1, \al_1)$. Define
\[ \bs{D}_\ell = \bec 1 & \tx{if} \;
\wh{\bs{\se}}_\ell \leq \udl{f} (n_\ell, \se_0, \se_1, \al_0, \ba_1),\\
2 & \tx{if} \; \ovl{g} (n_\ell, \se_0, \se_1, \al_0, \ba_1)
< \wh{\bs{\se}}_\ell \leq \udl{f} (n_\ell, \se_1, \se_2, \al_1, \ba_2),\\
3 & \tx{if} \; \wh{\bs{\se}}_\ell > \ovl{g} (n_\ell, \se_1,
\se_2, \al_1, \ba_2),\\
0 & \tx{else} \eec \]  for $\ell = 1, \cd, s$.   Then, the following
statements hold true.

(i) $\Pr \{ \tx{Reject} \; \mscr{H}_0  \mid \se \} \leq s \times
\max \{ \al_0, \al_1 \}$ for $\se \in {\Se}$ no greater than
$\se_0$. Moreover, $\Pr \{ \tx{Reject} \; \mscr{H}_0  \mid \se \}$
is non-decreasing with respect to $\se \in {\Se}$ no greater than
$\se_0$.

(ii) $\Pr \{ \tx{Reject} \; \mscr{H}_2  \mid \se \} \leq s \times
\max \{ \ba_1, \ba_2 \}$ for $\se \in {\Se}$ no less than $\se_2$.
Moreover, $\Pr \{ \tx{Reject} \; \mscr{H}_2  \mid \se \}$ is
non-increasing with respect to $\se \in {\Se}$ no less than $\se_2$.

(iii) $\Pr \{ \tx{Reject} \; \mscr{H}_1 \mid \se_1 \} \leq  s (\al_1
+ \ba_1)$.

\eeC

Applying Theorem \ref{Multi_Comp_Chernoff} to test hypotheses
{\small $\mcal{H}_0: \se \leq \f{ \se_0 + \se_1} {2}, \; \mcal{H}_1:
\f{ \se_0 + \se_1} {2} < \se \leq \f{ \se_1 + \se_2} {2}$} and
{\small $\mcal{H}_2: \se > \f{ \se_1 + \se_2} {2}$} with
indifference zone $(\se_0, \se_1) \cup (\se_1, \se_2)$, we have the
following results.

\beC \la{Triple_Chernoff}
 Let $\al_0, \al_1, \ba_1, \ba_2 \in (0, 1)$.  Suppose that \be \la{ggg}
 f_c (n, \se_1, \ba_1) \geq g_c (n, \se_0, \al_0), \qqu
 f_c (n, \se_2, \ba_2) \geq g_c (n, \se_1, \al_1) \ee
 if $n$ is sufficiently large.
 Suppose that the maximum sample size $n_s$ is no less than the minimum integer $n$
 such that (\ref{ggg}) is satisfied.  Suppose that
$\wh{\bs{\se}}_\ell$ is an unbiased and unimodal-likelihood
estimator of $\se$ for $\ell = 1, \cd, s$.  Define
\[ \bs{D}_\ell = \bec 1 & \tx{if} \;
\wh{\bs{\se}}_\ell \leq \udl{f}_c (n_\ell, \se_0, \se_1, \al_0, \ba_1),\\
2 & \tx{if} \; \ovl{g}_c (n_\ell, \se_0, \se_1, \al_0, \ba_1)
< \wh{\bs{\se}}_\ell \leq \udl{f}_c (n_\ell, \se_1, \se_2, \al_1, \ba_2),\\
3 & \tx{if} \; \wh{\bs{\se}}_\ell > \ovl{g}_c (n_\ell, \se_1,
\se_2, \al_1, \ba_2),\\
0 & \tx{else} \eec \] for $\ell = 1, \cd, s$.    Then, the same
conclusion as that of Corollary \ref{Exact_Triple} holds true.

\eeC

In Corollaries \ref{Exact_Triple} and \ref{Triple_Chernoff},  we have proposed different testing plans for testing the triple hypotheses. To
satisfy the risk requirements (\ref{triska}), (\ref{triskb}) and (\ref{triskc}), we need to perform risk tuning procedure.  Specifically, we
propose to choose $\al_0 = \ze a_0, \; \al_1 = \ze a_1, \; \ba_1 = \ze b_1, \; \ba_2 = \ze b_2$, where $\ze$ is the risk tuning parameter and
$a_0, a_1, b_1, b_2$ are weighting coefficients. Due to the monotonicity of the OC function, it suffices to ensure the risk requirements for
parametric values $\se_0, \; \se_1$ and $\se_2$. Consider the family of testing plans associated with $\al_0 = \ze a_0, \; \al_1 = \ze a_1, \;
\ba_1 = \ze b_1, \; \ba_2 = \ze b_2$. Define
\[
Q_0 = \f{ \de_0 } { \Pr \{ \tx{Reject} \; \mscr{H}_0  \mid \se_0  \} }, \qqu Q_1 = \f{ \de_1 } { \Pr \{ \tx{Reject} \; \mscr{H}_1 \mid \se_1 \}
}, \qqu Q_2 = \f{ \de_2 } { \Pr \{ \tx{Reject} \; \mscr{H}_2 \mid \se_2 \} }
\]
\[ Q = \max \li \{ Q_0, \; Q_1, \; Q_2  \ri \}, \qqu R = \min \li \{ Q_0, \; Q_1, \; Q_2 \ri \}
\]
as functions of $a_0, a_1, b_1,  b_2$ and $\ze$.  To obtain an efficient testing plan satisfying the risk requirements, we propose to determine
risk tuning parameter $\ze$ and weighting coefficients $a_0, a_1, b_1,  b_2$ such that $Q$ is minimized under the constraint that $R$ is no less
than $1$. Our computational experiences indicate the truth of the following intuition:

(i) $\Pr \{ \tx{Reject} \; \mscr{H}_0  \mid \se_0  \}$ is ``roughly'' increasing with respect to $a_0$.

(ii) $\Pr \{ \tx{Reject} \; \mscr{H}_1  \mid \se_1  \}$ is ``roughly'' increasing with respect to $a_1 + b_1$.

(iii) $\Pr \{ \tx{Reject} \; \mscr{H}_2  \mid \se_2  \}$ is ``roughly'' increasing with respect to $b_2$.

By virtue of the above intuition, we can adapt the iterative minimax optimization algorithm developed in Section \ref{minmaxcom} to solve our
present special problem.  The adapted algorithm is as follows.

\bed

\item [\underline{Step 1}]: Set the maximum number of iterations as
$k_{\mrm max}$.    Choose the initial values of weighting coefficients as $a_0 = \de_0, b_2 = \de_2$ and $a_1 = b_1 = \f{\de_1}{2}$. Let $\wh{Q}
\leftarrow \iy$ and $k \leftarrow 0$.

\item  [\underline{Step 2}]:  While $k \leq k_{\mrm max}$, do the following:

\bed

\item [\underline{Step 2-1}]: Based on $\al_0 = \ze a_0, \; \al_1 = \ze a_1, \; \ba_1 = \ze
b_1, \; \ba_2 = \ze b_2$,  use a bisection search method to determine $\ze > 0$ as large as possible such that the value of $R$ associated with
$a_0, \; a_1, \;  b_1, \; b_2$ and $\ze$ is no less than $1$. Let the value of $\ze$ obtained at this step be denoted by $\ze^*$. Let $Q^*, \;
R^*$ and $Q_0^*, \; Q_1^*, \; Q_2^*$ respectively denote the corresponding values of $Q, \; R$ and $Q_0, \; Q_1, \; Q_2$.

\item [\underline{Step 2-2}]: If $Q^* < \wh{Q}$, then let
$\wh{a}_0 \leftarrow \ze^* a_0, \; \wh{a}_1 \leftarrow \ze^* a_1, \;  \wh{b}_1  \leftarrow \ze^* b_1, \; \wh{b}_2 \leftarrow \ze^* b_2$ and
$\wh{Q} \leftarrow Q^*$.

\item [\underline{Step 2-3}]: If $Q_0^* = Q^*$, then let $a_0 \leftarrow \ze^* a_0 \li ( 1 + \f{ Q^* - 1 }{5} \ri )$.

\item [\underline{Step 2-4}]: If $Q_1^* = Q^*$, then let $a_1 \leftarrow \ze^* a_1 \li ( 1 + \f{ Q^* - 1 }{5} \ri )$
and $b_1 \leftarrow \ze^* b_1 \li ( 1 + \f{ Q^* - 1 }{5} \ri )$.

\item [\underline{Step 2-5}]: If $Q_2^* = Q^*$, then let $b_2 \leftarrow \ze^* b_2 \li ( 1 + \f{ Q^* - 1 }{5} \ri )$.

\item [\underline{Step 2-6}]:  Let $k \leftarrow k + 1$.

\eed

\item [\underline{Step 3}]: Return $\ze = 1$ as the desired risk tuning parameter
and $\wh{a}_0, \wh{a}_1, \wh{b}_1, \wh{b}_2$ as the desired weighting coefficients.

\eed

\subsection{Interval Tests} \la{secInterval Tests}

It is a frequent problem is to test hypothesis $\mscr{H}_0: \se \in [\se_1, \se_2]$ versus $\mscr{H}_1: \se \notin [\se_1, \se_2]$.  For risk
control purpose, it is typically required that, for two prescribed numbers $\al, \; \ba \in (0, 1)$, \bel &  & \Pr \li \{ \tx{Reject} \;
\mscr{H}_0 \mid \se \ri \} \leq \al \qu \tx{for $\se \in \Se$ such that $\se \in [\se_1^{\prime
\prime}, \se_2^\prime]$}, \la{intesta}\\
&  & \Pr \li \{ \tx{Accept} \; \mscr{H}_0 \mid \se \ri \} \leq \ba \qu \tx{for  $\se \in \Se$ such that $\se \notin (\se_1^\prime, \se_2^{\prime
\prime})$ },  \la{intestb} \eel where $\se_i^\prime, \se_i^{\prime \prime}$ are parametric values in $\Se$ such that $\se_1^\prime < \se_1 <
\se_1^{\prime \prime} < \se_2^\prime < \se_2 < \se_2^{\prime \prime}$. Since there is no requirement imposed on probabilities of committing
decision errors for $\se \in (\se_1^\prime, \se_1^{\prime \prime}) \cup (\se_2^\prime, \se_2^{\prime \prime})$, the union of intervals,
$(\se_1^\prime, \se_1^{\prime \prime}) \cup (\se_2^\prime, \se_2^{\prime \prime})$, is referred to as an indifference zone.

In view of the fact that the objective of the test is to decide
whether the parameter $\se$ falls into a specified interval, such a
test is called an ``interval test''.

Applying Theorem \ref{Multi_Comp_Exact} to test hypotheses
$\mcal{H}_0: \se \leq \se_1, \; \mcal{H}_1: \se_1 < \se \leq \se_2$
and $\mcal{H}_2: \se > \se_2$ with indifference zone $(\se_1^\prime,
\se_1^{\prime \prime}) \cup (\se_2^\prime, \se_2^{\prime \prime})$,
we  have  $\Pr \{ \tx{Reject} \; \mcal{H}_0 \; \tx{and} \;
\mcal{H}_2 \mid \se \} = \Pr \{ \tx{Accept} \; \mscr{H}_0  \}$ and
the following result follows immediately.

\beC \la{Exact_Interval} Let $\al_1, \al_2, \ba_1, \ba_2 \in (0,
1)$. Suppose that $\bs{\varphi}_n$ is a ULCE of $\se$ and that the
maximum sample size $n_s$ is no less than the minimum integer $n$
such that $f (n, \se_1^{\prime \prime}, \ba_1) \geq g (n,
\se_1^\prime, \al_1)$ and $f (n, \se_2^{\prime \prime}, \ba_2) \geq
g (n, \se_2^\prime, \al_2)$.  Define
\[ \bs{D}_\ell = \bec 1 & \tx{if} \; \ovl{g}
(n_\ell, \se_1^\prime, \se_1^{\prime \prime}, \al_1, \ba_1) <
\wh{\bs{\se}}_\ell \leq \udl{f} (n_\ell, \se_2^\prime, \se_2^{\prime
\prime}, \al_2, \ba_2),\\
2 & \tx{if} \; \wh{\bs{\se}}_\ell \leq \udl{f} (n_\ell,
\se_1^\prime, \se_1^{\prime \prime}, \al_1, \ba_1) \; \tx{or} \;
\wh{\bs{\se}}_\ell > \ovl{g} (n_\ell, \se_2^\prime,
\se_2^{\prime \prime}, \al_2, \ba_2),\\
0 & \tx{else} \eec \]  for $\ell = 1, \cd, s$.  Then, the following
statements hold true.

(i) $\Pr \{ \tx{Accept} \; \mscr{H}_0  \mid \se \} \leq s \times
\max \{ \al_1, \ba_2 \}$ for $\se \in \Se$ such that $\se \notin
(\se_1^\prime, \se_2^{\prime \prime})$.

(ii) $\Pr \{ \tx{Reject} \; \mscr{H}_0 \mid \se \} \leq  s (\al_2 +
\ba_1)$ for $\se \in \Se$ such that $\se \in [\se_1^{\prime \prime},
\se_2^\prime]$.

(iii) $\Pr \{ \tx{Accept} \; \mscr{H}_0  \mid \se \}$ is
non-decreasing with respect to $\se \in {\Se}$ no greater than
$\se_1^\prime$ and is non-increasing with respect to $\se \in {\Se}$
no less than $\se_2^{\prime \prime}$. Moreover, \bee & & \Pr \{
\tx{Reject} \; \mscr{H}_0  \mid \se \} \leq  \Pr \{ \tx{Reject} \;
\mscr{H}_0, \; \wh{\bs{\se}} \leq a \mid
a \} + \Pr \{ \tx{Reject} \; \mscr{H}_0, \; \wh{\bs{\se}} \geq b  \mid b \},\\
&  &  \Pr \{ \tx{Reject} \; \mscr{H}_0  \mid \se \} \geq \Pr \{
\tx{Reject} \; \mscr{H}_0, \; \wh{\bs{\se}} \leq a  \mid b \} + \Pr
\{ \tx{Reject} \; \mscr{H}_0, \; \wh{\bs{\se}} \geq b  \mid a \}
\eee for any $\se \in [a, b] \subseteq [\se_1^{\prime \prime},
\se_2^\prime] \cap \Se$.

\eeC

Applying Theorem \ref{Multi_Comp_Chernoff} to test hypotheses
$\mcal{H}_0: \se \leq \se_1, \; \mcal{H}_1: \se_1 < \se \leq \se_2$
and $\mcal{H}_2: \se > \se_2$ with indifference zone $(\se_1^\prime,
\se_1^{\prime \prime}) \cup (\se_2^\prime, \se_2^{\prime \prime})$,
we have the following results.

 \beC \la{Interval_Chernoff}

Let $\al_1, \al_2, \ba_1, \ba_2 \in (0, 1)$.  Suppose that \be
\la{net} f_c (n, \se_1^{\prime \prime}, \ba_1) \geq g_c (n,
\se_1^\prime, \al_1), \qqu f_c (n, \se_2^{\prime \prime}, \ba_2)
\geq g_c (n, \se_2^\prime, \al_2) \ee if $n$ is sufficiently large.
Suppose that the maximum sample size $n_s$ is no less than the
minimum integer $n$ such that (\ref{net}) is satisfied. Suppose that
$\wh{\bs{\se}}_\ell$ is an unbiased and unimodal-likelihood
estimator of $\se$ for $\ell = 1, \cd, s$. Define
\[ \bs{D}_\ell = \bec 1 & \tx{if} \;
\ovl{g}_c (n_\ell, \se_1^\prime, \se_1^{\prime \prime}, \al_1,
\ba_1) < \wh{\bs{\se}}_\ell \leq \udl{f}_c (n_\ell, \se_2^\prime,
\se_2^{\prime \prime}, \al_2, \ba_2),\\
2 & \tx{if} \; \wh{\bs{\se}}_\ell \leq \udl{f}_c (n_\ell,
\se_1^\prime, \se_1^{\prime \prime}, \al_1, \ba_1) \; \tx{or} \;
\wh{\bs{\se}}_\ell > \ovl{g}_c (n_\ell,
\se_2^\prime, \se_2^{\prime \prime}, \al_2, \ba_2),\\
0 & \tx{else} \eec \] for $\ell = 1, \cd, s$.  Then, the same
conclusion as that of Corollary \ref{Exact_Interval} holds true.

\eeC

In Corollaries \ref{Exact_Interval} and \ref{Interval_Chernoff},  we have proposed different testing plans for testing the hypotheses. To
satisfy the risk requirements (\ref{intesta}) and (\ref{intestb}),  we propose to choose $\al_1 = \ze a_1, \; \al_2 = \ze a_2, \; \ba_1 = \ze
b_1, \; \ba_2 = \ze b_2$, where $\ze$ is the risk tuning parameter and $a_1, a_2, b_1, b_2$ are weighting coefficients. If appropriate values
for weighting coefficients are available, then we can apply the bisection risk tuning procedure described in Section \ref{BRT} to determine $\ze
> 0$ as large as possible such that the risk requirements (\ref{intesta}) and (\ref{intestb}) are satisfied.

To determine appropriate values for the weighting coefficients, consider the family of testing plans associated with $\al_1 = \ze a_1, \; \al_2
= \ze a_2, \; \ba_1 = \ze b_1, \; \ba_2 = \ze b_2$. Define
\[
 A_1 = \f{ \al } { \Pr \{ \tx{Reject} \; \mscr{H}_0  \mid
\se_1^{\prime \prime}  \} }, \qqu A_2 = \f{ \al } { \Pr \{ \tx{Reject} \; \mscr{H}_0  \mid \se_2^{\prime}  \} },
\]
\[
 B_1 = \f{ \ba } { \Pr \{ \tx{Reject} \; \mscr{H}_1  \mid \se_1^{\prime}  \} }, \qqu B_2 = \f{ \ba } { \Pr \{ \tx{Reject} \; \mscr{H}_1
\mid \se_2^{\prime \prime}  \} }
\]
and
\[ Q = \max \li \{ A_1, \; A_2, \; B_1, \; B_2  \ri \}, \qqu R = \min \li \{ A_1, \; A_2, \; B_1, \; B_2 \ri \}
\]
as functions of $a_1, a_2, b_1,  b_2$ and $\ze$.  For purpose of efficiency, we propose to determine risk tuning parameter $\ze$ and weighting
coefficients $a_1, a_2, b_1,  b_2$ such that $Q$ is minimized under the constraint that $R$ is no less than $1$. Our computational experiences
indicate the truth of the following intuition:

(i) $\Pr \{ \tx{Reject} \; \mscr{H}_0  \mid \se_1^{\prime \prime}  \}$ is ``roughly'' increasing with respect to $b_1$.

(ii) $\Pr \{ \tx{Reject} \; \mscr{H}_0  \mid \se_2^{\prime}  \}$ is ``roughly'' increasing with respect to $a_2$.

(iii) $\Pr \{ \tx{Reject} \; \mscr{H}_1  \mid \se_1^{\prime}  \}$ is ``roughly'' increasing with respect to $a_1$.

(iv) $\Pr \{ \tx{Reject} \; \mscr{H}_1 \mid \se_2^{\prime \prime}  \}$ is ``roughly'' increasing with respect to $b_2$.

Making use of the above intuition, we can adapt the iterative minimax optimization algorithm developed in Section \ref{minmaxcom} to solve our
present special problem.  The adapted algorithm is as follows.

\bed

\item [\underline{Step 1}]: Set the maximum number of iterations as
$k_{\mrm max}$.    Choose the initial values of weighting coefficients as $a_1 = b_2 = \ba$ and $a_2 = b_1 = \al$. Let $\wh{Q} \leftarrow \iy$
and $k \leftarrow 0$.

\item  [\underline{Step 2}]:  While $k \leq k_{\mrm max}$, do the following:

\bed

\item [\underline{Step 2-1}]: Based on $\al_1 = \ze a_1, \; \al_2 = \ze a_2,
\; \ba_1 = \ze b_1, \; \ba_2 = \ze b_2$,  use a bisection search method to determine $\ze > 0$ as large as possible such that the value of $R$
associated with $a_1, \; a_2, \; b_1, \; b_2$ and $\ze$ is no less than $1$. Let the value of $\ze$ obtained at this step be denoted by $\ze^*$.
Let $Q^*, \; R^*$ and $A_1^*, \; A_2^*, \;  B_1^*, \; B_2^*$ respectively denote the corresponding values of $Q, \; R$ and $A_1, \; A_2, \; B_1,
\; B_2$.

\item [\underline{Step 2-2}]: If $Q^* < \wh{Q}$, then let
$\wh{a}_1 \leftarrow \ze^* a_1, \; \wh{a}_2 \leftarrow \ze^* a_2, \;  \wh{b}_1  \leftarrow \ze^* b_1, \; \wh{b}_2 \leftarrow \ze^* b_2$ and
$\wh{Q} \leftarrow Q^*$.

\item [\underline{Step 2-3}]: If $A_1^* = Q^*$, then let $b_1 \leftarrow \ze^* b_1 \li ( 1 + \f{ Q^* - 1 }{5} \ri )$.

\item [\underline{Step 2-4}]: If $A_2^* = Q^*$, then let $a_2 \leftarrow \ze^* a_2 \li ( 1 + \f{ Q^* - 1 }{5} \ri )$.

\item [\underline{Step 2-5}]: If $B_1^* = Q^*$, then let $a_1 \leftarrow \ze^* a_1 \li ( 1 + \f{ Q^* - 1 }{5} \ri )$.

\item [\underline{Step 2-6}]: If $B_2^* = Q^*$, then let $b_2 \leftarrow \ze^* b_2 \li ( 1 + \f{ Q^* - 1 }{5} \ri )$.

\item [\underline{Step 2-7}]:  Let $k \leftarrow k + 1$.

\eed

\item [\underline{Step 3}]: Return $\wh{a}_1, \wh{a}_2, \wh{b}_1, \wh{b}_2$ as the desired weighting coefficients.

\eed

Using the output of the above algorithm $\wh{a}_1, \wh{a}_2, \wh{b}_1, \wh{b}_2$ as the weighting coefficients to define $\al_1 = \ze \wh{a}_1,
\; \al_2 = \ze \wh{a}_2, \; \ba_1 = \ze \wh{b}_1, \; \ba_2 = \ze \wh{b}_2$, we can apply the bisection risk tuning technique described in
Section \ref{minmaxcom} to determine $\ze > 0$ as large as possible such that the risk requirements (\ref{intesta}) and (\ref{intestb}) are
satisfied.  Our computational experiments show that in many situations, the resultant $\ze$ is equal or very close to $1$.

\subsection{Applications} \la{appdenmo}

In this section, we shall demonstrate that the general principle
proposed above can be applied to develop specific test plans for
common important distributions.  To apply our general method, we
need to choose appropriate estimator $\bs{\varphi}_n = \varphi (X_1,
\cd, X_n)$ for $\se$ and investigate whether $\bs{\varphi}_n$ has
the following properties:

(i) $\bs{\varphi}_n$ is a ULE of $\se$;

(ii) $\bs{\varphi}_n$ converges in probability to $\se$;

(iii) $\bs{\varphi}_n$ is an unbiased estimator of $\se$;

(iv) The likelihood ratio is monotonically increasing with respect
to $\bs{\varphi}_n$;

(v) For $\se^{\prime} < \se^{\prime \prime}$ in ${\Se}$ and
$\de^{\prime}, \; \de^{\prime \prime} \in (0, 1)$, $f_c (n,
\se^{\prime \prime}, \de^{\prime \prime})$ is no less than $g_c (n,
\se^{\prime}, \de^{\prime})$ if $n$ is sufficiently large.

\subsubsection{Testing a Binomial Proportion}

Let $X$ be a Bernoulli random variable with distribution $\Pr \{ X =
1 \} = 1 - \Pr \{ X = 0 \} = p \in (0, 1)$. To test hypotheses
regarding $p$ based on i.i.d. samples $X_1, X_2, \cd$ of $X$, we
shall take $\bs{\varphi}_n = \varphi(X_1, \cd, X_n) = \f{\sum_{i =
1}^{n} X_i} { n }$ as an estimator of $p$.  With such a choice of
estimator, it can be shown that, for $n \in \bb{N}, p \in (0, 1),
\de \in (0, 1)$,  {\small \bee & & \wh{F} (n, p, \de) = \bec
\f{1}{n} \times \max \li \{ k \in \bb{Z}: \sum_{i = 0}^k \bi{n}{i}
p^i (1 - p)^{n - i} \leq \de, \; 0 \leq k \leq n \ri \}  & \tx{for}
\; n
\geq \f{ \ln (\de) } { \ln (1 - p) },\\
- \iy & \tx{for} \; n < \f{ \ln (\de) } { \ln (1 - p) } \eec\\
&  & \wh{G} (n, p, \de) = \bec \f{1}{n} \times \min \li \{ k \in
\bb{Z}: \sum_{i = k}^n \bi{n}{i} p^i (1 - p)^{n - i} \leq \de, \; 0
\leq k \leq n \ri \}  & \tx{for} \; n
\geq \f{ \ln (\de) } { \ln (p) },\\
\iy & \tx{for} \; n < \f{ \ln (\de) } { \ln (p) } \eec\\
&  & f (n, p, \de) = \bec \f{1}{n} \times \max \li \{ k \in \bb{Z}:
\sum_{i = 0}^k \bi{n}{i} p^i (1 - p)^{n - i} \leq \de, \; 0 \leq k
\leq np \ri \}  & \tx{for} \; n
\geq \f{ \ln (\de) } { \ln (1 - p) },\\
- \iy & \tx{for} \; n < \f{ \ln (\de) } { \ln (1 - p) } \eec\\
&  &  g (n, p, \de) = \bec \f{1}{n} \times \min \li \{ k \in \bb{Z}:
\sum_{i = k}^n \bi{n}{i} p^i (1 - p)^{n - i} \leq \de, \; np \leq k
\leq n \ri \}  & \tx{for} \; n
\geq \f{ \ln (\de) } { \ln (p) },\\
\iy & \tx{for} \; n < \f{ \ln (\de) } { \ln (p) } \eec \eee} and
{\small \bee &  & f_c (n, p, \de) = \bec \max \{ z \in [0, p]:
\mscr{M}_{\mrm{B}} (z, p) \leq \f{ \ln (\de) }{n} \} & \tx{for} \; n
\geq \f{ \ln (\de) } { \ln (1 - p) },\\
- \iy & \tx{for} \; n < \f{ \ln (\de) } { \ln (1 - p) } \eec\\
&  & g_c (n, p, \de) = \bec \min \{ z \in [p, 1]: \mscr{M}_{\mrm{B}}
(z, p) \leq \f{ \ln (\de) }{n} \}
& \tx{for}  \; n \geq \f{ \ln (\de) } { \ln (p) },\\
\iy & \tx{for}  \; n < \f{ \ln (\de) } { \ln (p) } \eec \eee} where
{\small \[ \mscr{M}_{\mrm{B}} (z,p) = \bec z \ln \f{p}{z} + (1 - z)
\ln \f{1 - p}{1 - z} &
\tx{for} \; z \in (0,1),\\
\ln(1-p) & \tx{for} \; z = 0,\\
\ln p &  \tx{for} \; z = 1. \eec
\]}
Moreover, it can be verified that the estimator $\bs{\varphi}_n$
possesses all properties described at the beginning of Section
\ref{appdenmo}.  This implies that all testing methods proposed in
previous sections are applicable.

\subsubsection{Testing the Proportion of a Finite Population} \la{Finite_Plan}

It is a frequent problem to test the proportion of a finite
population.  Consider a population of $N$ units, among which there
are $N p$ units having a certain attribute, where $p \in {\Se} = \{
\f{i}{N} : i = 0, 1, \cd, N \}$.  The procedure of sampling without
replacement can be described as follows:

 Each time a single unit is drawn without replacement from the remaining population so
that every unit of the remaining population has equal chance of being selected.

Such a sampling process can be exactly characterized by random
variables $X_1, \cd, X_N$ defined in a probability space $(\Omega,
\mscr{F}, \Pr)$ such that $X_i$ denotes the characteristics of the
$i$-th sample in the sense that $X_i = 1$ if the $i$-th sample has
the attribute and $X_i = 0$ otherwise.  By the nature of the
sampling procedure, it can be shown that \be \la{basis}
 \Pr \{ X_i =
x_i, \; i = 1, \cd, n \mid p \} = \bi{N p}{\sum_{i = 1}^n x_i} \bi{N
- N p}{n - \sum_{i = 1}^n x_i}
 \li \slash \li [ \bi{n}{\sum_{i = 1}^n x_i} \bi{N}{n} \ri. \ri ]
\ee for any $n \in \{1, \cd, N\}$ and any $x_i \in \{0, 1\}, \; i = 1, \cd, n$.  By virtue of (\ref{basis}), it can be shown that $\Pr \{ X_i =
1 \} = 1 - \Pr \{ X_i = 0 \} = p \in {\Se}$, which implies that $X_1, \cd, X_N$ can be treated as identical but dependent samples of a Bernoulli
random variable $X$ such that $\Pr \{ X = 1 \} = 1 - \Pr \{ X = 0 \} = p \in {\Se}$.  Recently, we have shown in \cite{Chen_EST} that, for any
$n \in \{1, \cd, N\}$, the sample mean {\small $\bs{\varphi}_n = \f{ \sum_{i=1}^n X_i }{n}$} is a ULE for $p \in {\Se}$. Clearly, {\small
$\bs{\varphi}_n$} is not a MLE for $p \in {\Se}$.  Hence, we can develop multistage testing plans in the framework outlined in Section
\ref{secSRI}. With the choice of {\small $\bs{\varphi}_n = \f{ \sum_{i=1}^n X_i }{n}$} as the estimator of $p$, it can be shown that {\small
\bee & & \wh{F} (n, p, \de) = \bec \f{1}{n} \times \max \li \{ k \in \bb{Z}: \sum_{i= 0 }^k {p N \choose i} {N - p N \choose n - i } \slash {N
\choose n}
\leq \de, \; 0 \leq k < n \ri \} & \tx{for} \; \bi{N - p N}{n} \leq \de \bi{N}{n},\\
- \iy &  \tx{for} \; \bi{N - p N}{n} > \de \bi{N}{n} \eec\\
&   & \wh{G} (n, p, \de) = \bec \f{1}{n} \times  \min \li \{ k \in
\bb{Z}: \sum_{i= k }^n {p N \choose i} {N - p N \choose n - i  }
\slash {N \choose n}
\leq \de, \; 0 < k \leq n \ri \} & \tx{for} \; \bi{p N}{n} \leq \de \bi{N}{n},\\
\iy &  \tx{for} \; \bi{p N}{n} > \de \bi{N}{n} \eec \\
&   & f (n, p, \de) = \bec \f{1}{n} \times  \max \li \{ k \in
\bb{Z}: \sum_{i= 0 }^k {p N \choose i} {N - p N \choose n - i }
\slash {N \choose n}
\leq \de, \; 0 \leq k \leq n p \ri \} & \tx{for} \; \bi{N - p N}{n} \leq \de \bi{N}{n},\\
- \iy &  \tx{for} \; \bi{N - p N}{n} > \de \bi{N}{n} \eec\\
&   & g (n, p, \de) = \bec \f{1}{n} \times  \min \li \{ k \in
\bb{Z}: \sum_{i= k }^n {p N \choose i} {N - p N \choose n - i  }
\slash {N \choose n}
\leq \de, \; n p \leq k \leq n \ri \} & \tx{for} \; \bi{p N}{n} \leq \de \bi{N}{n},\\
\iy &  \tx{for} \; \bi{p N}{n} > \de \bi{N}{n} \eec \eee} for $n \in
\{1, \cd, N \}, \; p \in {\Se}$ and $\de \in (0, 1)$.  Clearly,
$\bs{\varphi}_n$ converges in probability to $p$ and thus is a ULCE
of $p$.  Moreover, it can be verified that the likelihood ratio is
monotonically increasing with respect to $\bs{\varphi}_n$. This
implies that the general results described in the previous sections
can be useful.

In order to develop test plans with simple stopping boundary, we
define multivariate functions {\small \bee &  & f_c (n, p, \de) =
\bec \max \{ z \in I_{ \bs{\varphi}_n } : \mcal{C} (n, z, p) \leq
\de, \; z \leq p \} & \tx{if} \; \{ \mcal{C} (n, \bs{\varphi}_n,
p) \leq \de, \; \bs{\varphi}_n \leq p \} \neq \emptyset,\\
- \iy & \tx{otherwise}  \eec\\
&  & g_c (n, p, \de) = \bec \min \{ z \in I_{ \bs{\varphi}_n } :
\mcal{C} (n, z, p) \leq \de, \; z \geq p \} & \tx{if} \; \{ \mcal{C}
(n, \bs{\varphi}_n, p) \leq
\de, \; \bs{\varphi}_n \geq p \} \neq \emptyset,\\
\iy & \tx{otherwise}  \eec
 \eee}
for $n \in \bb{N}, p \in \Se, \de \in (0, 1)$, where \be \la{defc}
\mcal{C} (n, z, p) = \bec \f{ \bi{Np}{n} } { \bi{N}{n} } &
\tx{for} \; z = 1,\\
\f{ \bi{Np}{n z} \bi{N - Np}{n - n z} } { \bi{ \lf (N+1) z \rf }{n
z} \bi{N - \lf (N+1) z \rf }{n - n z} } & \tx{for} \; z \in \{
\f{k}{n} : k \in \bb{Z}, \; 0 \leq k < n \}.
 \eec
\ee Moreover, define {\small \bee & & \udl{f}_c (n, p^\prime,
p^{\prime \prime}, \de^\prime, \de^{\prime \prime} ) = \min \li \{
f_c (n, p^{\prime \prime}, \de^{\prime \prime}), \qu \f{1}{2} [ f_c
(n, p^{\prime \prime}, \de^{\prime \prime}) + g_c (n, p^{\prime},
\de^{\prime}) ] \ri \},\\
&  & \ovl{g}_c (n, p^\prime, p^{\prime \prime}, \de^\prime,
\de^{\prime \prime} ) = \max \li \{  g_c (n, p^{\prime},
\de^{\prime}), \qu \f{1}{2} [ f_c (n, p^{\prime \prime}, \de^{\prime
\prime}) + g_c (n, p^{\prime}, \de^{\prime}) ] \ri \} \eee} for
$p^\prime < p^{\prime \prime}$ in $\Se, \; \de^\prime, \de^{\prime
\prime} \in (0, 1)$ and $n \in \bb{N}$.

For the multi-hypothesis testing problem stated in the introduction
with $\se$ replaced by $p$, we have the following results.

 \beT
\la{Multi_Comp_Exact_Finite}

Let $\al_i, \ba_i \in (0, 1)$ for $i = 1, \cd, m - 1$ and $\al_m =
\ba_0 = 0$. Define $\ovl{\al}_i = \max \{ \al_j: i < j \leq m \}$
and $\ovl{\ba}_i = \max \{ \ba_j: 0 \leq j \leq i \}$ for $i = 0, 1,
\cd, m - 1$. Suppose that the maximum sample size $n_s$ is no less
the minimum integer $n$ such that $f_c (n, p_i^{\prime \prime},
\ba_i) \geq g_c (n, p_i^{\prime}, \al_i)$ for $i = 1, \cd, m-1$.
Define $f_{\ell, i} = \udl{f}_c (n_\ell, p_i^{\prime}, p_i^{\prime
\prime}, \al_i, \ba_i)$ and $g_{\ell, i} = \ovl{g}_c (n_\ell,
p_i^{\prime}, p_i^{\prime \prime}, \al_i, \ba_i)$ for $i = 1, \cd,
m-1$. Define {\small $\wh{\bs{p}}_\ell = \bs{\varphi}_{n_\ell} = \f{
\sum_{i=1}^{n_\ell} X_i }{n}$} and \be \la{defgoogfinite}
 \bs{D}_\ell = \bec 1 & \tx{if} \; \;
\wh{\bs{p}}_\ell \leq
f_{\ell ,1}, \\
i & \tx{if} \; \;  g_{\ell, i-1} < \wh{\bs{p}}_\ell \leq f_{\ell,i}
\; \tx{where} \; 2 \leq i \leq m - 1,\\
m  & \tx{if} \; \;  \wh{\bs{p}}_\ell > g_{\ell ,m-1},\\
0 & \tx{else} \eec \ee for $\ell = 1, \cd, s$.  The following
statements (I)-(V) hold true for $m \geq 2$.

(I) $\Pr \{ \tx{Reject} \; \mscr{H}_0 \mid p \}$ is non-decreasing
with respect to $p \in \varTheta_0$.

(II) $\Pr \{ \tx{Reject} \; \mscr{H}_{m-1} \mid p \} \; \tx{is
non-increasing with respect to} \; p \in \varTheta_{m-1}$.

(III)  $\Pr \{ \tx{Reject} \; \mscr{H}_i  \mid p \} \leq s (
\ovl{\al}_i  + \ovl{\ba}_i  )$ for any $p \in \varTheta_i$ and $i =
0, 1, \cd, m - 1$.

(IV) For $0 < i \leq m- 1$, $\Pr \{ \tx{Accept} \; \mscr{H}_i \mid p
\}$ is no greater than $s \al_i$ and is non-decreasing with respect
to $p \in {\Se}$ no greater than $p_i^\prime$.

(V) For $0 \leq i \leq m-2$, $\Pr \{ \tx{Accept} \; \mscr{H}_i \mid
p \}$ is no greater than $s \ba_{i + 1}$ and is non-increasing with
respect to $p \in {\Se}$ no less than $p_{i+1}^{\prime \prime}$.

 Moreover, the following statements (VI), (VII) and (VIII) hold true for $m \geq 3$.

(VI) \bee & & \Pr \{ \tx{Reject} \; \mscr{H}_i \mid p \} \leq \Pr \{
\tx{Reject} \; \mscr{H}_i, \; \wh{\bs{p}} \leq a \mid a \} +
\Pr \{ \tx{Reject} \; \mscr{H}_i, \; \wh{\bs{p}} \geq b \mid b \},\\
&  & \Pr \{ \tx{Reject} \; \mscr{H}_i \mid p \} \geq  \Pr \{
\tx{Reject} \; \mscr{H}_i, \; \wh{\bs{p}} \leq a \mid b \} + \Pr \{
\tx{Reject} \; \mscr{H}_i, \; \wh{\bs{p}} \geq b \mid a \} \eee for
any $p \in [a, b ] \subseteq \varTheta_i$ and $1 \leq i \leq m - 2$.

(VII) $\Pr \{ \tx{Reject} \; \mscr{H}_0 \; \tx{and} \;
\mscr{H}_{m-1} \mid p \}$ is non-decreasing with respect to $p \in
\varTheta_0$ and is non-increasing with respect to $p \in
\varTheta_{m-1}$.

(VIII) $\Pr \{ \tx{Reject} \; \mscr{H}_0 \; \tx{and} \;
\mscr{H}_{m-1} \mid p \}$ is no greater than $s \times \max \{
\al_i: 1 \leq i \leq m - 2 \}$ for $p \in \varTheta_0$ and is no
greater than $s \times \max \{ \ba_i: 2 \leq i \leq m - 1 \}$ for $p
\in \varTheta_{m-1}$.

\eeT

It should be noted that $p_i^\prime, \; p_i^{\prime \prime}$ in
Theorem \ref{Multi_Comp_Exact_Finite} play similar roles as
$\se_i^\prime, \; \se_i^{\prime \prime}$ in the introduction in
defining the requirement of risk control. Accordingly, $\varTheta_i$
in Theorem \ref{Multi_Comp_Exact_Finite} has the same notion as
$\varTheta_i$ in introduction with parameter $\se$ identified as
$p$.

Theorem \ref{Multi_Comp_Exact_Finite} can be shown by using a
similar argument as that for Theorem \ref{Multi_Comp_Exact} and the
following results obtained by Chen \cite{ChenR}, \bel &  & \Pr \li
\{ \f{ \sum_{i = 1}^n X_i }{n} \geq z \mid p \ri \} \leq \mcal{C}
(n, z, p) \qu \tx{for} \; z \in
\li \{ \f{k}{n}: k \in \bb{Z}, \; n p \leq k \leq n \ri \}, \la{hyp1}\\
&  & \Pr \li \{ \f{ \sum_{i = 1}^n X_i }{n}  \leq z \mid p \ri \}
\leq \mcal{C} (n, z, p) \qu \tx{for} \; z \in \li \{ \f{k}{n}: k \in
\bb{Z}, \; 0 \leq k \leq n p \ri \} \la{hyp2} \eel where $p \in \Se$
and $\mcal{C} (n, z, p)$ is defined by (\ref{defc}). Since $\sum_{i
= 1}^n X_i$ has a hypergeometric distribution, the above
inequalities (\ref{hyp1}) and (\ref{hyp2}) provide simple bounds for
the tail probabilities of hypergeometric distribution, which are
substaintially less conservative than Hoeffding's inequalities
\cite{Hoeffding}.

\subsubsection{Testing the Parameter of a Poisson Distribution}

Let $X$ be a Poisson variable of mean $\lm > 0$.  We shall consider
the test of hypotheses regarding $\lm$ based on i.i.d. random
samples $X_1, X_2, \cd$ of $X$.  Choosing {\small $\bs{\varphi}_n =
\f{ \sum_{i = 1}^n X_i } { n }$} as an estimator for $\lm$, we can
show that, for $n \in \bb{N}, \; \lm \in (0, \iy), \; \de \in (0,
1)$, {\small \bee & & \wh{F} (n, \lm, \de) = \bec \f{1}{n} \times
\max \li \{ k \in \bb{Z}: \sum_{i = 0}^k \f{ (n \lm)^i e^{-n \lm}
}{i!} \leq \de, \; k \geq 0 \ri \}  & \tx{for} \;
n \geq \f{ \ln (\de) } { - \lm },\\
- \iy & \tx{for} \; n < \f{ \ln (\de) } { - \lm } \eec\\
&  & \wh{G} (n, \lm, \de) = \f{1}{n} \times \min \li \{ k \in
\bb{Z}: \sum_{i = 0}^{k-1} \f{ (n \lm)^i
e^{- n \lm} }{i!} \geq 1 - \de, \; k \geq 1 \ri \} \\
&  & f (n, \lm, \de) = \bec \f{1}{n} \times \max \li \{ k \in
\bb{Z}: \sum_{i = 0}^k \f{ (n \lm)^i e^{-n \lm} }{i!} \leq \de, \; 0
\leq k \leq n \lm \ri \}  & \tx{for} \;
n \geq \f{ \ln (\de) } { - \lm },\\
- \iy & \tx{for} \; n < \f{ \ln (\de) } { - \lm } \eec\\
&  & g (n, \lm, \de) = \f{1}{n} \times \min \li \{ k \in \bb{Z}:
 \sum_{i = 0}^{k-1} \f{ (n \lm)^i
e^{- n \lm} }{i!} \geq 1 - \de, \; k \geq n \lm \ri \} \eee} and
{\small \bee &  &  f_c (n, \lm, \de) = \bec \max \{ z \in [0, \lm]:
\mscr{M}_{\mrm{P}} (z, \lm) \leq \f{ \ln (\de) }{n} \} & \tx{for}
\;  n \geq \f{ \ln (\de) } { - \lm },\\
- \iy & \tx{for} \;  n < \f{ \ln (\de) } { - \lm } \eec\\
& & g_c (n, \lm, \de) = \min \li \{ z \in [\lm, \iy):
\mscr{M}_{\mrm{P}} (z, \lm ) \leq \f{ \ln (\de) }{n}  \ri \} \eee}
where
\[
\mscr{M}_{\mrm{P}} (z, \lm) =
\bec z - \lm + z \ln \li ( \f{\lm}{z} \ri ) & \tx{for} \; z > 0,\\
- \lm & \tx{for} \; z = 0.  \eec
\]
Moreover, it can be verified that the estimator $\bs{\varphi}_n$
possesses all properties described at the beginning of Section
\ref{appdenmo}.  This implies that all testing methods proposed in
previous sections are applicable.

\subsubsection{Testing the Mean of a Normal Distribution with Known Variance}

It is an important problem to test the mean, $\mu$, of a Gaussian
random variable $X$ with known variance $\si^2$ based on i.i.d.
random samples $X_1, X_2, \cd$ of $X$. Choosing {\small
$\bs{\varphi}_n = \f{\sum_{i=1}^n X_i}{n}$} as an estimator of
$\mu$, we have
\[
\wh{F} (n, \mu, \de) = f (n, \mu, \de) = \mu - \si \; \f{
\mcal{Z}_{\de} }{\sq{n}}, \qqu \wh{G} (n, \mu, \de) = g (n, \mu,
\de) = \mu + \si \; \f{ \mcal{Z}_{\de} }{\sq{n}}
\]
for $n \in \bb{N}, \; \mu \in (- \iy, \iy), \; \de \in (0,
\f{1}{2})$. It can be shown that the estimator $\bs{\varphi}_n$
possesses all properties described at the beginning of Section
\ref{appdenmo} and consequently, all testing methods proposed in
previous sections can be used.

\subsubsection{Testing the Variance of a Normal Distribution}

Let $X$ be a Gaussian random variable with mean $\mu$ and variance
$\si^2$.  In many applications, it is important to test the variance
based on i.i.d. random samples $X_1, X_2, \cd$ of $X$.

In situations that the mean $\mu$ of the Gaussian variable $X$ is
known, we shall use {\small $\bs{\varphi}_n = \sq{ \f{1}{n}
\sum_{i=1}^{n} (X_i - \mu)^2}$} as an estimator of $\si$.  It can be
verified that {\small \bee & & \wh{F} (n, \si, \de) = \si \; \sq{
\f{\chi_{n, \de}^2}{n}}, \qu \qqu \qqu \qqu \qqu \qu \wh{G} (n,
\si, \de) = \si \; \sq{ \f{ \chi_{n, 1 - \de}^2}{n} }, \\
&   & f (n, \si, \de) = \si \; \min \li \{ 1, \; \sq{ \f{\chi_{n,
\de}^2}{n}} \ri \}, \qu \qqu \qqu  g (n, \si, \de) = \si \; \max \li
\{ 1, \; \sq{ \f{ \chi_{n, 1 - \de}^2}{n} } \ri \} \eee} for $n \in
\bb{N}, \; \si \in (0, \iy), \; \de \in (0, 1)$.  Moreover, it can
be verified that the estimator $\bs{\varphi}_n$ possesses all
properties described at the beginning of Section \ref{appdenmo}.
This implies that all testing methods proposed in previous sections
are applicable.

In situations that the mean $\mu$ of the Gaussian variable $X$ is
unknown, we shall use {\small $\bs{\varphi}_n = \sq{ \f{1}{n}
\sum_{i=1}^{n} (X_i - \ovl{X}_{n})^2}$}, where {\small $\ovl{X}_{n}
= \f{ \sum_{i = 1}^{n} X_i } { n }$}, as an estimator of $\si$.  To
design multistage sampling schemes  for testing $\si$, we shall make
use of the observation that $\bs{\varphi}_n$ is a ULCE of $\si$ and
relevant results described in previous sections.   By the definition
of $\bs{\varphi}_n$, it can be readily shown that {\small
\[
f (n, \si, \de) = \si \; \min \li \{ 1, \; \sq{ \f{\chi_{n - 1,
\de}^2}{n}} \ri \}, \qqu \qqu  g (n, \si, \de) = \si \; \max \li \{
1, \; \sq{ \f{ \chi_{n - 1, 1 - \de}^2}{n} } \ri \}
\]}
for $n \in \bb{N}, \; \si \in (0, \iy), \; \de \in (0, 1)$.  Let
$\bs{\al} = O(\ze) \in (0, 1), \; \bs{\ba} = O(\ze) \in (0, 1)$ and
$0 < \si^\prime < \si^{\prime \prime}$.  Let $\ovl{n} (\ze)$ be the
minimum integer $n$ such that $f(n, \si^{\prime \prime}, \bs{\ba} )
\geq g(n, \si^\prime, \bs{\al} )$.  We can show that \be \la{toshow}
\ovl{n} (\ze) \leq \max \li \{ \f{ 2 \ln \bs{\al}  } { 1 -
\f{\si^{\prime \prime}}{\si^\prime}  + \ln \f{\si^{\prime
\prime}}{\si^\prime} } + 1, \; \f{2  \ln \bs{\ba} } { 1 -
\f{\si^\prime}{\si^{\prime \prime}} + \ln \f{\si^\prime}{\si^{\prime
\prime}} } + 1, \; \f{1}{ 1 - \f{\si^\prime}{\si^{\prime \prime}} }
\ri \} = O \li ( \ln \f{1}{\ze} \ri ). \ee To show (\ref{toshow}),
note that $f(n, \si^{\prime \prime}, \bs{\ba} ) \geq g(n,
\si^\prime, \bs{\al} )$ is equivalent to \be \la{need} \max \{ n,
\chi_{n- 1, 1 - \bs{\al}}^2 \} \leq \li ( \f{\si^{\prime
\prime}}{\si^\prime}
 \ri )^2 \min \{ n, \chi_{n- 1, \bs{\ba}}^2 \}.
\ee Let $Z$ be a chi-square variable of $n - 1$ degrees of freedom.
Then, $\Pr \{  Z \geq \chi_{n- 1, 1 - \bs{\al}}^2 \} = \bs{\al}$ and
$\Pr \{ Z \leq \chi_{n- 1, \bs{\ba}}^2 \} = \bs{\ba}$.  By Lemma
\ref{lemCH} in Appendix \ref{App_Normal_Mean_unknown_variance}, we
have
\[
\Pr \li \{ Z \geq (n-1) \li ( \f{\si^{\prime \prime}}{\si^\prime}
 \ri ) \ri \} \leq \li [ \li ( \f{\si^{\prime \prime}}{\si^\prime}
 \ri ) \exp \li ( 1 - \f{\si^{\prime \prime}}{\si^\prime} \ri ) \ri ]^{(n - 1) \sh 2}
 \leq \bs{\al}
\]
and thus $\chi_{n- 1, 1 - \bs{\al}}^2 \leq (n-1) \li (
\f{\si^{\prime \prime}}{\si^\prime} \ri ) < n \li ( \f{\si^{\prime
\prime}}{\si^\prime} \ri )^2$ provided that {\small $\f{n - 1}{2}
\geq \f{ \ln \bs{\al}  } { 1 - \f{\si^{\prime \prime}}{\si^\prime} +
\ln \f{\si^{\prime \prime}}{\si^\prime} }$}.  Similarly, by Lemma
\ref{lemCH} in Appendix \ref{App_Normal_Mean_unknown_variance}, we
have
\[
\Pr \li \{ Z \leq (n-1) \li ( \f{\si^\prime}{\si^{\prime \prime}}
 \ri ) \ri \} \leq \li [ \li ( \f{\si^\prime}{\si^{\prime \prime}}
 \ri ) \exp \li ( 1 - \f{\si^\prime}{\si^{\prime \prime}} \ri ) \ri ]^{(n - 1) \sh 2}
 \leq \bs{\ba}
\]
and thus $\chi_{n- 1, \bs{\ba}}^2 \geq (n-1) \li (
\f{\si^\prime}{\si^{\prime \prime}}
 \ri ) > n \li ( \f{\si^\prime}{\si^{\prime \prime}}  \ri )^2$ provided that
\[
\f{n - 1}{2} \geq \f{ \ln \bs{\ba}  } { 1 -
\f{\si^\prime}{\si^{\prime \prime}}  + \ln
\f{\si^\prime}{\si^{\prime \prime}} }, \qqu n > \f{1}{ 1 -
\f{\si^\prime}{\si^{\prime \prime}} }.
\]
It can be seen that a sufficient condition for (\ref{need}) is
\[
n \geq \max \li \{ \f{ 2 \ln \bs{\al}  } { 1 - \f{\si^{\prime
\prime}}{\si^\prime}  + \ln \f{\si^{\prime \prime}}{\si^\prime} } +
1, \; \f{2  \ln \bs{\ba} } { 1 - \f{\si^\prime}{\si^{\prime \prime}}
+ \ln \f{\si^\prime}{\si^{\prime \prime}} } + 1, \; \f{1}{ 1 -
\f{\si^\prime}{\si^{\prime \prime}} } \ri \}.
\]
It follows immediately that (\ref{toshow}) is true.  Making use of
(\ref{toshow}), we can show that, in the context of testing multiple
hypotheses regarding $\si$ with our proposed multistage testing
plan, the risk of making wrong decisions can be made arbitrarily
small by choosing a sufficiently small $\ze > 0$.  Specifically, if
we identify parameter $\se$ in Theorem \ref{Multi_Comp_Exact} as
$\si$, using (\ref{toshow}), we can show that $\lim_{\ze \to 0} \Pr
\{ \tx{Reject} \; \mscr{H}_i \mid \se \} = 0$ for any $\se \in
\varTheta_i$ and $i = 0, 1, \cd, m - 1$.

Our method for the exact computation of the OC function
$\Pr \{ \tx{Accept} \; \mscr{H}_0 \mid \si \}$ is
described as follows. Since $\Pr \{ \tx{Accept} \; \mscr{H}_0 \mid \si \}
= 1 - \Pr \{ \tx{Reject} \; \mscr{H}_0
\mid \si \}$, it suffices to compute $\Pr
 \{ \tx{Reject} \; \mscr{H}_0 \mid \si \}$. By the definition
of the testing plan, we have \be \la{good289test} \Pr \li \{
\tx{Reject} \; \mscr{H}_0 \mid \si \ri \} = \sum_{\ell = 1}^s \Pr
\li \{ \bs{\varphi}_{n_\ell} > b_\ell, \; a_j \leq
\bs{\varphi}_{n_j} \leq b_j, \; 1 \leq j < \ell \mid \si \ri \}. \ee
If we choose the sample sizes to be even numbers $n_\ell = 2 k_\ell,
\; \ell = 1, \cd, s$ for the case of known variance and odd numbers
$n_\ell = 2 k_\ell + 1, \; \ell = 1, \cd, s$ for the case of unknown
variance, we can rewrite (\ref{good289test}) as {\small \be
\la{ssen28test} \Pr \li \{ \tx{Reject} \; \mscr{H}_0 \mid \si \ri \}
= \sum_{\ell = 1}^s \Pr \li \{ \sum_{q = 1}^{k_\ell} Z_q \geq
\f{n_\ell}{2} \li ( \f{ b_\ell }{\si} \ri )^2, \; \f{n_j}{2} \li (
\f{ a_j }{\si} \ri )^2 \leq \sum_{q = 1}^{k_j} Z_q \leq \f{n_j}{2}
\li ( \f{ b_j }{\si} \ri )^2 \; \tx{for} \; 1 \leq j < \ell \mid \si
\ri \}, \qu \ee} where $Z_1, Z_2, \cd$ are i.i.d. exponential random
variables with common mean unity. To compute the probabilities in
the right-hand side of (\ref{ssen28test}), we can make use of the
following results established by Chen \cite{Chen_EST}.

\beT \la{lemchen} Let $1 = k_0 < k_1 < k_2 < \cd$ be a sequence of
positive integers.  Let $0 = z_0 < z_1 < z_2 < \cd$ be a sequence of
positive numbers. Define $w(0, 1) = 1$ and
\[
w(\ell, 1) = 1, \qu w(\ell, q) = \sum_{i = 1}^{k_r} \f{ w(r, i) \;
(z_{\ell} - z_r)^{q - i} } { (q - i)!  }, \qu k_r < q \leq k_{r +
1}, \qu r = 0, 1, \cd, \ell - 1
\]
for $\ell = 1, 2, \cd$.  Let $Z_1, Z_2, \cd$ be i.i.d. exponential
random variables with common mean unity.
Then,
\[
\Pr \li \{\sum_{q = 1}^{k_j} Z_q > z_j \; \tx{for} \; j = 1, \cd,
\ell \ri \} =  e^{- z_\ell} \sum_{q = 1}^{k_\ell} w(\ell, q)
\]
for $\ell = 1, 2, \cd$.  Moreover, the following statements hold true.

(I)  {\small \bee &   & \Pr \li \{ a_j < \sum_{q = 1}^{k_j} Z_q <
b_j \; \tx{for} \; j = 1, \cd, \ell \ri \} \\
& = & \li [ \sum_{i = 1}^{2^{\ell - 1}} \Pr \li \{ \sum_{q =
1}^{k_j} Z_q
> [A_\ell]_{i, j} \; \tx{for} \; j = 1, \cd, \ell \ri \} \ri ] -
\li [ \sum_{i = 1}^{2^{\ell - 1}} \Pr \li \{ \sum_{q = 1}^{k_j} Z_q
> [B_\ell]_{i, j} \; \tx{for} \; j = 1, \cd, \ell \ri \} \ri ],
\eee} where $A_1 = [ a_1 ], \;  B_1 = [ b_1 ]$ and
\[ A_{r + 1} = \bem A_r & a_{r + 1} I_{2^{r - 1} \times 1}\\
B_r & b_{r + 1} I_{2^{r - 1} \times 1} \eem, \qqu
B_{r + 1} = \bem B_r & a_{r + 1} I_{2^{r - 1} \times 1}\\
A_r & b_{r + 1} I_{2^{r - 1} \times 1} \eem, \qqu r = 1, 2, \cd,
\] where $I_{2^{r - 1} \times 1}$ represents a column matrix with all $2^{r - 1}$
elements assuming value $1$.

(II)

{\small \bee &  & \Pr \li \{ a_j < \sum_{q = 1}^{k_j} Z_q < b_j \;
\tx{for} \; j = 1, \cd, \ell, \; \sum_{q = 1}^{k_{\ell + 1}} Z_q >
b_{\ell + 1} \ri
\}\\
 & = & \li [ \sum_{i = 1}^{2^{\ell-1}} \Pr \li \{ \sum_{q =
1}^{k_j} Z_q > [E]_{i, j} \; \tx{for} \; j = 1, \cd, \ell + 1 \ri \}
\ri ]  - \li [ \sum_{i = 1}^{2^{\ell-1}} \Pr \li \{ \sum_{q =
1}^{k_j} Z_q > [F]_{i, j} \; \tx{for} \; j = 1, \cd, \ell + 1 \ri \}
\ri ], \eee} where $E = \bem A_\ell & b_{\ell + 1} I_{2^{\ell - 1}
\times 1} \eem$ and $F = \bem B_\ell  & b_{\ell + 1} I_{2^{\ell - 1}
\times 1} \eem$.

(III) {\small \bee &  & \Pr \li \{ a_j < \sum_{q = 1}^{k_j} Z_q <
b_j \; \tx{for} \; j = 1, \cd, \ell, \; \sum_{q = 1}^{k_{\ell + 1}}
Z_q <
b_{\ell + 1} \ri \}\\
 & = & \Pr \li \{ a_j < \sum_{q = 1}^{k_j} Z_q < b_j
\; \tx{for} \; j = 1, \cd, \ell \ri \} - \Pr \li \{ a_j < \sum_{q =
1}^{k_j} Z_q < b_j \; \tx{for} \; j = 1, \cd, \ell, \; \sum_{q =
1}^{k_{\ell + 1}} Z_q > b_{\ell + 1} \ri \}. \eee}

 \eeT

\subsubsection{Testing the Parameter of an Exponential Distribution}

Let $X$ be a random variable with density function $f(x) =
\f{1}{\se} e^{- \f{x}{\se} }$ for $0 < x < \infty$, where $\se$ is a
parameter.  In many applications, it is important to test the
parameter $\se$ based on i.i.d. random samples $X_1, X_2, \cd$ of
$X$.  We shall use  {\small $\bs{\varphi}_n = \f{\sum_{i=1}^{n} X_i
}{n}$} as an estimator for $\se$.   Accordingly, for $\ell = 1, \cd,
s$, the estimator of $\se$ at the $\ell$-th stage is {\small
$\wh{\bs{\se}}_\ell = \bs{\varphi}_{n_\ell} = \f{\sum_{i=1}^{n_\ell}
X_i }{n_\ell}$}.  It can be shown that {\small \bee &   & \wh{F} (n,
\se, \de) = \f{\se \chi_{2 n, \de}^2}{2 n}, \qu \qqu \qqu \qqu \qqu
\wh{G} (n, \se,
\de) = \f{\se \chi_{2 n, 1 - \de}^2}{2 n},\\
&   &  f (n, \se, \de) = \se \; \min \li \{ 1, \; \f{\chi_{2 n,
\de}^2}{2 n} \ri \}, \qu  \qqu \qu g ( n, \se, \de) = \se \; \max
\li \{ 1, \; \f{ \chi_{2 n, 1 - \de}^2}{2 n} \ri \} \eee} for $n \in
\bb{N}, \; \se \in (0, \iy), \; \de \in (0, 1)$.   Since the
estimator $\bs{\varphi}_n$ possesses all properties described at the
beginning of Section \ref{appdenmo},  all testing methods proposed
in previous sections are applicable.  Moreover,  it is possible to
exactly  compute the OC function $\Pr \{ \tx{Accept} \; \mscr{H}_0
\mid \se \}$. Since $\Pr \{ \tx{Accept} \; \mscr{H}_0 \mid \se \} =
1 - \Pr \{ \tx{Reject} \; \mscr{H}_0 \mid \se \}$, it suffices to
compute $\Pr \li \{ \tx{Reject} \; \mscr{H}_0 \mid \se \ri \}$. By
the definition of the stopping rule, we have  \be
\la{good289exptest} \Pr \li \{ \tx{Reject} \; \mscr{H}_0 \mid \se
\ri \} = \sum_{\ell = 1}^s \Pr \li \{ \wh{\bs{\se}}_\ell > b_\ell,
\; a_j \leq \wh{\bs{\se}}_\ell \leq b_j, \; 1 \leq j < \ell \mid \se
\ri \}. \ee Let $Z_1, Z_2, \cd$ be i.i.d. exponential random
variables with common mean unity. Then, we can rewrite
(\ref{good289exptest}) as {\small \be \la{ssen28exptest} \Pr \li \{
\tx{Reject} \; \mscr{H}_0 \mid \se \ri \} = \sum_{\ell = 1}^s \Pr
\li \{ \sum_{q = 1}^{n_\ell} Z_q \geq n_\ell \li ( \f{ b_\ell }{\se}
\ri ), \; n_j \li ( \f{ a_j }{\se} \ri ) \leq \sum_{q = 1}^{n_j} Z_q
\leq n_j \li ( \f{ b_j }{\se} \ri ) \; \tx{for} \; 1 \leq j < \ell
\mid \se \ri \}. \ee} To evaluate the probabilities in the
right-hand side of (\ref{ssen28exptest}), we can make use of the
results in Theorem \ref{lemchen}.

\subsubsection{Testing the Scale Parameter of a Gamma Distribution}

In probability theory and statistics, a random variable $Y$ is said
to have a gamma distribution if its density function is of the form
\[
f(y) = \frac{y^{k - 1}} { \Gamma(k) \se ^{ k }    }  \exp \li ( -
\frac{y}{\se} \ri ) \;\;\; {\rm for} \;\;\; 0 < y < \infty
\]
where $\se > 0, \; k > 0$ are referred to as the scale parameter and
shape parameter respectively.  To test the scale parameter, $\se$,
of a Gamma distribution, consider random variable $X = \f{Y}{k}$.
Let $Y_1, Y_2, \cd$ be i.i.d. samples of $Y$ and $X_i = \f{Y_i}{k}$
for $i = 1, 2, \cd$. Define $\bs{\varphi}_n = \f{ \sum_{i = 1}^n X_i
}{n}$.  Then, $\bs{\varphi}_n$ is an unbiased and unimodal
likelihood estimator of $\se$ for all positive integer $n$.  It
follows that we can apply the theory and techniques in Section 2 to
test the multiple hypotheses like (\ref{mainpr}).

\subsubsection{Life Testing}

In this section, we shall consider the problem of life testing using
the classical exponential model \cite{Epstein}.  Suppose the lengths
of life of all components to be tested can be modeled as i.i.d.
random variables with common probability density function $f_T (t) =
\lm \exp \li ( - \lm t \ri )$,  where the parameter $\lm > 0$ is
referred to as the {\it failure rate} and its reciprocal is referred
to as the  {\it mean time between failures}. In reliability
engineering, it is a central issue to test the failure rate $\lm$
based on i.i.d. random samples $X_1, X_2, \cd$ of $X$.

In practice, for purpose of efficiency, multiple components are
initially placed on test. The test can be done with or without
replacement whenever a component fails. The decision of rejecting,
or accepting hypotheses or continuing test is based on the number of
failures and the {\it accumulated test time}. Here it should be
emphasized  that the accumulated test time is referred to as the
total running time of all components placed on test instead of the
real time.

The main idea of existing life-testing plans is to check how much test time has been accumulated whenever a failure occurs.   The test plans are
designed by truncating the sequential probability ratio tests (SPRT). There are several drawbacks with existing test plans. First, the existing
test plans are limited by the number of hypotheses. Currently, there is no highly effect methods for testing more than two hypotheses. Second,
when the indifference zone is narrow, the required accumulated test time may be very long. Third, the specified level of power may not be
satisfied due to the truncation of SPRT. Four, the administrative cost may be very high in the situations of high failure rate, since it
requires to check the status of test whenever a component fails. To overcome such drawbacks,  we tackle the life testing problem in the general
framework of testing $m$ mutually exclusive and exhaustive composite hypotheses: \[ \mscr{H}_0: 0 < \lm \leq \lm_1;  \qqu \mscr{H}_i: \lm_i <
\lm \leq \lm_{i+1},  \;\; i = 1, \cd, m - 2; \qqu \mscr{H}_{m - 1}:  \lm > \lm_{m - 1}, \] where $\lm_1 < \lm_2 < \cd < \lm_{m - 1}$. To control
the probabilities of making wrong decisions, it is typically required that, for pre-specified numbers $\de_i \in (0, 1), \; i = 0, 1, \cd, m-1$,
\bee &  & \Pr \{ \tx{Accept} \; \mscr{H}_0 \mid \lm \} \geq 1 - \de_0
\qqu \tx{for} \; 0 < \lm \leq \lm_1^{\prime}, \\
&  & \Pr \{ \tx{Accept} \; \mscr{H}_i \mid \lm \} \geq 1 - \de_i
\qqu \tx{for} \; \lm_i^{\prime \prime} \leq \lm \leq \lm_{i+1}^{\prime}  \; \tx{and} \;  i = 1, \cd, m - 2, \\
&  & \Pr \{ \tx{Accept} \; \mscr{H}_{m-1} \mid \lm \} \geq 1 -
\de_{m-1} \qqu \tx{for} \; \lm \geq \lm_{m - 1}^{\prime \prime}
 \eee
where $\lm_i^\prime, \; \lm_i^{\prime \prime}$ are parametric values such that $0 < \lm_1^\prime < \lm_1, \; \lm_{m-1}^{\prime \prime} >
\lm_{m-1}$ and $\lm_{i - 1} < \lm_{i - 1}^{\prime \prime} \leq \lm_i^\prime < \lm_i < \lm_i^{\prime \prime} \leq \lm_{i + 1}^\prime < \lm_{i +
1}$ for $i = 2, \cd, m - 2$.   This problem can be addressed by the general principle described in previous sections. Specifically, we proceed
as follows.

Let $\vDe$ be a positive number. Let $Z$ be the number of attempted
connections in a time interval of length $\vDe$. Then, $Z$ is a
Poisson variable of mean value $\lm \vDe$. Define $X = \f{Z}{\vDe}$.
The distribution of $X$ is determined as
\[
\Pr \li \{  X = \f{k}{\vDe} \ri \} = \f{ (\lm \vDe)^k e^{-\lm \vDe }
} { k! }, \qqu k = 0,1,2, \cd.
\]
Let $X_i = \f{Z_i}{\vDe}$, where $Z_i$ is the number of attempts in
time interval $[(i-1) \vDe, \; i \vDe)$ for $i = 1, 2, \cd$.  It
follows that $X_1, X_2, \cd$ are i.i.d. samples of $X$. Therefore,
the life testing problem can be cast in our general framework of
multistage hypothesis tests with sample sizes $n_1 < n_2 < \cd <
n_s$. Accordingly, the testing time is $t_\ell = n_\ell \vDe, \;
\ell = 1, \cd, s$.  For $\ell = 1, \cd, s$,  we propose to define
the estimator for $\lm$ at the $\ell$-th stage  as
\[
\wh{\bs{\lm}}_\ell = \varphi(X_1, \cd, X_{n_\ell}) = \f{
\sum_{i=1}^{n_\ell} X_i }{ n_\ell  }  = \f{ \sum_{i=1}^{n_\ell} Z_i
}{ n_\ell \vDe } = \f{ \tx{Number of arrivals in} \; [0, t_\ell) }{
t_\ell }.
\]
Clearly, $\bs{\varphi}_n$ is a ULCE of $\lm$; $\bs{\varphi}_n$ is an
unbiased estimator of $\lm$; the likelihood ratio is monotonically
increasing with respect to $\bs{\varphi}_n$. Hence, the estimator
$\bs{\varphi}_n$ possesses all the properties described at the
beginning of Section \ref{appdenmo}. This implies that all testing
methods proposed in previous sections are applicable.

It can be seen that all tests described above depend on, $\vDe$, the
unit of time used to convert the continuous time process to a
discrete time process. In applications, it may be preferred to use
the test derived by letting $\vDe \to 0$.   In this direction, we
have established such limiting procedure as follows. The testing
process is divided into $s$ stages with testing time
 $t_1 < t_2 < \cd < t_s$. For $\ell = 1, \cd, s$, at the $\ell$-th
stage, a decision variable $\bs{D}_\ell$ is defined based on
estimator $\wh{\bs{\lm}}_\ell = \f{ \tx{Number of arrivals in} \;
[0, t_\ell) }{ t_\ell }$ for $\lm$ such that the sampling process is
continued if $\bs{D}_\ell = 0$ and that hypothesis $\mscr{H}_i$ is
accepted if $\bs{D}_\ell = i + 1$, where $i \in \{ 0, 1, \cd, m - 1
\}$. Define multivariate functions {\small \bee & & f (t, \lm, \de)
= \bec \f{1}{t} \times \max \li \{ k \in \bb{Z}: \sum_{i = 0}^k \f{
(t \lm)^i e^{-t \lm} }{i!} \leq \de, \; 0 \leq k \leq t \lm \ri \} &
\tx{for} \;
t \geq \f{ \ln (\de) } { - \lm },\\
- \iy & \tx{for} \; t < \f{ \ln (\de) } { - \lm } \eec\\
&  & g (t, \lm, \de) = \f{1}{t} \times \min \li \{ k \in \bb{Z}:
 \sum_{i = 0}^{k-1} \f{ (t \lm)^i e^{- t \lm} }{i!} \geq 1 - \de, \; k \geq t \lm \ri \} \eee}
for $t > 0, \lm > 0, \de \in (0, 1)$ and multivariate functions
{\small \bee & & \udl{f} (t, \lm^\prime, \lm^{\prime \prime},
\de^\prime, \de^{\prime \prime} ) = \min \li \{ f (t, \lm^{\prime
\prime}, \de^{\prime \prime}), \; \f{1}{2} [ f (t, \lm^{\prime
\prime}, \de^{\prime \prime}) + g (t, \lm^{\prime}, \de^{\prime}) ] \ri \}, \\
 &  & \ovl{g} (t, \lm^\prime, \lm^{\prime \prime},
\de^\prime, \de^{\prime \prime} ) = \max \li \{ g (t, \lm^{\prime},
\de^{\prime}), \; \f{1}{2} [ f (t, \lm^{\prime \prime}, \de^{\prime
\prime}) + g (t, \lm^{\prime}, \de^{\prime}) ] \ri \}  \eee} for $0
< \lm^\prime < \lm^{\prime \prime}$ and $\de^\prime, \de^{\prime
\prime} \in (0, 1)$.  Let $\al_i = O(\ze) \in (0, 1), \ba_i = O(\ze)
\in (0, 1)$ for $i = 1, \cd, m - 1$ and $\al_m = \ba_0 = 0$.  Under
the assumptions  that  the maximum testing time $t_s$ is no less
than  the minimum positive number $t$ such that $f (t, \lm_i^{\prime
\prime}, \ba_i) \geq g (t, \lm_i^{\prime}, \al_i)$ for $i = 1, \cd,
m-1$,  We  propose to define the decision variables as \be
\la{defgoog}
 \bs{D}_\ell = \bec 1 & \tx{if} \; \;
\wh{\bs{\lm}}_\ell \leq
f_{\ell ,1}, \\
i & \tx{if} \; \;  g_{\ell, i-1} < \wh{\bs{\lm}}_\ell \leq
f_{\ell,i}
\; \tx{where} \; 2 \leq i \leq m - 1,\\
m  & \tx{if} \; \;  \wh{\bs{\lm}}_\ell > g_{\ell ,m-1},\\
0 & \tx{else} \eec \ee for $\ell = 1, \cd, s$, where $f_{\ell, i} =
\udl{f} (t_\ell, \lm_i^{\prime}, \lm_i^{\prime \prime}, \al_i,
\ba_i)$ and $g_{\ell, i} = \ovl{g} (t_\ell, \lm_i^{\prime},
\lm_i^{\prime \prime}, \al_i, \ba_i)$ for $i = 1, \cd, m-1$.

In order to simply the stopping boundary of the testing plans,
define multivariate functions {\small \bee &  &  f_c (t, \lm, \de) =
\bec \max \{ z \in (0, \lm]:  z - \lm + z \ln \f{\lm}{z}  \leq \f{
\ln (\de) }{t} \} & \tx{for}
\;  t \geq \f{ \ln (\de) } { - \lm },\\
- \iy & \tx{for} \;  t < \f{ \ln (\de) } { - \lm } \eec\\
& & g_c (t, \lm, \de) = \min \li \{ z \in [\lm, \iy): z - \lm + z
\ln \f{\lm}{z} \leq \f{ \ln (\de) }{t}  \ri \} \eee}  for $t > 0,
 \lm > 0, \de \in (0, 1)$.  Moreover, define {\small \bee &
& \udl{f}_c (t, \lm^\prime, \lm^{\prime \prime}, \de^\prime,
\de^{\prime \prime} ) = \min \li \{  f_c (t, \lm^{\prime \prime},
\de^{\prime \prime}), \; \f{1}{2} [ f_c (t, \lm^{\prime \prime},
\de^{\prime \prime}) + g_c (t, \lm^{\prime}, \de^{\prime}) ] \ri \}, \\
&  & \ovl{g}_c (t, \lm^\prime, \lm^{\prime \prime}, \de^\prime,
\de^{\prime \prime} ) = \max \li \{ g_c (t, \lm^{\prime},
\de^{\prime}), \; \f{1}{2} [ f_c (t, \lm^{\prime \prime},
\de^{\prime \prime}) + g_c (t, \lm^{\prime}, \de^{\prime}) ] \ri \}
\eee} for $0 < \lm^\prime < \lm^{\prime \prime}; \; \de^\prime,
\de^{\prime \prime} \in (0, 1)$ and $t
> 0$.  Under the assumption that the maximum testing time
$t_s$ is no less than the minimum positive number $t$ such that $f_c
(t, \lm_i^{\prime \prime}, \ba_i) \geq g_c (t, \lm_i^{\prime},
\al_i)$ for $i = 1, \cd, m-1$, we propose to define decision
variable $\bs{D}_\ell$ by (\ref{defgoog}) for $\ell = 1, \cd, s$
with $f_{\ell, i} = \udl{f}_c (t_\ell, \lm_i^{\prime}, \lm_i^{\prime
\prime}, \al_i, \ba_i)$ and $g_{\ell, i} = \ovl{g}_c (t_\ell,
\lm_i^{\prime}, \lm_i^{\prime \prime}, \al_i, \ba_i)$ for $i = 1,
\cd, m-1$.   We have established that the same conclusion as
described by statements (I)--(IX) of Theorem 1 holds true.

Clearly, once the limits of testing time are determined, we have a
multistage test plan which depends on the risk tuning parameter
$\ze$. We can evaluate the risk of such a limiting test plan. If the
risk requirement is not satisfied, then we can change $\ze$ and find
the corresponding limiting test plan.  This process can be repeated
until a satisfactory test plan is found.

In this section, we only consider the general problem of testing multiple hypotheses. The approach can be readily adapted to special problems
such as testing one-sided hypotheses, two-sided hypotheses, triple hypotheses, and interval test, etc. Specific procedures can be developed by
mimicking the techniques presented in Sections 3.

\section{Tests for the Mean of a Normal Distribution with Unknown Variance}

In this section, we shall focus on tests for the mean, $\mu$, of a
Gaussian variable $X$ with unknown variance $\si^2$ based on i.i.d.
samples $X_1, X_2, \cd$ of $X$. Our objective is to develop
multistage sampling schemes for testing hypotheses regarding $\se =
\f{\mu}{\si}$, which is the ratio of the mean to the standard
deviation.

\subsection{General Principle}

A general problem regarding $\se = \f{\mu}{\si}$ is to test  $m$
mutually exclusive and exhaustive composite hypotheses: \[
\mscr{H}_0: \se \in {\Se}_0, \qu \mscr{H}_1: \se \in {\Se}_1, \qu
\ldots, \qu \mscr{H}_{m - 1}: \se \in {\Se}_{m - 1}, \] where $\Se_0
= ( - \iy, \se_1], \; \Se_{m-1} = (\se_{m - 1}, \iy )$ and $\Se_i =
(\se_i, \se_{i+1} ], \; i = 1, \cd, m - 2$ with $\se_1 < \se_2 < \cd
< \se_{m - 1}$. To control the probabilities of making wrong
decisions, it is typically required that, for pre-specified numbers
$\de_i \in (0, 1)$, \bee &   & \Pr \li \{ \tx{Accept} \; \mscr{H}_i
\mid \se \ri \} \geq 1 - \de_i  \qu \fa \se \in \varTheta_i , \qqu i
= 0, 1, \cd, m - 1 \eee with $\varTheta_0 = ( - \iy,
\se_1^{\prime}], \; \varTheta_{m-1} = [ \se_{m - 1}^{\prime \prime},
\iy)$ and $\varTheta_i = [ \se_i^{\prime \prime}, \se_{i+1}^{\prime}
]$ for $i = 1, \cd, m - 2$, where $\se_i^\prime, \; \se_i^{\prime
\prime}$ satisfy $\se_1^\prime < \se_1, \; \se_{m-1}^{\prime \prime}
> \se_{m-1}$ and $\se_{i - 1} < \se_{i - 1}^{\prime \prime} \leq
\se_i^\prime < \se_i < \se_i^{\prime \prime} \leq \se_{i + 1}^\prime < \se_{i + 1}$ for $i = 2, \cd, m - 2$.

\beT \la{Normal_Mean_unknown_variance} Suppose that $\al_i = O(\ze)
\in (0, 1)$ and $\ba_i = O(\ze) \in (0, 1)$ for $i = 1, \cd, m - 1$.
 Let $2 \leq n_1 < n_2 < \cd <
n_s$ be the sample sizes such that the largest sample size $n_s$ is
no less than the minimum integer $n$ guaranteeing $(\se_i^{\prime
\prime} - \se_i^{\prime}) \sq{n - 1} \geq t_{n - 1, \al_i} + t_{n -
1, \ba_i}$ for $i = 1, \cd, m - 1$. Define
\[
f_{\ell, i} = \min \li \{ \se_{i}^{\prime \prime} - \f{ t_{n_\ell -
1, \ba_i} } { \sq{n_\ell - 1} }, \qu  \f{ \se_{i}^{\prime} +
\se_{i}^{\prime \prime} }{2} + \f{ t_{n_\ell - 1,   \al_i} -
t_{n_\ell - 1,   \ba_i} } {2 \sq{n_\ell - 1} } \ri \} ,
\]
\[
g_{\ell, i} = \max \li \{ \se_{i}^{\prime} + \f{ t_{n_\ell - 1,
\al_i} } { \sq{n_\ell - 1} }, \qu \f{ \se_{i}^{\prime} +
\se_{i}^{\prime \prime} }{2} + \f{ t_{n_\ell - 1,   \al_i} -
t_{n_\ell - 1,   \ba_i} } {2 \sq{n_\ell - 1} } \ri \}
\]
for $i = 1, \cd, m - 1$.  Define {\small \[ \ovl{X}_{n_\ell} = \f{
\sum_{i = 1}^{n_\ell} X_i } { n_\ell }, \qqu \wt{\si}_{n_\ell} =
\sq{ \f{\sum_{i=1}^{n_\ell} (X_i - \ovl{X}_{n_\ell})^2}{n_\ell} },
\qqu \wh{\bs{\se}}_\ell = \f{ \ovl{X}_{n_\ell}}{\wt{\si}_{n_\ell}},
\]
\be \la{defgoognormal}
\bs{D}_\ell = \bec 1 & \tx{for} \; \wh{\bs{\se}}_\ell \leq f_{\ell,1},\\
i & \tx{for} \; g_{\ell, i - 1} < \wh{\bs{\se}}_\ell
\leq f_{\ell, i} \; \tx{where} \; 2 \leq i \leq m - 1,\\
m & \tx{for} \; \wh{\bs{\se}}_\ell > g_{\ell,m-1},\\
0 & \tx{else} \eec \ee} for $\ell = 1, \cd, s$.  Then,  $\Pr \{
\tx{Reject} \; \mscr{H}_i  \mid \se \}$ is no greater than $\de_i$
for any $\se \in \varTheta_i$ and $i = 0, 1, \cd, m - 1$ provided
that $\ze$ is sufficiently small.

\eeT

See Appendix \ref{App_Normal_Mean_unknown_variance} for a proof.   It should be noted that the probabilities like $\Pr \{ \tx{Reject} \;
\mscr{H}_i \mid \se \}$ can be evaluated by a Monte Carlo method.  Making use of the Monte Carlo method for estimating risks, the idea of
bisection risk tuning described in Section \ref{BRT}  and the iterative minimax optimization algorithm proposed in Section \ref{minmaxcom}, we
can develop efficient testing plans satisfying the risk requirement by determining appropriate risk tuning parameter $\ze$ and weighting
coefficients.

\subsection{Applications}

In this section, we shall study the applications of Theorem
\ref{Normal_Mean_unknown_variance} to specific testing problems.
Specially, the following Sections \ref{secA}, \ref{secB} and
\ref{secC} are devoted to the discussion of hypotheses concerned
with the comparison of the mean $\mu$ of Gaussian variable $X$ with
a pre-specified number $\ga$. Such issues can be formulated as
problems of testing hypotheses regarding $\vse = \f{ \mu -
\ga}{\si}$. To develop concrete testing plans, we make use of the
following statistics
\[
\ovl{X}_{n_\ell} = \f{ \sum_{i = 1}^{n_\ell} X_i } { n_\ell }, \qu
\wh{\si}_{n_\ell} = \sq{ \f{\sum_{i=1}^{n_\ell} (X_i -
\ovl{X}_{n_\ell})^2}{n_\ell - 1} }, \qu  \wh{T}_\ell = \f{
\sq{n_\ell} (\ovl{X}_{n_\ell} - \ga)}{\wh{\si}_{n_\ell}}
\]
for $\ell = 1, \cd, s$.

\subsubsection{One-sided Tests} \la{secA}

In many situations, it is an important problem to test hypotheses
$\mscr{H}_0 : \vse < 0$ versus $\mscr{H}_1 : \vse > 0$.  To control
the risks of committing decision errors, it is typically required
that, for prescribed numbers $\al, \; \ba \in (0, 1)$, \bee & & \Pr
\li \{ \tx{Accept} \; \mscr{H}_0 \mid \vse \ri \} > 1 -
\al \qu \tx{for  $\vse \leq - \vep$,}\\
&  & \Pr \li \{ \tx{Accept} \; \mscr{H}_1 \mid \vse \ri \} > 1 - \ba
\qu \tx{for  $\vse \geq \vep$,} \eee where the indifference zone is
$(- \vep, \vep)$.  Applying Theorem
\ref{Normal_Mean_unknown_variance} to the special case of $m = 2$,
we have the following results.

\beC \la{Normal_Mean_unknown_variance_One} Let $\bs{\al} = O(\ze)
\in (0, 1)$ and $\bs{\ba} = O(\ze) \in (0, 1)$.  Let $2 \leq n_1 <
n_2 < \cd < n_s$ be the sample sizes such that the largest sample
size $n_s$ is no less than  the minimum integer $n$ guaranteeing
{\small $t_{n - 1, \bs{\al}} + t_{n - 1, \bs{\ba}} \leq 2 \vep \sq{n
- 1}$}. Define $a_\ell = \vep \sq{n_\ell - 1} - t_{n_\ell - 1,
\bs{\ba}}, \; b_\ell = t_{n_\ell - 1, \bs{\al}} - \vep \sq{n_\ell -
1}$ for $\ell = 1, \cd, s - 1$, and {\small $a_s = b_s = \f{ t_{n_s
- 1, \bs{\al}} - t_{n_s - 1, \bs{\ba}} }{2}$}. Define {\small
\[  \bs{D}_\ell = \bec 1 &
\tx{for} \; \wh{T}_\ell \leq a_\ell,\\
2 & \tx{for} \; \wh{T}_\ell > b_\ell,\\
 0 & \tx{else} \eec
\]}
for $\ell = 1, \cd, s$.  Then,  the following statements hold true.

(i) $\Pr \li \{ \tx{Accept} \; \mscr{H}_0 \mid \vse \ri \}$ is less
than $\ba$ for  $\vse$ no less than $\vep$ if $\ze
> 0$ is sufficiently small.

(ii) $\Pr \li \{ \tx{Reject} \; \mscr{H}_0 \mid \vse \ri \}$ is less
than $\al$ for $\vse$ no greater than $- \vep$ if $\ze
> 0$ is sufficiently small.

\eeC

To develop an efficient testing plan satisfying the risk requirement, we need to determine appropriate risk tuning parameter $\ze$ and weighting
coefficients. This can be accomplished by making use of the Monte Carlo method for estimating risks, the idea of bisection risk tuning described
in Section \ref{BRT} and the iterative minimax optimization algorithm proposed in Section \ref{secones}.

\subsubsection{Two-sided Tests} \la{secB}

It is a frequent problem to test hypotheses $\mscr{H}_0 : \vse = 0$
versus $\mscr{H}_1 : \vse \neq 0$.  To control the risks of
committing decision errors, it is typically required that, for
prescribed numbers $\al, \; \ba \in (0, 1)$, \bee &  & \Pr \li \{
\tx{Accept} \; \mscr{H}_0 \mid \vse \ri \} > 1 -
\al \qu \tx{for  $\vse = 0$,}\\
&  & \Pr \li \{ \tx{Accept} \; \mscr{H}_1 \mid \vse \ri \} > 1 - \ba
\qu \tx{for  $|\vse| \geq \vep$}, \eee where the indifference zone
is $(- \vep, 0) \cup (0, \vep)$.  Applying Theorem
\ref{Normal_Mean_unknown_variance} to test hypotheses $\mcal{H}_0:
\vse \leq -\f{\vep}{2}, \; \mcal{H}_1: - \f{\vep}{2} < \vse \leq
\f{\vep}{2}$ and $\mcal{H}_2: \vse > \f{\vep}{2}$ with indifference
zone $(- \vep, 0) \cup (0, \vep)$, we have $\Pr \{ \tx{Reject} \;
\mcal{H}_0 \; \tx{and} \; \mcal{H}_2 \mid \vse \} = \Pr \{
\tx{Accept} \; \mscr{H}_0  \mid \vse \}$ and the following results
follow immediately.

\beC \la{Normal_Mean_unknown_variance_T} Let $\bs{\al} = O(\ze) \in
(0, 1)$ and $\bs{\ba} = O( \ze ) \in (0, 1)$.  Let $2 \leq n_1 < n_2
< \cd < n_s$ be the sample sizes such that the largest sample size
$n_s$ is no less than the minimum integer $n$ guaranteeing {\small
$t_{n - 1, \bs{\al}  }  + t_{n - 1, \bs{\ba}  } \leq \vep \sq{n -
1}$}. Define $a_\ell = \vep \sq{n_\ell - 1} - t_{n_\ell - 1,
\bs{\ba}  }, \; b_\ell = t_{n_\ell - 1, \bs{\al}  }$ for $\ell = 1,
\cd, s - 1$, and {\small $a_s = b_s = \f{ t_{n_s - 1, \bs{\al}  } -
t_{n_s - 1, \bs{\ba}  } }{2} + \f{\vep}{2} \sq{n_s - 1}$}. Define
{\small
\[ \bs{D}_\ell = \bec 1 &
\tx{for} \; | \wh{T}_\ell | \leq a_\ell,\\
2 & \tx{for} \; | \wh{T}_\ell | > b_\ell,\\
 0 & \tx{else} \eec
\]}
for $\ell = 1, \cd, s$.  Then,  the following statements hold true.

(i) $\Pr \li \{ \tx{Accept} \; \mscr{H}_0 \mid \vse \ri \}$ is less
than $\ba$ for any $\vse \in (-\iy, - \vep] \cup [ \vep, \iy)$ if
$\ze > 0$ is sufficiently small.

(ii) $\Pr \li \{ \tx{Reject} \; \mscr{H}_0 \mid \vse \ri \}$ is less
than $\al$ for $\vse = 0$ if $\ze > 0$ is sufficiently small.

\eeC

By virtue of the Monte Carlo method for estimating risks, the idea of bisection risk tuning described in Section \ref{BRT}  and the iterative
minimax optimization algorithm proposed in Section \ref{sec Twos}, we can determine appropriate risk tuning parameter $\ze$ and weighting
coefficients such that the risk requirement is satisfied.

\subsubsection{Tests of Triple Hypotheses} \la{secC}

In many applications, it is desirable to test three hypotheses
$\mscr{H}_0 : \vse < 0, \; \mscr{H}_1 : \vse = 0, \; \mscr{H}_2 :
\vse > 0$.  To control the risks of committing decision errors, it
is typically required that, for prescribed numbers $\al, \; \ba \in
(0, 1)$, \bee & & \Pr \li \{ \tx{Accept} \; \mscr{H}_0 \mid \vse \ri
\}
> 1 -
\ba \qu \tx{for  $\vse \leq - \vep$,}\\
&  & \Pr \li \{ \tx{Accept} \;
\mscr{H}_1 \mid \vse \ri \} > 1 - \al \qu \tx{for  $\vse = 0$},\\
&  & \Pr \li \{ \tx{Accept} \; \mscr{H}_2 \mid \vse \ri \} > 1 - \ba
\qu \tx{for  $\vse \geq \vep$,} \eee where the indifference zone is
$(- \vep, 0) \cup (0, \vep)$.  Applying Theorem
\ref{Normal_Mean_unknown_variance} to test hypotheses $\mcal{H}_0:
\vse \leq -\f{\vep}{2}, \; \mcal{H}_1: - \f{\vep}{2} < \vse \leq
\f{\vep}{2}$ and $\mcal{H}_2: \vse > \f{\vep}{2}$ with indifference
zone $(- \vep, 0) \cup (0, \vep)$, we have the following results.

\beC \la{Normal_Mean_unknown_variance_Trip} Let $\bs{\al} = O(\ze)
\in (0, 1)$ and $\bs{\ba} = O( \ze ) \in (0, 1)$.  Let $2 \leq n_1 <
n_2 < \cd < n_s$ be the sample sizes such that the largest sample
size $n_s$ is no less than the minimum integer $n$ guaranteeing
 {\small $t_{n - 1, \bs{\al}  } + t_{n - 1, \bs{\ba}  } \leq \vep
\sq{n - 1}$}.  Define $a_\ell = \vep \sq{n_\ell - 1} - t_{n_\ell -
1, \bs{\ba}  }, \; b_\ell = t_{n_\ell - 1, \bs{\al}  }$ for $\ell =
1, \cd, s - 1$, and {\small $a_s = b_s = \f{ t_{n_s - 1, \bs{\al} }
- t_{n_s - 1, \bs{\ba}  } }{2} + \f{\vep}{2} \sq{n_s - 1}$}. Define
{\small
\[
\bs{D}_\ell = \bec 1 & \tx{for} \;  \wh{T}_\ell  < - b_\ell,\\
2 & \tx{for} \; | \wh{T}_\ell | \leq a_\ell,\\
3 & \tx{for} \;  \wh{T}_\ell  > b_\ell,\\
 0 & \tx{else} \eec
\]}
for $\ell = 1, \cd, s$.  Then,  the following statements hold true.

(i) $\Pr \li \{ \tx{Accept} \; \mscr{H}_0 \mid \vse \ri \}$ is
greater than $1 - \ba$ for any $\vse \in (-\iy, - \vep ]$ if $\ze
> 0$ is sufficiently small.

(ii) $\Pr \li \{ \tx{Accept} \; \mscr{H}_2 \mid \vse \ri \}$ is
greater than $1 - \ba$ for any $\vse \in [ \vep , \iy)$ if $\ze
> 0$ is sufficiently small.

(iii) $\Pr \li \{ \tx{Accept} \; \mscr{H}_1 \mid \vse \ri \}$  for
$\vse = 0$ is greater than $1 - \al$ if $\ze > 0$ is sufficiently
small.

\eeC

Applying the Monte Carlo method for estimating risks, the idea of bisection risk tuning described in Section \ref{BRT}  and the iterative
minimax optimization algorithm proposed in Section \ref{secTriple}, we can determine appropriate risk tuning parameter $\ze$ and weighting
coefficients such that the risk requirement is satisfied.

\subsubsection{Interval Tests}

In some situations, it is desirable to test hypothesis $\mscr{H}_0:
\se \in [\se_1, \se_2]$ versus $\mscr{H}_1: \se \notin [\se_1,
\se_2]$. For risk control purpose, it is typically required that,
for two prescribed numbers $\al, \; \ba \in (0, 1)$, \bee &  & \Pr
\li \{ \tx{Reject} \; \mscr{H}_0 \mid \se \ri \} \leq \al \qu
\tx{for $\se \in \Se$ such that $\se \in [\se_1^{\prime
\prime}, \se_2^\prime]$}, \\
&  & \Pr \li \{ \tx{Accept} \; \mscr{H}_0 \mid \se \ri \} \leq \ba
\qu \tx{for  $\se \in \Se$ such that $\se \notin (\se_1^\prime,
\se_2^{\prime \prime})$ },  \eee where $\se_1^\prime < \se_1 <
\se_1^{\prime \prime} < \se_2^\prime < \se_2 < \se_2^{\prime
\prime}$. Since there is no requirement imposed on probabilities of
committing decision errors for $\se \in (\se_1^\prime, \se_1^{\prime
\prime}) \cup (\se_2^\prime, \se_2^{\prime \prime})$, the union of
intervals,  $(\se_1^\prime, \se_1^{\prime \prime}) \cup
(\se_2^\prime, \se_2^{\prime \prime})$, is referred to as an
indifference zone.

Applying Theorem \ref{Normal_Mean_unknown_variance} to test
hypotheses $\mcal{H}_0: \se \leq \se_1, \; \mcal{H}_1: \se_1 < \se
\leq \se_2$ and $\mcal{H}_2: \se
> \se_2$ with indifference zone $(\se_1^\prime, \se_1^{\prime \prime} )
 \cup (\se_2^\prime, \se_2^{\prime \prime}
)$, we have $\Pr \{ \tx{Reject} \; \mcal{H}_0 \; \tx{and} \;
\mcal{H}_2 \mid \se \} = \Pr \{ \tx{Accept} \; \mscr{H}_0  \mid \se
\}$ and the following results.

\beC  Let $\al_i = O(\ze) \in (0, 1)$ and $\ba_i = O(\ze) \in (0,
1)$ for $i = 1, 2$.  Let $2 \leq n_1 < n_2 < \cd < n_s$ be the
sample sizes such that the largest sample size $n_s$ is no less than
the minimum integer $n$ guaranteeing  $(\se_i^{\prime \prime} -
\se_i^{\prime}) \sq{n - 1} \geq t_{n - 1, \al_i}  + t_{n - 1,
\ba_i}$ for $i = 1, 2$.   Define
\[
f_{\ell, i} = \min \li \{  \se_i^{\prime \prime} - \f{ t_{n_\ell -
1, \ba_i} } { \sq{n_\ell - 1} }, \qu \f{\se_i^{\prime \prime} +
\se_i^{\prime} }{2} + \f{ t_{n_\ell - 1, \al_i} - t_{n_\ell - 1,
\ba_i} } {2 \sq{n_\ell - 1} } \ri \},
\]
\[
g_{\ell, i} = \max \li \{ \se_i^{\prime} + \f{ t_{n_\ell - 1, \al_i}
} { \sq{n_\ell - 1} }, \qu \f{\se_i^{\prime \prime} + \se_i^{\prime}
}{2} + \f{ t_{n_\ell - 1, \al_i} - t_{n_\ell - 1, \ba_i} } {2
\sq{n_\ell - 1} } \ri \}
\]
for $i = 1, 2$.  Define
\[ \bs{D}_\ell = \bec 1 & \tx{if} \; g_{\ell, 1} <
\wh{\bs{\se}}_\ell \leq f_{\ell, 2},\\
2 & \tx{if} \; \wh{\bs{\se}}_\ell \leq f_{\ell, 1} \; \tx{or} \;
\wh{\bs{\se}}_\ell > g_{\ell, 2},\\
0 & \tx{else} \eec \]  for $\ell = 1, \cd, s$.  Then, the following
statements hold true.

(i) $\Pr \{ \tx{Accept} \; \mscr{H}_0  \mid \se \} \leq \ba$ for
$\se \notin (\se_1^\prime, \se_2^{\prime \prime})$ if $\ze$ is
sufficiently small.

(ii) $\Pr \{ \tx{Reject} \; \mscr{H}_0 \mid \se \} \leq  \al$ for
 $\se \in [\se_1^{\prime \prime},
\se_2^\prime]$ if $\ze$ is sufficiently small.

\eeC

Making use of the Monte Carlo method for estimating risks, the idea of bisection risk tuning described in Section \ref{BRT}  and the iterative
minimax optimization algorithm proposed in Section \ref{secInterval Tests}, we can determine appropriate risk tuning parameter $\ze$ and
weighting coefficients such that the risk requirement is satisfied.

\subsubsection{Tests of ``Simple'' Hypotheses}

In some situations, it may be interesting to test multiple simple
hypotheses $\mscr{H}_i : \se = \se_i$ for $i = 0, 1, \cd, m - 1$.
For risk control purpose, it is typically required that, for
prescribed numbers $\de_i \in (0, 1)$,
\[
\Pr \li \{ \tx{Accept} \; \mscr{H}_i \mid \se_i \ri \} \geq 1 -
\de_i, \qqu i = 0, 1, \cd, m - 1.
\]
Applying Theorem \ref{Normal_Mean_unknown_variance} to test the
following hypotheses
\[
\mcal{H}_0: \se \leq \vse_1, \qu \mcal{H}_1: \vse_1 < \se \leq
\vse_2, \qu \ldots, \qu \mcal{H}_{m-2}: \vse_{m-2} < \se \leq
\vse_{m-1}, \qu \mcal{H}_{m-1}: \se > \vse_{m-1}
\]
with $\vse_i = \f{ \se_{i - 1} + \se_i} {2}, \; i = 1, \cd, m-1$ and
indifference zone $\cup_{i=1}^{m-1} (\se_{i-1}, \se_i)$, we have the
following results.

\beC Let $\al_i = O(\ze) \in (0, 1)$ and $\ba_i = O(\ze) \in (0, 1)$
for $i = 1, \cd, m - 1$. Let $2 \leq n_1 < n_2 < \cd < n_s$ be the
sample sizes such that the largest sample size $n_s$ is no less than
the minimum integer $n$ guaranteeing $(\se_i - \se_{i-1}) \sq{n - 1}
\geq t_{n - 1, \al_i} + t_{n - 1, \ba_i }$ for $i = 1, \cd, m - 1$.
Define \bee &  & f_{\ell, i} = \min \li \{  \se_{i} - \f{ t_{n_\ell
- 1, {\ba_i}} } { \sq{n_\ell - 1} }, \qu \f{\se_{i} + \se_{i-1} }{2}
+ \f{ t_{n_\ell - 1, {\al_i}} - t_{n_\ell - 1, {\ba_i}} } {2
\sq{n_\ell - 1} } \ri \},\\
&  & g_{\ell, i} = \max \li \{ \se_{i-1} + \f{ t_{n_\ell - 1,
{\al_i}} } { \sq{n_\ell - 1} }, \qu \f{\se_{i} + \se_{i-1} }{2} +
\f{ t_{n_\ell - 1, {\al_i}} - t_{n_\ell - 1, {\ba_i}} } {2
\sq{n_\ell - 1} } \ri \} \eee for $i = 1, \cd, m - 1$. Define
$f_{\ell, i} = \udl{f} (n_\ell, \se_{i-1}, \se_i, \al_i, \ba_i)$ and
$g_{\ell, i} = \ovl{g} (n_\ell, \se_{i-1}, \se_i, \al_i, \ba_i)$ for
$i = 1, \cd, m-1$. Define
 decision variable $\bs{D}_\ell$ by (\ref{defgoognormal}) for $\ell = 1, \cd, s$.
Then, $\Pr \{ \tx{Reject} \; \mscr{H}_i \mid \se_i \} \leq \de_i$
for $i = 0, 1, \cd, m - 1$ if $\ze$ is sufficiently small. \eeC

Making use of the Monte Carlo method for estimating risks, the idea of bisection risk tuning described in Section \ref{BRT}  and the iterative
minimax optimization algorithm proposed in Section \ref{minmaxcomsimple}, we can determine appropriate risk tuning parameter $\ze$ and weighting
coefficients such that the risk requirement is satisfied.

\section{Tests for the Ratio of Variances of Two Normal Distributions}

Let $X$ be a random variable possessing a normal distribution  with
mean $\mu_X$ and variance $\si_X^2$.  Let $Y$ be a random variable
possessing a normal distribution  with mean $\mu_Y$ and variance
$\si_Y^2$. Define $\se = \f{ \si_X^2 }{ \si_Y^2 }$.  A general
problem regarding $\se$ is to test  $m$ mutually exclusive and
exhaustive composite hypotheses: $\mscr{H}_0: \se \in {\Se}_0, \qu
\mscr{H}_1: \se \in {\Se}_1, \qu \ldots, \qu \mscr{H}_{m - 1}: \se
\in {\Se}_{m - 1}$,   where $\Se_0 = ( 0, \se_1], \; \Se_{m-1} =
(\se_{m - 1}, \iy )$ and $\Se_i = (\se_i, \se_{i+1} ], \; i = 1,
\cd, m - 2$ with $\se_1 < \se_2 < \cd < \se_{m - 1}$. To control the
probabilities of making wrong decisions, it is typically required
that, for pre-specified numbers $\de_i \in (0, 1)$, \bee & & \Pr \li
\{ \tx{Accept} \; \mscr{H}_i \mid \se \ri \} \geq 1 - \de_i  \qu \fa
\se \in \varTheta_i , \qqu i = 0, 1, \cd, m - 1 \eee with
$\varTheta_0 = ( 0, \se_1^{\prime}], \; \varTheta_{m-1} = [ \se_{m -
1}^{\prime \prime}, \iy)$ and $\varTheta_i = [ \se_i^{\prime
\prime}, \se_{i+1}^{\prime} ]$ for $i = 1, \cd, m - 2$, where
$\se_i^\prime, \; \se_i^{\prime \prime}$ satisfy $\se_1^\prime <
\se_1, \; \se_{m-1}^{\prime \prime}
> \se_{m-1}$ and $\se_{i - 1} < \se_{i - 1}^{\prime \prime} \leq
\se_i^\prime < \se_i < \se_i^{\prime \prime} \leq \se_{i + 1}^\prime < \se_{i + 1}$ for $i = 2, \cd, m - 2$.  We shall address this problem for
the case that the mean values are known and the case that the mean values are unknown.   The tests will be defined based on i.i.d. samples $X_1,
X_2, \cd$ of $X$ and i.i.d samples $Y_1, Y_2, \cd$ of $Y$.  It is assumed that $X, Y$ and their samples are mutually independent.

\subsection{Tests with Known Means}

Let $\Up(d_1, d_2, \al)$ denote the $100 \al \%$ quantile of an
$F$-distribution of $d_1$ and $d_2$ degrees of freedom.  That is,
for a chi-square variable, $U$, of $d_1$ degrees of freedom and a
chi-square variable, $V$, of $d_2$ degrees of freedom, $\Pr \li \{
\f{U d_2}{V d_1} \leq  \Up(d_1, d_2, \al) \ri \} = \al$, where $\al
\in (0, 1)$. In the case that the mean values $\mu_X$ and $\mu_Y$
are known, we propose to design multistage plans as follows.

\beT \la{Ratio_known_mean} Suppose that $\al_i = O(\ze) \in (0, 1)$
and $\ba_i = O(\ze) \in (0, 1)$ for $i = 1, \cd, m - 1$. Let $2 \leq
n_1^X < n_2^X < \cd < n_s^X$ and $2 \leq n_1^Y < n_2^Y < \cd <
n_s^Y$ be the sample sizes for variable $X$ and $Y$ respectively.
For $\ell = 1, \cd, s$, define
\[ u_{\ell,i} = \max \{1, \Up (n_\ell^X, n_\ell^Y, 1 - \al_i) \}, \qu  v_{\ell,
i} = \min \{ 1, \Up (n_\ell^X, n_\ell^Y, \ba_i) \}, \qu
\wh{\bs{\se}}_\ell = \f{ n_\ell^Y \sum_{i=1}^{n_\ell^X} (X_i -
\mu_X)^2}{n_\ell^X \sum_{i=1}^{n_\ell^Y} (Y_i - \mu_Y)^2}. \]
Suppose that the maximum sample sizes $n_s^X$ and $n_s^Y$ satisfy
$\se_{i}^{\prime \prime} v_{s, i} \geq \se_{i}^{\prime} u_{s, i}$
for $i = 1, \cd, m - 1$. Define  {\small \bee f_{\ell, i} = \min \li
\{ \se_{i}^{\prime \prime} v_{\ell, i}, \; \f{1}{2} \li (
\se_i^\prime u_{\ell, i} + \se_{i}^{\prime \prime} v_{\ell, i} \ri )
\ri \}, \qqu  g_{\ell, i} = \max \li \{  \se_{i}^{\prime} u_{\ell,
i}, \; \f{1}{2} \li ( \se_i^\prime u_{\ell, i} + \se_{i}^{\prime
\prime} v_{\ell, i} \ri ) \ri \} \eee} for $i = 1, \cd, m - 1$ and
$\ell = 1, \cd, s$.  Define {\small \be \la{defgoognormal}
\bs{D}_\ell = \bec 1 & \tx{for} \; \wh{\bs{\se}}_\ell \leq f_{\ell,1},\\
i & \tx{for} \; g_{\ell, i - 1} < \wh{\bs{\se}}_\ell
\leq f_{\ell, i} \; \tx{where} \; 2 \leq i \leq m - 1,\\
m & \tx{for} \; \wh{\bs{\se}}_\ell > g_{\ell,m-1},\\
0 & \tx{else} \eec \ee} for $\ell = 1, \cd, s$.  Then,  $\Pr \{
\tx{Reject} \; \mscr{H}_i  \mid \se \}$ is no greater than $\de_i$
for any $\se \in \varTheta_i$ and $i = 0, 1, \cd, m - 1$ provided
that $\ze$ is sufficiently small.

 \eeT

See Appendix \ref{proof ratio} for a proof.

\subsection{Tests with Unknown Means}

In the case that the mean values $\mu_X$ and $\mu_Y$ are unknown, we
propose to design multistage plans as follows.

\beT \la{Ratio_unknown_mean} Suppose that $\al_i = O(\ze) \in (0,
1)$ and $\ba_i = O(\ze) \in (0, 1)$ for $i = 1, \cd, m - 1$. Let $2
\leq n_1^X < n_2^X < \cd < n_s^X$ and $2 \leq n_1^Y < n_2^Y < \cd <
n_s^Y$ be the sample sizes for variable $X$ and $Y$ respectively.
For $\ell = 1, \cd, s$, define
\[ u_{\ell,i} = \max
\{1, \Up (n_\ell^X - 1, n_\ell^Y - 1, 1 - \al_i) \}, \qu  v_{\ell,
i} = \min \{1, \Up (n_\ell^X - 1, n_\ell^Y - 1, \ba_i) \} \] and
$\wh{\bs{\se}}_\ell = \f{ (n_\ell^Y - 1) \sum_{i=1}^{n_\ell^X} (X_i
- \ovl{X}_{n_\ell^X} )^2}{ (n_\ell^X - 1) \sum_{i=1}^{n_\ell^Y} (Y_i
- \ovl{Y}_{n_\ell^Y} )^2}$,  where $\ovl{X}_{n_\ell^X}  = \f{
\sum_{i = 1}^{n_\ell^X} X_i } { n_\ell^X  }$ and $\ovl{Y}_{n_\ell^Y}
= \f{ \sum_{i = 1}^{n_\ell^Y} Y_i  } { n_\ell^Y }$.  Suppose that
the maximum sample sizes $n_s^X$ and $n_s^Y$ satisfy
$\se_{i}^{\prime \prime} v_{s, i} \geq \se_{i}^{\prime} u_{s, i}$
for $i = 1, \cd, m - 1$. Define {\small \bee f_{\ell, i} = \min \li
\{ \se_{i}^{\prime \prime} v_{\ell, i}, \; \f{1}{2} \li (
\se_i^\prime u_{\ell, i} + \se_{i}^{\prime \prime} v_{\ell, i} \ri )
\ri \}, \qqu g_{\ell, i} = \max \li \{  \se_{i}^{\prime} u_{\ell,
i}, \; \f{1}{2} \li ( \se_i^\prime u_{\ell, i} + \se_{i}^{\prime
\prime} v_{\ell, i} \ri ) \ri \} \eee} for $i = 1, \cd, m - 1$ and
$\ell = 1, \cd, s$.  Define decision variables $\bs{D}_\ell$ by
(\ref{defgoognormal}) for $\ell = 1, \cd, s$.  Then, $\Pr \{
\tx{Reject} \; \mscr{H}_i  \mid \se \}$ is no greater than $\de_i$
for any $\se \in \varTheta_i$ and $i = 0, 1, \cd, m - 1$ provided
that $\ze$ is sufficiently small.

\eeT

See Appendix \ref{proof ratio} for a proof.  To determine appropriate risk tuning parameter $\ze$ and weighting coefficients such that the risk
requirements are satisfied for the testing plans proposed in Theorems \ref{Ratio_known_mean} and \ref{Ratio_unknown_mean}, we can make use of
the Monte Carlo method for estimating risks, the idea of bisection risk tuning described in Section \ref{BRT} and the iterative minimax
optimization algorithm proposed in Section \ref{minmaxcom}.  In this section, we only consider the general problem of testing multiple
hypotheses.  The general method presented here can be easily adapted to special problems such as testing one-sided hypotheses, two-sided
hypotheses, triple hypotheses, and interval test, etc. Concrete procedures can be worked out by mimicking the techniques presented in Sections 3
and 4.

\sect{Exact Performance Evaluation of SPRT}

To demonstrate the advantages of the present methods upon existing
methods, we shall compare them with the sequential probability ratio
test (SPRT) developed by Wald \cite{Wald}.

Consider an exponential family which consists of density functions of the form  $f_X(x, \se) = h(x) \exp ( \eta (\se) x - \psi(\se) )$, where
$\eta^\prime (\se) > 0$ and $\f{ \psi^\prime (\se) } { \eta^\prime (\se) } = \se$. Let $X_1, X_2, \cd$ be i.i.d. samples of $X$ with a density
function $f_X(x, \se)$ in the exponential family. Let $k_0 \in (0, 1)$ and $k_1 \in (1, \iy)$ be two numbers used to define the stopping rule of
SPRT. The SPRT for testing $\mscr{H}_0: \se \leq \se_0$ versus $\mscr{H}_1: \se \geq \se_1$, where $\se_0 < \se_1$ are parametric values, can be
described as follows:

(i) Reject $\mscr{H}_0$  if $\f{ \prod_{i = 1}^n f_X ( X_i, \se_0) }
{ \prod_{i = 1}^n f_X ( X_i, \se_1) } \leq k_0$;

(ii) Accept $\mscr{H}_0$  if $\f{ \prod_{i = 1}^n f_X ( X_i, \se_0)
} { \prod_{i = 1}^n f_X ( X_i, \se_1) } \geq k_1$;

(iii) Continue sampling if $k_0 < \f{ \prod_{i = 1}^n f_X ( X_i,
\se_0) } { \prod_{i = 1}^n f_X ( X_i, \se_1) } < k_1$.

\bsk

By virtue of the expression of the density function of the
exponential family, the SPRT can be simplified as follows:

(a) Reject $\mscr{H}_0$ if $ [ \eta (\se_0) - \eta (\se_1) ]
(\sum_{i=1}^n X_i) - n [ \psi (\se_0) - \psi(\se_1) ] \leq \ln k_0
$;

(b) Accept $\mscr{H}_0$ if $ [ \eta (\se_0) - \eta (\se_1) ]
(\sum_{i=1}^n X_i)  - n [ \psi (\se_0) - \psi(\se_1) ] \geq \ln
k_1$;

(c) Continue sampling if $\ln k_0  < [ \eta (\se_0) - \eta (\se_1) ]
( \sum_{i=1}^n X_i )  - n [ \psi (\se_0) - \psi(\se_1) ] < \ln k_1$.

\bsk

For exact computation of the OC function and average sample number
of the SPRT, we have established the following results.

\beT

\la{Exact_Computation_SPRT}

Define {\small $u (n) =  \f{  \psi (\se_1) - \psi(\se_0)  - \f{1}{n}
\ln k_1  } { \eta (\se_1) - \eta (\se_0) }$} and {\small $v (n) =
\f{ \psi (\se_1) - \psi(\se_0)  - \f{1}{n} \ln k_0 } { \eta (\se_1)
- \eta (\se_0) }$}  for $n = 1, 2, \cd$. Let $\bf{n}$ denote the
number of samples at the termination of the sampling process of
SPRT. For $\ep > 0$,  the following statements (i) and (ii) hold
true:

(i)  \[ \Pr \{ \tx{Accept $\mscr{H}_0$}, \; \mbf{n} \leq m \mid \se
\} \leq \Pr \{ \tx{Accept $\mscr{H}_0$} \mid \se \} \leq \Pr \{
\tx{Accept $\mscr{H}_0$}, \; \mbf{n} \leq m \mid \se \} + \ep,
\]
\[
\Pr \{ \tx{Accept $\mscr{H}_1$}, \; \mbf{n} \leq m   \mid \se \}
\leq \Pr \{ \tx{Accept $\mscr{H}_1$} \mid \se \} \leq \Pr \{
\tx{Accept $\mscr{H}_1$}, \; \mbf{n} \leq m \mid \se \} + \ep
\]
provided that $\se < u, \; \li [ \f{ \exp ( \eta (\se) u  - \psi (\se) ) } { \exp ( \eta (u) u - \psi (u) ) } \ri ]^m  < \ep$ or $\se > v, \;
\li [ \f{ \exp ( \eta (\se) v - \psi (\se) ) } { \exp ( \eta (v) v - \psi (v) ) } \ri ]^m  < \ep$.

(ii) $\sum_{n = 1}^{m - 1} \Pr \{  \mbf{n} > n  \} \leq \bb{E} [
\mbf{n} ] \leq \sum_{n = 1}^{m - 1} \Pr \{  \mbf{n} > n  \} + \ep$
provided that
\[
\se < u, \qqu \li [ \f{ \exp ( \eta (\se) u - \psi (\se) ) } { \exp ( \eta (u) u - \psi (u) ) } \ri ]^m   < \ep \li [ 1 - \f{ \exp ( \eta (\se)
u - \psi (\se) ) } { \exp ( \eta (u) u - \psi (u) ) } \ri ]
\]
or
\[
\se > v, \qqu \li [ \f{ \exp ( \eta (\se) v - \psi (\se) ) } { \exp ( \eta (v) v - \psi (v) ) } \ri ]^m <  \ep \li [ 1 - \f{ \exp ( \eta (\se) v
- \psi (\se) ) } { \exp ( \eta (v) v - \psi (v) ) } \ri ].
\]

\eeT \bsk

See Appendix \ref{Exact_Computation_SPRT_app} for a proof.

For simplicity of notations, in the above statements (i) and (ii), the abbreviations $u$ and $v$ have been used for $u(m)$ and $v(m)$
respectively. Based on the above bounds, we can apply recursive algorithms to compute the risks and average sample number of SPRT and compare
them with the adaptive methods presented in preceding sections.

\sect{Exact Computation of Boundary Crossing Probabilities}

The problem of evaluating the risk of making incorrect decisions in
multistage hypothesis testing can be addressed in the following
general framework of computing the boundary crossing probabilities.

Let $Y_1, Y_2, \cd, Y_s$ be random variables such that the
increments between consecutive random variables are mutually
independent. Namely, if we define $Z_\ell = Y_{\ell + 1} - Y_{\ell},
\; \ell = 1, \cd, s - 1$, then $Z_1, Z_2, \cd, Z_{s-1}$ are mutually
independent random variables.  It is a frequent problem to compute
the boundary crossing probability $\Pr \li \{  Y_\ell \in
\mscr{Y}_\ell, \; \ell = 1, \cd, s \ri \}$, where $\mscr{Y}_\ell$ is
a subset of the support of $Y_\ell$.   For this purpose, it suffices
to consider $G_\ell (y) \DEF \Pr \li \{ Y_\ell < y; \; Y_i \in
\mscr{Y}_i, \; i = 1, \cd, \ell - 1 \ri \}$ and $g_\ell (y) \DEF
\f{d}{dy} G_\ell (y)$ for $\ell = 1, \cd, s$. Clearly, $g_1 (y)$ is
equal to the probability density function of $Y_1$. Hence, the main
problem is to recursively compute
\[
g_{\ell+1} (y) = \f{d}{dy} \Pr \li \{  Y_{\ell + 1} < y; \; Y_i \in
\mscr{Y}_i, \; i = 1, \cd, \ell \ri \}
\]
for $\ell = 1, \cd, s - 1$.  Let $f_{Z_\ell} (z)$ denote the
probability density function of $Z_\ell$ for $\ell = 1, \cd, s - 1$.
By the independence of $Y_\ell$ and $Z_\ell$, it can be readily
shown that \bee \Pr \li \{ Y_{\ell + 1} < y; \; Y_i \in \mscr{Y}_i,
\; i = 1, \cd, \ell \ri \} =  \int_{u \in \mscr{Y}_\ell } \Pr \li \{
Z_\ell  < y - u \ri \} g_\ell (u) d u, \qqu \fa y \in \mscr{Y}_{\ell
+ 1} \eee or equivalently, \be \la{reck}
 g_{\ell+1} (y) = \int_{u
\in \mscr{Y}_\ell } f_{Z_\ell} (y - u) g_\ell (u) d u, \qqu \fa y
\in \mscr{Y}_{\ell + 1} \ee for $\ell = 1, \cd, s - 1$.  Based on
formula (\ref{reck}), recursive methods have been developed in the
literature for computing the boundary crossing probability. However,
existing methods fail to rigorously control the approximation error,
which is mainly due to the finite partition of the integration
domains.  To overcome the limitation of existing methods, we have
established a new recursive method in the sequel.

Assume that $f_{Z_\ell} (z)$ is differentiable for  for $\ell = 1,
\cd, s - 1$.    Let $0< \vep_1 < \vep_2 < \cd < \vep_s < 1$.  For
$\ell = 1, \cd, s$, we propose to determine $m_\ell$ intervals
$[A_{\ell, i}, B_{\ell, i}], \; i = 1, 2, \cd, m_\ell$ to cover
$\mscr{Y}_\ell$ such that \be \la{imp8} \cup_{i = 1}^{m_\ell}
[A_{\ell, i}, B_{\ell, i}] = \mscr{Y}_\ell, \qqu B_{\ell, i} \leq
A_{\ell, i + 1}, \qu i = 1, \cd, m_\ell - 1 \ee and that \be
\la{imp88} (1 - \vep_\ell) h_{\ell, i} < g_\ell (y) \leq (1 +
\vep_\ell) h_{\ell, i}, \qu \fa y \in [A_{\ell, i}, B_{\ell, i}],
\qu i = 1, \cd, m_\ell. \ee Once this can be accomplished, we have
\[
(1 - \vep_s) \sum_{i = 1}^{m_s} (B_{s,i} - A_{s, i}) h_{s, i}  < \Pr
\li \{  Y_\ell \in \mscr{Y}_\ell, \; \ell = 1, \cd, s \ri \} < (1 +
\vep_s) \sum_{i = 1}^{m_s} (B_{s,i} - A_{s, i}) h_{s, i}.
\]
So, $\sum_{i = 1}^{m_s} (B_{s,i} - A_{s, i}) h_{s, i}$ is an
estimate of $\Pr \li \{  Y_\ell \in \mscr{Y}_\ell, \; \ell = 1, \cd,
s \ri \}$.  The relative precision of such an estimate can be
controlled by $\vep_s$.

The desired intervals for covering $\mscr{Y}_\ell, \; \ell = 1, \cd,
s$ can be constructed recursively.  First, it is not difficult to
determine $m_1$ intervals $[A_{1, i}, B_{1, i}], \; i = 1, 2, \cd,
m_1$ to cover $\mscr{Y}_1$, since $g_1(y)$ is equal to the
probability density function of $Y_1$.  Given that $m_\ell$
intervals $[A_{\ell, i}, B_{\ell, i}], \; i = 1, 2, \cd, m_\ell$
have been determined to cover $\mscr{Y}_\ell$,  we can determine
$[A_{\ell + 1, i}, B_{\ell + 1, i}], \; i = 1, 2, \cd, m_{\ell+1}$
 to cover $\mscr{Y}_{\ell+1}$ by virtue of the following result.

\beT  \la{recursive}  Suppose that $\mscr{Y}_\ell$ is covered by
$m_\ell$ intervals $[A_{\ell, i}, B_{\ell, i}], \; i = 1, 2, \cd,
m_\ell$ such that (\ref{imp8}) and (\ref{imp88}) are satisfied.  Let
$\udl{I} (a, b, c, d)$ and $\ovl{I} (a, b, c, d)$ be multivariate
functions such that $\udl{I} (a, b, c, d) \leq \int_a^b f_{Z_\ell}
(y - u) du \leq \ovl{I} (a, b, c, d)$ for any $y \in [c, d]$ and
that $\ovl{I} (a, b, c, d) - \udl{I} (a, b, c, d) \to 0$ as $d - c
\to 0$.  Let $\eta_\ell$ and $\ga_\ell$ be positive numbers such
that $(1 + \ga_\ell) (1 - \eta_\ell)
> 1 + \eta_\ell$ and
\be \la{assign}
 \vep_{\ell} \leq \f{2}{  1 + \f{ 1 - \vep_{\ell + 1}
}{ 1 + \vep_{\ell + 1}  } \f{ ( 1 + \eta_\ell ) ( 1 + \ga_\ell ) } {
1 - \eta_\ell } } - 1. \ee Let $\ovl{L}_\ell, \ovl{U}_\ell,
\udl{L}_\ell$ and $\udl{U}_\ell$ be positive real numbers such that
\bel & & \ovl{L}_\ell < \sum_{i= 1}^{m_\ell} h_{\ell, i} \; \ovl{I}
(A_{\ell,i}, B_{\ell, i}, C, D) < \ovl{U}_\ell < \f{1 + \eta_\ell}{1
- \eta_\ell} \ovl{L}_\ell , \la{con1}\\
&  & \udl{U}_\ell > \sum_{i= 1}^{m_\ell} h_{\ell, i} \; \udl{I}
(A_{\ell,i}, B_{\ell, i}, C, D) > \udl{L}_\ell > \f{1 - \eta_\ell}{1
+ \eta_\ell} \udl{U}_\ell. \la{con2} \eel Define \be \la{con3}
 h_{\ell + 1} =
\f{1}{4} \li [ ( 1 + \vep_\ell ) (1 + \eta_\ell ) \li ( \ovl{L}_\ell
+ \ovl{U}_\ell \ri ) + ( 1 - \vep_\ell ) (1 - \eta_\ell) \li (
\udl{L}_\ell + \udl{U}_\ell \ri ) \ri ]. \ee Then, $(1 -
\vep_{\ell+1}) h_{\ell+1} < g_{\ell+1} (y) < (1 + \vep_{\ell+1})
h_{\ell+1}$ for any $y \in [C, \; D] \subseteq \mscr{Y}_{\ell + 1}$
provided that \be \la{con4}
 (1 + \ga_\ell) (\udl{L}_\ell + \udl{U}_\ell)
>  \ovl{L}_\ell + \ovl{U}_\ell,
\ee which can be satisfied if $D - C$ is sufficiently small.  \eeT

See Appendix \ref{recursive_app} for a proof.  In Theorem \ref{recursive}, for simplicity, one can take $\ga_\ell = \ga < 1, \;  \eta_\ell =
\f{\ga}{3},  \; \ell = 1, \cd, s - 1$ and determine $0 < \vep_1 < \vep_2 < \cd < \vep_s < 1$ such that \[ \vep_{\ell} = \f{2}{  1 + \f{ 1 -
\vep_{\ell + 1} }{ 1 + \vep_{\ell + 1}  } \f{ ( 1 + \eta_\ell ) ( 1 + \ga_\ell ) } { 1 - \eta_\ell } } - 1, \qqu \ell = 1, \cd, s - 1.
\]
Let $a < b$ and $c < d$.  The multivariate functions $\udl{I} (a, b,
c, d)$ and $\ovl{I} (a, b, c, d)$ in Theorem \ref{recursive} can
 be readily constructed in many situations.  Under the assumption that $f(u)$ is non-decreasing with
respect to $u \leq q$ and non-increasing with respect to $u \geq q$,
we have established that $\udl{I} (a, b, c, d) < \int_a^b f (x - u )
d u < \ovl{I} (a, b, c, d)$ for any $x \in [c, d]$, where {\small
\[ \ovl{I} (a, b, c, d) = \bec   \int^{d-a}_{d-b} f (  u  )
d u  & \tx{for} \;  c - q < d - q < a < b,\\
(d - q  - a) f (q) + \int^{q}_{d-b} f ( u  ) d u &  \tx{for}
\; c - q < a \leq d - q \leq b,\\
(b - a) f (q) & \tx{for} \; c - q < a < b < d - q,\\
(d - c) f (q) + \int^{c-a}_{q} f ( u  ) d u + \int^{q}_{d-b} f ( u
) d u  & \tx{for} \; a \leq c - q  < d - q \leq b,\\
(b - c + q) f (q) + \int^{c-a}_{q} f ( u  ) d u & \tx{for} \;
a \leq c - q \leq b < d - q,\\
\int^{c-a}_{c-b} f (u  ) d u & \tx{for} \; a < b < c - q < d - q
\eec
\]}
and {\small \bee  \udl{I} (a, b,c,d) &  = & \bec \int^{c-a}_{c-b} f
(u) d u  & \tx{for} \;  c - q < d - q < a < b,\\
(d - q - a) \times \min \li \{ f (d - a  ), \; f (c - d + q) \ri \}
& \tx{for}
\; c - q < a \leq d - q \leq b,\\
+ \int^{c-d + q}_{c-b} f (u  ) d u  & \\
(b - a) \times \min \li \{ f (d - a  ), \; f (c - b ) \ri \}
& \tx{for} \; c - q < a < b < d - q,\\
(d - c) \times \min \li \{ f (d - c + q ), \; f (c - d + q) \ri \}
& \tx{for} \; a \leq c - q < d - q \leq b,\\
+ \int^{d-a}_{d-c+q} f (u  ) d u + \int^{c-d+q}_{c-b} f (u ) d u &
\\
(b - c + q) \times \min \li \{ f (d - c + q), \; f (c - b ) \ri \}
& \tx{for} \;
a \leq c - q \leq b < d - q,\\
 + \int^{d - a}_{d - c+q} f (u  ) d u & \\
 \int^{d-a}_{d-b} f (u  ) d u
& \tx{for} \; a < b < c - q < d - q \eec \qqu \qqu \qqu \qu \eee}
Moreover, $\ovl{I} (a, b, c, d) - \udl{I} (a, b, c, d) \to 0$ for
any $x \in [c, d]$ as $d - c \to 0$.

Recall that, given that $\mscr{Y}_\ell$ is covered by $m_\ell$
intervals $[A_{\ell, i}, B_{\ell, i}], \; i = 1, 2, \cd, m_\ell$
such that (\ref{imp8}) and (\ref{imp88}) are satisfied, our
objective is to construct intervals $[A_{\ell+1, i}, B_{\ell+1, i}],
i = 1, 2, \cd, m_{\ell +1}$ to cover $\mscr{Y}_{\ell + 1}$ such that
$\cup_{i = 1}^{m_{\ell + 1}} [A_{{\ell + 1}, i}, B_{{\ell + 1}, i}]
= \mscr{Y}_{\ell + 1}, \;  B_{{\ell + 1}, i} \leq A_{{\ell + 1}, i +
1}, \;  i = 1, \cd, m_{\ell + 1} - 1$ and that $(1 - \vep_{\ell +
1}) h_{{\ell + 1}, i} < g_{\ell + 1} (y) \leq (1 + \vep_{\ell + 1})
h_{{\ell + 1}, i}, \;  \fa y \in [A_{{\ell + 1}, i}, B_{{\ell + 1},
i}], \;  i = 1, \cd, m_{\ell + 1}$.  This can be accomplished by
virtue of Theorem \ref{recursive} as follows.

For simplicity of illustration, we focus on the special case that
$\mscr{Y}_{\ell + 1} = [\udl{y}, \ovl{y}]$.  The general case that
$\mscr{Y}_{\ell + 1}$ consists of multiple subintervals like
$[\udl{y}, \ovl{y}]$ can be addressed by repeatedly applying the
method described in the sequel to each subinterval.

Clearly, for each subinterval $[A_{\ell+1, i}, B_{\ell+1, i}]$,
there exist a lower bound $\ovl{L}_{\ell,i}$ and an upper bound
$\ovl{U}_{\ell,i}$ of $\sum_{j= 1}^{m_\ell} h_{\ell, j} \; \ovl{I}
(A_{\ell,j}, B_{\ell, j}, A_{\ell+1, i}, B_{\ell+1, i})$ such that
\[
\ovl{L}_{\ell, i} < \sum_{j= 1}^{m_\ell} h_{\ell, j} \; \ovl{I}
(A_{\ell,j}, B_{\ell, j}, A_{\ell+1, i}, B_{\ell+1, i}) <
\ovl{U}_{\ell, i} < \f{1 + \eta_\ell}{1 - \eta_\ell} \ovl{L}_{\ell,
i}.
\]
Similarly, there exist a lower bound $\udl{L}_{\ell,i}$ and an upper
bound $\udl{U}_{\ell,i}$ of {\small $\sum_{j= 1}^{m_\ell} h_{\ell,
j} \; \udl{I} (A_{\ell,j}, B_{\ell, j}, A_{\ell+1, i}, B_{\ell+1,
i})$} such that
\[
\udl{U}_{\ell,i} > \sum_{j= 1}^{m_\ell} h_{\ell, j} \; \udl{I}
(A_{\ell,j}, B_{\ell, j}, A_{\ell+1, i}, B_{\ell+1, i})
> \udl{L}_{\ell, i} > \f{1 - \eta_\ell}{1 + \eta_\ell}
\udl{U}_{\ell, i}.
\]
Actually, the bounds $\ovl{L}_{\ell,i}, \; \ovl{U}_{\ell,i}, \;
\udl{L}_{\ell,i}, \; \udl{U}_{\ell,i}$ are multivariate functions of
$A_{\ell+1, i}, B_{\ell+1, i}$ and $\eta_\ell, \; \ga_\ell$. Such
bounds can be calculated by a computer program.   Starting from the
left endpoint of interval $[\udl{y}, \ovl{y}]$, we determine an
initial $[A_{\ell+1, 1}, B_{\ell +1, 1}]$ with $A_{\ell +1, 1} =
\udl{y}$ such that $(1 + \ga_\ell) (\udl{L}_{\ell,1} +
\udl{U}_{\ell, 1})
>  \ovl{L}_{\ell, 1} + \ovl{U}_{\ell, 1}$.  Then, we determine next subinterval $[A_{\ell+1, 2},
B_{\ell +1, 2}]$ as the form
\[ A_{\ell+1, 2} = A_{\ell+1, 1}, \qqu B_{\ell+1, 2} = \min \{\ovl{y}, B_{\ell+1, 1} + (B_{\ell+1, 1} - A_{\ell+1, 1}) 2^{j} \},
\]
with $j$ taken as the largest integer no greater than $1$ to ensure
 $(1 + \ga_\ell) (\udl{L}_{\ell, 2} + \udl{U}_{\ell, 2})
>  \ovl{L}_{\ell, 2} + \ovl{U}_{\ell, 2}$.  For $i > 1$, given interval $[A_{\ell+1, i}, B_{\ell +1, i}]$,  we
determine next subinterval $[A_{\ell+1, {i+1}}, B_{\ell +1, {i+1}}]$
as the form
\[ A_{\ell+1, {i+1}} = B_{\ell +1, i}, \qqu B_{\ell +1, {i+1}} = \min \{\ovl{y}, B_{\ell +1, i} + (B_{\ell +1, i} - A_{\ell +1, i}) 2^{j} \},
\]
with $j$ taken as the largest integer no greater than $1$ to ensure
that $(1 + \ga_\ell) (\udl{L}_{\ell, i+1} + \udl{U}_{\ell, i+1})
>  \ovl{L}_{\ell, i+1} + \ovl{U}_{\ell, i+1}$.  We repeat this process until $B_{\ell +1, i} = \ovl{y}$ for some $i$,
which is taken as the number of intervals $m_{\ell + 1}$.

In the above procedure of constructing intervals to cover
$\mscr{Y}_\ell, \; \ell = 1, \cd, s$, a critical step is to
determine lower and upper bounds for quantities {\small $\sum_{j=
1}^{m_\ell} h_{\ell, j} \; \ovl{I} (A_{\ell,j}, B_{\ell, j},
A_{\ell+1, i}, B_{\ell+1, i})$}  and {\small $\sum_{j= 1}^{m_\ell}
h_{\ell, j} \; \udl{I} (A_{\ell,j}, B_{\ell, j}, A_{\ell+1, i},
B_{\ell+1, i})$} to ensure certain relative precision requirements.
Such quantities can be expressed as the following general form
\[
Q = w_0 + \sum_{i = 1}^m w_i \int_{a_i}^{b_i} f(x) dx,
\]
where $w_0, w_1, \cd, w_m$ are constants.  In the context of
coverage construction for $\mscr{Y}_\ell, \; \ell = 1, \cd, s$, the
number $m$ is very large and the width of each interval $[a_i, b_i]$
is very small.  Hence, there is no need to partition each interval
$[a_i, b_i]$ as many subintervals for purpose of evaluating
$\int_{a_i}^{b_i} f(x) dx$.   Under the assumption that $f(x)$ is
either convex or concave in each $[a_i, b_i]$ (i.e. $f^{\prime
\prime} (x)$ has the same sign in $[a_i, b_i]$),  we propose a {\it
Globally Adaptive Splitting} method for fast computing $Q$ as
follows.

As a consequence of the assumption of convexity on $f(x)$, we have
{\small \be \la{VIP} \f{1}{2} [f(a) + f(b) ] (b - a) +  \min \{0,
\vDe \} \leq \int_a^b f(x) dx \leq \f{1}{2} [f(a) + f(b) ] (b - a) +
\max \{0, \vDe \}, \ee} where $\vDe = \f{1}{8} [ f^\prime(a) -
f^\prime(b) ] (b - a)^2 $.   Applying (\ref{VIP}),  we have $w_0 +
\sum_{i = 1}^m w_i \udl{q}_i < Q < w_0 + \sum_{i = 1}^m w_i
\ovl{q}_i$, where
\[
\udl{q}_i = J_i + \min \{ 0, \vDe_i \}, \qqu  \ovl{q}_i = J_i + \max
\{ 0, \vDe_i \}
\]
with \[ J_i = \f{1}{2} [f(a_i) + f(b_i) ] (b_i - a_i), \qqu \vDe_i =
\f{1}{8} [ f^\prime(a_i) - f^\prime(b_i) ] (b_i - a_i)^2.
\]
Now we find the index $j \in \{1, \cd, m \}$ such that $w_j
\ovl{q}_j - w_j \udl{q}_j = \min \{ w_i \ovl{q}_i - w_i \udl{q}_i: i
= 1, \cd, m \}$.  Then, we split the interval associated with index
$j$ as two subintervals $[a_j, \f{a_j + b_j}{2}]$ and $[\f{a_j +
b_j}{2}, b_j]$ and bound the integrals $\int_{a_j}^{ (a_j + b_j) \sh
2 } f (x) dx$ and $\int^{b_j}_{ (a_j + b_j) \sh 2 } f (x) dx$ by
(\ref{VIP}).  Namely, apply (\ref{VIP}) to determine bounds
$\udl{q}_{j,1}, \; \ovl{q}_{j,1}$ and $\udl{q}_{j,2}, \;
\ovl{q}_{j,2}$ such that
\[
\udl{q}_{j,1} < \int_{a_j}^{ (a_j + b_j) \sh 2 } f (x) dx <
\ovl{q}_{j,1},  \qqu  \udl{q}_{j,2} < \int^{b_j}_{ (a_j + b_j) \sh 2
} f (x) dx < \ovl{q}_{j,2}.
\]
Once theses bounds are computed, we replace the term $w_j \udl{q}_j$
in $w_0 + \sum_{i = 1}^m w_i \udl{q}_i$ by two terms $w_j
\udl{q}_{j,1}$ and $w_j \udl{q}_{j,2}$.  Similarly, we replace the
term $w_j \ovl{q}_j$ in $w_0 + \sum_{i = 1}^m w_i \ovl{q}_i$ by two
terms $w_j \ovl{q}_{j,1}$ and $w_j \ovl{q}_{j,2}$.  Therefore, we
have new lower and upper bounds for $Q$, which can still be
expressed in the form $w_0 + \sum_{i = 1}^m w_i \udl{q}_i < Q < w_0
+ \sum_{i = 1}^m w_i \ovl{q}_i$, where the number $m$, the terms for
summation, and the associated intervals  have been updated.
Repeatedly apply the above splitting technique to the lower and
upper bounds of $Q$. As the splitting process goes on, the gap
between the lower and upper bounds of $Q$ decreases and eventually
we obtain a lower bound $L$ in the form $w_0 + \sum_{i = 1}^m w_i
\udl{q}_i$, and an upper bound $U$ in the form $w_0 + \sum_{i = 1}^m
w_i \ovl{q}_i$ such that
\[
L < Q < U, \qqu \f{U}{1 + \eta} < \f{L}{1 - \eta}
\]
for some pre-specified $\eta \in (0, 1)$.

As can be seen from the above description, the computational
complexity of our recursive method for computing the boundary
crossing probability $\Pr \li \{ Y_\ell \in \mscr{Y}_\ell, \; \ell =
1, \cd, s \ri \}$ depends on the partition of the sets
$\mscr{Y}_\ell, \; \ell = 1, \cd, s$.  For purpose of reducing  the
computational complexity, we can apply the truncation technique to
reduce the domain of integration.  To illustrate, consider a typical
problem of computing $\Pr \{ Y_{\ell + 1} < a_{\ell + 1}; \; a_i <
Y_i < b_i, \; i = 1, \cd, \ell \}$ and $\Pr  \{ Y_{\ell + 1} >
b_{\ell + 1}; \; a_i < Y_i < b_i, \; i = 1, \cd, \ell \}$, which is
frequently encountered in the context of multistage hypothesis
testing.

Let $\ep \in (0, 1)$.  Let $\udl{z}_\ell$ be a number such that $\Pr
\{ Y_{\ell + 1} - Y_\ell > \udl{z}_\ell \}
> 1 - \ep$.  We can show that  \bee & & \Pr \li \{ Y_{\ell + 1} < a_{\ell + 1};
\; a_i < Y_i <
b_i, \; i = 1, \cd, \ell \ri \}\\
&  < &  \Pr \li \{ \udl{z}_\ell + a_\ell < Y_{\ell + 1} < a_{\ell +
1}; \; a_i < Y_i < b_i, \; i = 1, \cd, \ell \ri \} + \ep  \eee and
\bee &  & \Pr \li \{ Y_{\ell + 1} < a_{\ell + 1}; \; a_i < Y_i <
b_i, \; i = 1, \cd, \ell \ri \}\\
&  > &  \Pr \li \{ \udl{z}_\ell + a_\ell < Y_{\ell + 1} < a_{\ell +
1}; \; a_i < Y_i < b_i, \; i = 1, \cd, \ell \ri \} - \ep.  \eee On
the other hand, letting $\ovl{z}_\ell$ be a number such that $\Pr \{
Y_{\ell + 1} - Y_\ell < \ovl{z}_\ell \} > 1 - \ep$, we can establish
that \bee & & \Pr \li \{ Y_{\ell + 1} > b_{\ell + 1}; \; a_i < Y_i <
b_i, \; i = 1, \cd, \ell \ri \}\\
&  < &  \Pr \li \{  b_{\ell + 1} < Y_{\ell + 1} < \ovl{z}_\ell +
b_\ell; \; a_i < Y_i < b_i, \; i = 1, \cd, \ell \ri \} + \ep  \eee
and \bee &  & \Pr \li \{ Y_{\ell + 1} > b_{\ell + 1}; \; a_i < Y_i <
b_i, \; i = 1, \cd, \ell \ri \}\\
& > &  \Pr \li \{  b_{\ell + 1} < Y_{\ell + 1} < \ovl{z}_\ell +
b_\ell; \; a_i < Y_i < b_i, \; i = 1, \cd, \ell \ri \} - \ep. \eee
Applying Theorem \ref {Trun_THM5} and the truncation method
described in Section \ref{secDT}, we can further reduce the
 complexity for computing $\Pr \{ \udl{z}_\ell + a_\ell < Y_{\ell + 1} < a_{\ell +
1}; \; a_i < Y_i < b_i, \; i = 1, \cd, \ell \}$ and $\Pr \{ b_{\ell
+ 1} < Y_{\ell + 1} < \ovl{z}_\ell + b_\ell; \; a_i < Y_i < b_i, \;
i = 1, \cd, \ell \}$.

After employing the truncation technique to reduce the domain of
integration, one can use our recursive method to compute the
relevant boundary crossing probabilities.

 \sect{Conclusion}

In this paper, we have established a new framework of multistage
hypothesis tests which applies to arbitrary number of mutually
exclusive and exhaustive composite hypotheses.  Specific testing
plans for common problems have also been developed. Our test plans
have several important advantages upon existing tests. First, our
tests are more efficient. Second, our tests always guarantee
prescribed requirement of power. Third, the sample number or test
time of our tests are absolutely bounded. Such advantages have been
achieved by means of new structure of testing plans and powerful
computational machinery.

\bsk

\appendix

\sect{Preliminary Results}

We need some preliminary results. The following Lemmas
\ref{ProbTrans} and \ref{ULE_Basic} have been established in
\cite{Chen_EST}.

\beL \la{ProbTrans} $\Pr \{ F_Z(Z) \leq \al \} \leq \al$ and $\Pr \{
G_Z(Z) \leq \al \} \leq \al$ for any random variable $Z$ and
positive number $\al$. \eeL

\beL

\la{ULE_Basic}

Let $\mscr{E}$ be an event dependent only on random tuple $(X_1,
\cd, X_{\mbf{r}})$. Let $\varphi(X_1, \cd, X_{\mbf{r}})$ be a ULE of
$\se$. Then,

(i) $\Pr \{ \mscr{E} \mid \se \}$ is non-increasing with respect to
$\se \in {\Se}$ no less than $z$ provided that $\mscr{E} \subseteq
\{ \varphi(X_1, \cd, X_{\mbf{r}}) \leq z \}$.

(ii) $\Pr \{ \mscr{E} \mid \se \}$ is non-decreasing with respect to
$\se \in {\Se}$ no greater than $z$ provided that $\mscr{E}
\subseteq \{ \varphi(X_1, \cd, X_{\mbf{r}}) \geq z \}$.

\eeL

\beL \la{Unify_CH} Let $X$ be a random variable parameterized by its
mean $\bb{E} [ X ] = \se \in \Se$.  Suppose that $X$ is a ULE of
$\se$. Let $\ovl{X}_n = \f{ \sum_{i = 1}^n X_i } {n}$, where $X_1,
\cd, X_n$ are  i.i.d. samples of random variable $X$.  Then, \bee &
& \Pr
\{ \ovl{X}_n \leq z \} \leq [ \mscr{C} (z, \se) ]^n, \qqu \fa z \leq \se\\
&  & \Pr \{ \ovl{X}_n \geq z \} \leq [ \mscr{C} (z, \se) ]^n, \qqu
\fa z \geq \se. \eee Moreover, $\mscr{C} (z, \se)$ is non-decreasing
with respect to $\se$ no greater than $z$ and is non-increasing with
respect to $\se$ no less than $z$.  Similarly,  $\mscr{C} (z, \se)$
is non-decreasing with respect to $z$ no greater than $\se$ and is
non-increasing with respect to $z$ no less than $\se$.

\eeL

\bpf By the convexity of function $e^x$ and Jensen's inequality, we
have $\inf_{\ro > 0 } \bb{E} [e^{\ro (X - z)}] \geq \inf_{\ro > 0 }
e^{\ro \bb{E} [ X - z]} \geq 1$ for $\se \geq z$. In view of
$\inf_{\ro \leq 0 } \bb{E} [e^{\ro (X - z)}] \leq 1$, we have
$\mscr{C} (z, \se) = \inf_{\ro \leq 0 } \bb{E} [e^{\ro (X - z)}]$
for $\se \geq z$. Clearly, $\mscr{C} (z, \se) = \inf_{\ro \leq 0 }
e^{- \ro z} \bb{E} [e^{\ro X}]$ is non-decreasing with respect to
$z$ less than $\se$. Since $X$ is a ULE of $\se$, we have that
$\bb{E} [e^{\ro (X - z)}] = e^{- \ro z} \bb{E} [e^{\ro X} ] = e^{-
\ro z} \int_{u = 0}^\iy \Pr \{ e^{\ro X}
> u \} du$ is non-increasing with respect to $\se \geq z$ for $\ro
\leq 0$ and thus $\mscr{C} (z, \se)$ is non-increasing with respect
to $\se$ greater than $z$.

Observing that $\inf_{\ro \geq 0 } \bb{E} [e^{\ro (X - z)}] \leq 1$
and that $\inf_{\ro < 0 } \bb{E} [e^{\ro (X - z)}] \geq \inf_{\ro <
0 } e^{\ro \bb{E} [ X - z]} \geq 1$ for $\se < z$, we have $\mscr{C}
(z, \se) = \inf_{\ro \geq 0 } \bb{E} [e^{\ro (X - z)}]$ for $\se <
z$. Clearly, $\mscr{C} (z, \se) = \inf_{\ro \geq 0 } e^{- \ro z}
\bb{E} [e^{\ro X}]$ is non-increasing with respect to $z$ greater
than $\se$. Since $X$ is a ULE of $\se$, we have that $\bb{E}
[e^{\ro (X - z)}] = e^{- \ro z} \int_{u = 0}^\iy \Pr \{ e^{\ro X} >
u \} du$ is non-decreasing with respect to $\se$ for $\ro > 0$ and
consequently,  $\mscr{C} (z, \se)$ is non-decreasing with respect to
$\se$ smaller than $z$.

Making use of the established fact $\inf_{\ro \leq 0 } \bb{E}
[e^{\ro (X - z)}] = \mscr{C} (z, \se)$ and the Chernoff bound $\Pr
\{ \ovl{X}_n \leq z \}  \leq \li [ \inf_{\ro \leq 0 } \bb{E} [e^{\ro
(X - z)}] \ri ]^n$ (see, \cite{Chernoff}), we have $\Pr \{ \ovl{X}_n
\leq z \} \leq \li [ \mscr{C} (z, \se) \ri ]^n$ for $z \leq \se$.
Making use of the established fact $\inf_{\ro \geq 0 } \bb{E}
[e^{\ro (X - z)}] = \mscr{C} (z, \se)$ and the Chernoff bound $\Pr
\{ \ovl{X}_n \geq z \} \leq \li [ \inf_{\ro \geq 0 } \bb{E} [e^{\ro
(X - z)}] \ri ]^n$, we have $\Pr \{ \ovl{X}_n \geq z \} \leq [
\mscr{C} (z, \se) ]^n$ for $z \geq \se$. This concludes the proof of
Lemma \ref{Unify_CH}. \epf

\sect{Proof of Theorem \ref{Multi_Comp_Exact} }
\la{Multi_Comp_Exact_Ap}

For arbitrary parametric values $\se_0 < \se_1$ in $\Se$, by the
assumption that $\bs{\varphi}_n$ converges in probability to $\se$,
we have that $\Pr \{ \bs{\varphi}_n \geq \f{\se_0 + \se_1}{2} \mid
\se_0 \} \leq \Pr \{ | \bs{\varphi}_n - \se_0 | \geq \f{\se_1 -
\se_0}{2} \mid \se_0 \} \to 0$ and $\Pr \{ \bs{\varphi}_n \leq
\f{\se_0 + \se_1}{2} \mid \se_1 \} \leq \Pr \{ | \bs{\varphi}_n -
\se_1 | \geq \f{\se_1 - \se_0}{2} \mid \se_1 \} \to 0$ as $n \to
\iy$.  This shows that $\ovl{n}$ exists and is finite.

Since $F_{\wh{\bs{\se}}_\ell} (z, \se) = \Pr \{ \wh{\bs{\se}}_\ell
\leq z \mid \se \} = 1 - \Pr \{ \wh{\bs{\se}}_\ell > z \mid \se \}
$, making use of Lemma \ref{ULE_Basic} and the assumption that
$\wh{\bs{\se}}_\ell$ is a ULE of $\se$, we have that
$F_{\wh{\bs{\se}}_\ell} (z, \se)$ is non-increasing with respect to
$\se \in {\Se}$.  Similarly, since $G_{\wh{\bs{\se}}_\ell} (z, \se)
= \Pr \{ \wh{\bs{\se}}_\ell \geq z \mid \se \} = 1 - \Pr \{
\wh{\bs{\se}}_\ell < z \mid \se \} $, making use of Lemma
\ref{ULE_Basic} and the assumption that $\wh{\bs{\se}}_\ell$ is a
ULE of $\se$, we have that $G_{\wh{\bs{\se}}_\ell} (z, \se)$ is
non-decreasing with respect to $\se \in {\Se}$.

To show statement (I), notice that $\{ \tx{Reject} \; \mscr{H}_0 \}
\subseteq \{  \wh{\bs{\se}} \geq \se_1^\prime \}$ as a consequence
of the definition of the test plan.  Hence, statement (I) is proved
by virtue of Lemma \ref{ULE_Basic}.

To show statement (II), notice that $\{ \tx{Reject} \;
\mscr{H}_{m-1} \} \subseteq \{  \wh{\bs{\se}} \leq \se_{m-1}^{\prime
\prime} \}$ as a consequence of the definition of the test plan.
Hence, statement (II) is proved by virtue of Lemma \ref{ULE_Basic}.

To show statement (III), we first claim that $\Pr \{ 1 \leq
\bs{D}_\ell \leq i \mid \se \} \leq \ovl{\ba}_i$ for $0 \leq i \leq
m - 1$ and $\se \in \varTheta_i$.  Clearly, $\{ \wh{\bs{\se}}_\ell
\leq f_{\ell, j}  \} = \{ \wh{\bs{\se}}_\ell \leq \udl{f} (n_\ell,
\se_j^\prime, \se_j^{\prime \prime}, \al_j, \ba_j) \} \subseteq \{
\wh{\bs{\se}}_\ell \leq f (n_\ell, \se_j^{\prime \prime}, \ba_j) \}$
for $1 \leq j \leq i$.  Since $F_{\wh{\bs{\se}}_\ell} (z, \se)$ is
non-decreasing with respect to $z$, we have $\{ \wh{\bs{\se}}_\ell
\leq f (n_\ell, \se_j^{\prime \prime}, \ba_j) \} \subseteq \{
\wh{\bs{\se}}_\ell \leq \se_j^{\prime \prime},
F_{\wh{\bs{\se}}_\ell} (\wh{\bs{\se}}_\ell, \se_j^{\prime \prime})
\leq \ba_j \} \subseteq \{ F_{\wh{\bs{\se}}_\ell}
(\wh{\bs{\se}}_\ell, \se_j^{\prime \prime}) \leq \ba_j \}$ for $1
\leq j \leq i$.  Recalling that $F_{\wh{\bs{\se}}_\ell} (z, \se)$ is
non-increasing with respect to $\se \in {\Se}$ and invoking Lemma
\ref{ProbTrans}, we have \be \la{ying}
 \Pr \{ \wh{\bs{\se}}_\ell
\leq f_{\ell, j} \mid \se \}  \leq  \Pr \{ F_{\wh{\bs{\se}}_\ell}
(\wh{\bs{\se}}_\ell, \se_j^{\prime \prime}) \leq \ba_j \mid \se \}
\leq \Pr \{ F_{\wh{\bs{\se}}_\ell} (\wh{\bs{\se}}_\ell, \se) \leq
\ba_j \mid \se \} \leq \ba_j \leq \ovl{\ba}_i  \ee for $1 \leq j
\leq i$ and $\se \in \varTheta_i$.  For $i = 0$, it is clear that
$\Pr \{ 1 \leq \bs{D}_\ell \leq i \mid \se \} = 0 \leq \ovl{\ba}_0$
for $\se \in \varTheta_0$.  For $i = 1$, by virtue of (\ref{ying}),
we have $\Pr \{ 1 \leq \bs{D}_\ell \leq i \mid \se \} = \Pr \{
\wh{\bs{\se}}_\ell \leq f_{\ell, 1} \mid \se \} \leq \ovl{\ba}_1$
for $\se \in \varTheta_1$.  For $2 \leq i \leq m - 1$, define $S =
\{ j : g_{\ell, j-1} < f_{\ell, j}, \; 2 \leq j \leq i \}$ and let
$r$ be an integer such that $r$ assumes value $1$ if $S$ is empty
and that $r \in S, \; f_{\ell, r} = \max \{f_{\ell, j}: j \in S \}$
if $S$ is not empty. It follows from (\ref{ying}) that $\Pr \{ 1
\leq \bs{D}_\ell \leq i \mid \se \} \leq \Pr \{ \wh{\bs{\se}}_\ell
\leq f_{\ell, r} \mid \se \} \leq \ovl{\ba}_i$ for $2 \leq i \leq m
- 1$ and $\se \in \varTheta_i$. This proves our first claim.  Next,
we claim that $\Pr \{ i + 2 \leq \bs{D}_\ell \leq m \mid \se \} \leq
\ovl{\al}_i$ for $0 \leq i \leq m - 1$ and $\se \in \varTheta_i$.
Clearly, $\{ \wh{\bs{\se}}_\ell > g_{\ell, j}  \} = \{
\wh{\bs{\se}}_\ell > \ovl{g} (n_\ell, \se_j^\prime, \se_j^{\prime
\prime}, \al_j, \ba_j) \} \subseteq \{ \wh{\bs{\se}}_\ell \geq g
(n_\ell, \se_j^\prime, \al_j) \}$ for $i < j \leq m-1$. Since
$G_{\wh{\bs{\se}}_\ell} (z, \se)$ is non-increasing with respect to
$z$, we have $\{ \wh{\bs{\se}}_\ell \geq g (n_\ell, \se_j^\prime,
\al_j) \} \subseteq \{ \wh{\bs{\se}}_\ell \geq \se_j^\prime,
G_{\wh{\bs{\se}}_\ell} (\wh{\bs{\se}}_\ell, \se_j^\prime) \leq \al_j
\} \subseteq \{ G_{\wh{\bs{\se}}_\ell} (\wh{\bs{\se}}_\ell,
\se_j^\prime) \leq \al_j \}$ for $i < j \leq m-1$.  Recalling that
$G_{\wh{\bs{\se}}_\ell} (z, \se)$ is non-decreasing with respect to
$\se \in {\Se}$ and invoking Lemma \ref{ProbTrans}, we have \be
\la{yingchen} \Pr \{ \wh{\bs{\se}}_\ell > g_{\ell, j} \mid \se \}
\leq  \Pr \{ G_{\wh{\bs{\se}}_\ell} (\wh{\bs{\se}}_\ell,
\se_j^\prime) \leq \al_j \mid \se \} \leq \Pr \{
G_{\wh{\bs{\se}}_\ell} (\wh{\bs{\se}}_\ell, \se) \leq \al_j \mid \se
\} \leq \al_j \leq \ovl{\al}_i  \ee for $i < j \leq m-1$ and $\se
\in \varTheta_i$.  For $i = m-1$,  it is evident that $\Pr \{ i + 2
\leq \bs{D}_\ell \leq m \mid \se \} = 0 \leq \ovl{\al}_{m-1}$ for
$\se \in \varTheta_{m-1}$. For $i = m-2$, making use of
(\ref{yingchen}), we have $\Pr \{ i + 2 \leq \bs{D}_\ell \leq m \mid
\se \} = \Pr \{ \wh{\bs{\se}}_\ell > g_{\ell, m-1} \mid \se \} \leq
\ovl{\al}_{m-2}$ for $\se \in \varTheta_{m-2}$.  For $0 \leq i \leq
m - 3$, define $S = \{ j : g_{\ell, j-1} < f_{\ell, j}, \; i + 2
\leq j \leq m-1 \}$, and let $r$ be an integer such that $r$ assumes
value $m-1$ if $S$ is empty and that $r \in S, \; g_{\ell, r - 1} =
\min \{g_{\ell, j - 1}: j \in S \}$ if $S$ is not empty. It follows
from (\ref{yingchen}) that $\Pr \{ i + 2 \leq \bs{D}_\ell \leq m
\mid \se \} \leq \Pr \{ \wh{\bs{\se}}_\ell > g_{\ell, r - 1} \mid
\se \} \leq \ovl{\al}_i$ for $0 \leq i \leq m - 3$ and $\se \in
\varTheta_i$. This proves our second claim. Making use of these two
established claims, we have \be \la{simp8} \Pr \{ \tx{Reject} \;
\mscr{H}_i, \; \bs{l} = \ell  \mid \se \} \leq \Pr \{ 1 \leq
\bs{D}_\ell \leq i \mid \se \} + \Pr \{ i + 2 \leq \bs{D}_\ell \leq
m \mid \se \} \leq \ovl{\al}_i + \ovl{\ba}_i \ee for $i = 0, 1, \cd,
m - 1, \; \se \in \varTheta_i$ and  $\ell = 1, \cd, s$. It follows
that $\Pr \{ \tx{Reject} \; \mscr{H}_i \mid \se \} \leq \sum_{\ell =
1}^s [ \Pr \{ 1 \leq \bs{D}_\ell \leq i \mid \se \} + \Pr \{ i + 2
\leq \bs{D}_\ell \leq m \mid \se \} ] \leq \sum_{\ell = 1}^s
(\ovl{\al}_i + \ovl{\ba}_i)$ for $i = 0, 1, \cd, m - 1$ and $\se \in
\varTheta_i$.  This establishes statement (III).

Statements (IV) and (V) can be shown by virtue of Lemma
\ref{ULE_Basic} and the observation that $\{ \tx{Accept} \;
\mscr{H}_i \} \subseteq \{ \se_i^\prime \leq \wh{\bs{\se}} \leq
\se_{i+1}^{\prime \prime} \}$ and that $\{ \tx{Accept} \; \mscr{H}_i
\}$ is determined by the random tuple $(X_1, \cd, X_{\mbf{n}})$ as a
consequence of the definition of the testing plan.

We now want to show statement (VI).  Observing that $G_{\varphi_n}
(z, \se)$ is non-increasing with respect to $z$, we have that  $g(n,
\se_i^\prime, \al_i) \leq \f{\se_i^\prime + \se_i^{\prime
\prime}}{2}$ if $G_{\varphi_n} ( \f{\se_i^\prime + \se_i^{\prime
\prime}}{2}, \se_i^\prime ) \leq \al_i$.  Since $\varphi_n =
\f{\sum_{i = 1}^n X_i}{n}$ is an unbiased ULE for $\se$, it follows
from Lemma \ref{Unify_CH} that
\[
G_{\varphi_n} \li ( \f{\se_i^\prime + \se_i^{\prime \prime}}{2},
\se_i^\prime \ri ) = \Pr \li \{ \varphi_n \geq \f{ \se_i^{\prime
\prime} + \se_i^\prime }{2}  \mid \se_i^\prime \ri \} \leq \li [
\mcal{C} \li ( \f{\se_i^\prime + \se_i^{\prime \prime}}{2},
\se_i^\prime \ri ) \ri ]^n \leq \al_i
\]
if {\small $n \geq \f{ \ln (\al_i) } { \ln \mcal{C} (
\f{\se_i^\prime + \se_i^{\prime \prime}}{2}, \se_i^\prime ) }$}.  On
the other hand, observing that $F_{\varphi_n} (z, \se)$ is
non-decreasing with respect to $z$, we have that $f(n, \se_i^{\prime
\prime}, \ba_i) \geq \f{\se_i^\prime + \se_i^{\prime \prime}}{2}$ if
$F_{\varphi_n} ( \f{\se_i^\prime + \se_i^{\prime \prime}}{2},
\se_i^{\prime \prime} ) \leq \ba_i$. Since $\varphi_n$ is an
unbiased ULE for $\se$, it follows from Lemma \ref{Unify_CH} that
\[
F_{\varphi_n} \li ( \f{\se_i^\prime + \se_i^{\prime \prime}}{2},
\se_i^{\prime \prime} \ri ) = \Pr \li \{ \varphi_n \leq \f{
\se_i^{\prime \prime} + \se_i^\prime }{2}  \mid \se_i^{\prime
\prime} \ri \} \leq \li [ \mcal{C} \li ( \f{\se_i^\prime +
\se_i^{\prime \prime}}{2}, \se_i^{\prime \prime} \ri ) \ri ]^n \leq
\ba_i
\]
if {\small $n \geq \f{ \ln (\ba_i) } { \ln \mcal{C} (
\f{\se_i^\prime + \se_i^{\prime \prime}}{2}, \se_i^{\prime \prime} )
}$}. Therefore, $f(n, \se_i^{\prime \prime}, \ba_i) \geq g(n,
\se_i^\prime, \al_i)$ if
\[
n \geq \max \li \{  \f{ \ln (\al_i) } { \ln \mcal{C} (
\f{\se_i^\prime + \se_i^{\prime \prime}}{2}, \se_i^\prime ) },  \f{
\ln (\ba_i) } { \ln \mcal{C} ( \f{\se_i^\prime + \se_i^{\prime
\prime}}{2}, \se_i^{\prime \prime} ) } \ri \}.
\]
Let $\ovl{n}$ be the minimum integer $n$ such that $f (n,
\se_i^{\prime \prime}, \ba_i) \geq g (n, \se_i^{\prime}, \al_i)$ for
$i = 1, \cd, m-1$.  Then, $\{ \bs{l} \leq \ovl{n} \}$ is a sure
event and
\[
\ovl{n} \leq \max_{i \in \{1, \cd, m - 1 \} } \max \li \{  \f{ \ln
(\al_i) } { \ln \mcal{C} ( \f{\se_i^\prime + \se_i^{\prime
\prime}}{2}, \se_i^\prime ) },  \f{ \ln (\ba_i) } { \ln \mcal{C} (
\f{\se_i^\prime + \se_i^{\prime \prime}}{2}, \se_i^{\prime \prime} )
} \ri \} = O \li ( \ln \f{1}{\ze} \ri ).
\]
Noting that \bee \Pr \{ \tx{Reject} \; \mscr{H}_i  \mid \se \} =
\sum_{\ell = 1}^{\min \{ s, \ovl{n} \}} \Pr \{ \tx{Reject} \;
\mscr{H}_i, \; \bs{l} = \ell \mid \se \} \eee and making use of
(\ref{simp8}), we have that, as $\ze \to 0$,
\[
\Pr \{ \tx{Reject} \; \mscr{H}_i  \mid \se \} \leq \ovl{n} \; (
\ovl{\al}_i  + \ovl{\ba}_i  ) = O \li ( \ln \f{1}{\ze} \ri ) O( \ze)
\to 0
\]
 for any $\se \in \varTheta_i$ and
$i = 0, 1, \cd, m - 1$.  This proves statement (VI).

To show statement (VII), by the definition of the test plan, we have
that $\{ \tx{Reject} \; \mscr{H}_i \}$ is determined by the random
tuple $(X_1, \cd, X_{\mbf{n}})$.  Moreover, for any numbers $a$ and
$b$ such that $\se_i^{\prime \prime} \leq a < b \leq
\se_{i+1}^\prime$, we have that $\{ \tx{Reject} \; \mscr{H}_i \} =
\{ \tx{Reject} \; \mscr{H}_i, \; \wh{\bs{\se}}  \leq a \} \cup \{
\tx{Reject} \; \mscr{H}_i, \; \wh{\bs{\se}} \geq b \}$ and $\{
\tx{Reject} \; \mscr{H}_i, \; \wh{\bs{\se}}  \leq a \} \cap \{
\tx{Reject} \; \mscr{H}_i, \; \wh{\bs{\se}} \geq b \} = \emptyset$,
which imply that $\Pr \{ \tx{Reject} \; \mscr{H}_i \mid \se \} = \Pr
\{ \tx{Reject} \; \mscr{H}_i, \; \wh{\bs{\se}}  \leq a \mid \se \} +
\Pr \{ \tx{Reject} \; \mscr{H}_i, \; \wh{\bs{\se}} \geq b \mid \se
\}$.  By Lemma \ref{ULE_Basic}, we have that $\Pr \{ \tx{Reject} \;
\mscr{H}_i, \; \wh{\bs{\se}}  \leq a \mid \se \}$ is non-increasing
with respect to $\se \in {\Se}$ no less than $a$ and that $\Pr \{
\tx{Reject} \; \mscr{H}_i, \; \wh{\bs{\se}}  \geq b \mid \se \}$ is
non-decreasing with respect to $\se \in {\Se}$ no greater than $b$.
This leads to the upper and lower bounds of $\Pr \{ \tx{Reject} \;
\mscr{H}_i \mid \se \}$ in statement (VII).

Statement (VIII) can be shown by virtue of Lemma \ref{ULE_Basic}
based on the observation that

$\{ \tx{Reject} \; \mscr{H}_0 \; \tx{and} \; \mscr{H}_{m-1} \}
\subseteq \{  \se_1^\prime \leq \wh{\bs{\se}} \leq \se_{m-1}^{\prime
\prime} \}$ and that $\{ \tx{Reject} \; \mscr{H}_0 \; \tx{and} \;
\mscr{H}_{m-1} \}$ is determined by the random tuple $(X_1, \cd,
X_{\mbf{n}})$ as a consequence of the definition of the test plan.

Finally,  we shall show statement (IX). Note that $\Pr \{
\tx{Reject} \; \mscr{H}_0 \; \tx{and} \;  \mscr{H}_{m-1} \mid \se \}
\leq \sum_{\ell = 1}^s \Pr \{ 2 \leq \bs{D}_\ell \leq m -1 \mid \se
\}$. Define $S = \{ j : g_{\ell, j-1} < f_{\ell, j}, \; 2 \leq j
\leq m-1\}$.  In the case that $S$ is empty, $\Pr \{ 2 \leq
\bs{D}_\ell \leq m-1 \mid \se \} = 0$. In the case that $S$ is not
empty, let $r \in S$ be an integer such $f_{\ell, r} = \max
\{f_{\ell, j}: j \in S \}$. Then, $\Pr \{ 2 \leq \bs{D}_\ell \leq
m-1 \mid \se \} \leq \Pr \{ \wh{\bs{\se}}_\ell \leq f_{\ell, r} \mid
\se \} \leq \max \{ \ba_j: 2 \leq j \leq m -1 \}$ for $\se \in
\varTheta_{m-1}$. On the other hand, if we let $r \in S$ be an
integer such that $g_{\ell, r - 1} = \min \{g_{\ell, j - 1}: j \in S
\}$, then $\Pr \{ 2 \leq \bs{D}_\ell \leq m-1 \mid \se \} \leq \Pr
\{ \wh{\bs{\se}}_\ell \geq g_{\ell, r - 1} \mid \se \} \leq \max \{
\al_j: 1 \leq j \leq m - 2 \}$ for $\se \in \varTheta_0$.  This
proves statement (IX) and concludes the proof of the theorem.

\sect{Proof of Theorem \ref{boundCDFLR} } \la{boundCDFLR_app}

For simplicity of notations, define $F (z, \se) = \Pr \{ \varphi_n \leq z \mid \se \}$ and $G (z, \se) = \Pr \{ \varphi_n \geq z \mid \se \}$.
By the assumption of the theorem, {\small $\f{ f_n( X_1, \cd, X_n; \se) }{f_n( X_1, \cd, X_n; \wh{\se}_n)}  = \Lm (\varphi_n , \wh{\se}_n,
\se)$}.  By virtue of Theorem 1 in page $3$ of the $4$-th version of our paper \cite{ChenR} published in arXiv, we have  \bee &  & \Pr \li \{
\f{ f_n( X_1, \cd, X_n; \se) }{f_n( X_1, \cd, X_n; \wh{\se}_n)} \leq \f{\al}{2}, \; \wh{\se}_n \leq \se, \; \mid \se \ri \} = \Pr \li \{ \Lm
(\varphi_n , \wh{\se}_n, \se) \leq
\f{\al}{2}, \; \wh{\se}_n \leq \se \mid \se \ri \}\\
&  & =  \Pr \li \{ \Lm (\varphi_n , \wh{\se}_n, \se)  \leq \f{\al}{2}, \; \wh{\se}_n \leq \se \mid \se \ri \}
\leq  \Pr \li \{  F (\varphi_n, \se)  \leq \f{\al}{2}, \; \wh{\se}_n \leq \se \mid \se \ri \}\\
&  & \leq  \Pr \li \{ F (\varphi_n, \se)  \leq \f{\al}{2} \mid \se \ri \} \leq \f{\al}{2} \eee for any $\se \in \Se$. This proves (\ref{RB1}).
Similarly, for any $\se \in \Se$, \bee &  & \Pr \li \{   \f{ f_n( X_1, \cd, X_n; \se) }{f_n( X_1, \cd, X_n; \wh{\se}_n)}  \leq \f{\al}{2}, \;
\wh{\se}_n \geq \se \mid \se \ri \} = \Pr \li \{ \Lm (\varphi_n , \wh{\se}_n, \se) \leq
\f{\al}{2}, \; \wh{\se}_n \geq \se \mid \se \ri \}\\
&  & =   \Pr \li \{ \Lm (\varphi_n , \wh{\se}_n, \se)  \leq \f{\al}{2}, \; \wh{\se}_n \geq \se \mid \se \ri \} \leq \Pr \li \{ G (\varphi_n,
\se)  \leq \f{\al}{2}, \; \wh{\se}_n \geq \se \mid \se \ri \}\\
&  & \leq  \Pr \li \{ G (\varphi_n, \se)  \leq \f{\al}{2} \mid \se \ri \} \leq \f{\al}{2}, \eee which establishes (\ref{RB2}).  To show
(\ref{RB3}), making use of (\ref{RB1}) and (\ref{RB2}), we have
 \bee &  & \Pr \li \{ \f{ f_n( X_1, \cd, X_n; \se) }{f_n( X_1, \cd, X_n; \wh{\se}_n)} \leq \f{\al}{2} \mid \se \ri \}\\
&  & =  \Pr \li \{  \f{ f_n( X_1, \cd, X_n; \se) }{f_n( X_1, \cd, X_n; \wh{\se}_n)} \leq \f{\al}{2}, \; \wh{\se}_n \leq \se \mid \se \ri \}
 + \Pr \li \{  \f{ f_n( X_1, \cd, X_n; \se) }{f_n( X_1, \cd, X_n; \wh{\se}_n)} \leq \f{\al}{2}, \; \wh{\se}_n \geq \se \mid \se \ri \}\\
&  & \leq  \f{\al}{2} + \f{\al}{2} = \al \eee for any $\se \in \Se$. To show (\ref{LBA}), making use of (\ref{RB1}), we have  that \bee &  & \Pr
\li \{  \f{ \sup_{\vse \in \mscr{S}} f_n( X_1, \cd, X_n; \vse) }{\sup_{\vse \in \Se} f_n( X_1, \cd, X_n; \vse)} \leq \f{\al}{2}, \; \wh{\se}_n
\leq \inf
\mscr{S} \mid \se \ri \}\\
&  & \leq \Pr \li \{  \f{ \sup_{\vse \in \mscr{S}} f_n( X_1, \cd, X_n; \vse) }{\sup_{\vse \in \Se} f_n( X_1, \cd, X_n; \vse)} \leq \f{\al}{2},
\; \wh{\se}_n \leq \se \mid \se \ri \}\\
&  & \leq \Pr \li \{  \f{ f_n( X_1, \cd, X_n; \se) }{\sup_{\vse \in \Se} f_n( X_1, \cd, X_n; \vse)} \leq \f{\al}{2}, \; \wh{\se}_n \leq \se \mid
\se \ri \}\\
&  & = \Pr \li \{  \f{ f_n( X_1, \cd, X_n; \se) }{f_n( X_1, \cd, X_n; \wh{\se}_n)} \leq \f{\al}{2}, \; \wh{\se}_n \leq \se \mid \se \ri \} \leq
\f{\al}{2} \eee for any $\se \in \mscr{S}$.   To show (\ref{LBB}),  making use of (\ref{RB2}), we have that \bee &  & \Pr \li \{  \f{ \sup_{\vse
\in \mscr{S}} f_n( X_1, \cd, X_n; \vse) }{\sup_{\vse \in \Se} f_n( X_1, \cd, X_n; \vse)} \leq \f{\al}{2}, \; \wh{\se}_n \geq \sup
\mscr{S} \mid \se \ri \}\\
&  & \leq \Pr \li \{  \f{ \sup_{\vse \in \mscr{S}} f_n( X_1, \cd, X_n; \vse) }{\sup_{\vse \in \Se} f_n( X_1, \cd, X_n; \vse)} \leq \f{\al}{2},
\; \wh{\se}_n \geq \se \mid \se \ri \}\\
&  & \leq \Pr \li \{  \f{ f_n( X_1, \cd, X_n; \se) }{\sup_{\vse \in \Se} f_n( X_1, \cd, X_n; \vse)} \leq \f{\al}{2}, \; \wh{\se}_n \geq \se \mid
\se \ri \}\\
&  & = \Pr \li \{  \f{ f_n( X_1, \cd, X_n; \se) }{f_n( X_1, \cd, X_n; \wh{\se}_n)} \leq \f{\al}{2}, \; \wh{\se}_n \geq \se \mid \se \ri \} \leq
\f{\al}{2} \eee for any $\se \in \mscr{S}$.  To show (\ref{LBC}), we use (\ref{RB3}) to conclude that \bee \Pr \li \{  \f{ \sup_{\vse \in
\mscr{S}} f_n( X_1, \cd, X_n; \vse) }{\sup_{\vse \in \Se} f_n( X_1, \cd, X_n; \vse)} \leq \f{\al}{2} \mid \se \ri \}  &  \leq & \Pr \li \{  \f{
f_n( X_1, \cd, X_n; \se) }{\sup_{\vse \in \Se} f_n( X_1, \cd,
X_n; \vse)} \leq \f{\al}{2} \mid \se \ri \} \\
&  =  & \Pr \li \{  \f{ f_n( X_1, \cd, X_n; \se) }{f_n( X_1, \cd, X_n; \wh{\se}_n)} \leq \f{\al}{2} \mid \se \ri \} \leq \al  \eee for any $\se
\in \mscr{S}$.  This completes the proof of the theorem.

\section{Proof of Recursive Formula for Multistage Sampling Without Replacement} \la{multsamwithout_app}

Consider a multistage sampling scheme of $s$ stages, based on sampling without replacement,  from a population of $N$ units, among which there
are $p N$ units having a certain attribute, where $p \in \{ \f{i}{N}: i = 1, \cd, N \}$.  Let the sample sizes be deterministic numbers $n_1 <
n_2 < \cd < n_s$.  For $\ell = 1, \cd, s$, let $K_\ell$ be the number of units having the attribute accumulated up to the $\ell$-th stage.  For
$\ell = 1, \cd, s$, let $\mscr{K}_\ell$ be a subset of the support of $K_\ell$.  In many applications, it is crucial to compute probabilities
like $\Pr \{ K_i \in \mscr{K}_i, \; i = 1, \cd, \ell \}, \; \ell = 1, \cd, s$.  For this purpose, we need to establish the following recursive
formula: \bel & & \Pr \{ K_i \in \mscr{K}_i, \; i
= 1, \cd, \ell; \; K_{\ell + 1} = k_{\ell + 1} \} \nonumber\\
&   & = \sum_{k_\ell \in \mscr{K}_\ell} \Pr \{ K_i \in \mscr{K}_i, \; i = 1, \cd, \ell - 1; \; K_\ell = k_\ell  \} \times \f{ \bi{p N -
k_\ell}{k_{\ell + 1} - k_\ell} \bi{N - p N - n_\ell +  k_\ell} { n_{\ell + 1} - n_\ell - k_{\ell + 1} + k_\ell } } { \bi{N - n_\ell}{n_{\ell +
1} - n_\ell} }  \la{recur18} \eel for $k_{\ell + 1} \in \mscr{K}_{\ell + 1}$ and $\ell = 1, \cd, s - 1$.  In the sequel, we shall provide a
rigorous justification for (\ref{recur18}) based on the notion of probability space.

Note that {\small \bee &  & \Pr \{ K_i \in \mscr{K}_i, \; i = 1, \cd, \ell; \; K_{\ell + 1} = k_{\ell + 1} \} = \sum_{k_\ell \in \mscr{K}_\ell}
\Pr \{ (
K_1, \cd, K_{\ell - 1}) \in \fra{D}, \; K_\ell = k_\ell, \; K_{\ell + 1} = k_{\ell + 1} \},\\
&  & \Pr \{ K_i \in \mscr{K}_i, \; i = 1, \cd, \ell - 1; \; K_\ell = k_\ell  \} = \Pr \{ ( K_1, \cd, K_{\ell - 1}) \in \fra{D}, \; K_\ell =
k_\ell \}, \eee} where $\fra{D} = \{ (k_1, \cd, k_{\ell - 1} ): k_i \in  \mscr{K}_i, \; i = 1, \cd, \ell - 1 \}$.  Hence, to show
(\ref{recur18}), it suffices to show \bel &   & \Pr \{ ( K_1, \cd, K_{\ell - 1}) \in \fra{D}, \; K_\ell =
k_\ell, \; K_{\ell + 1} = k_{\ell + 1} \} \nonumber\\
& = & \Pr \{ ( K_1, \cd, K_{\ell - 1}) \in \fra{D}, \; K_\ell = k_\ell \} \times \f{ \bi{p N - k_\ell}{k_{\ell + 1} - k_\ell} \bi{N - n_\ell - p
N + k_\ell} { n_{\ell + 1} - n_\ell - k_{\ell + 1} + k_\ell } } { \bi{N - n_\ell}{n_{\ell + 1} - n_\ell} }. \la{want8} \eel

We enumerate all units of the population as $U_1, U_2, \cd, U_N$.  Accordingly, $\{U_1, U_2, \cd, U_N \}$ can be partitioned as two exclusive
subsets $A$ and $B$ such that $A$ consists of $p N$ units having the attribute and that $B$ consists of $N - p N$ units without the attribute.
Let $\bs{\fra{U}}$ denote the set of all permutations of $U_1, U_2, \cd, U_N$. Define
\[
\bs{\fra{V}}_\ell \DEF \li \{ \overline{e_1 \; e_2 \; \cd \; e_N} \in \bs{\fra{U}}: \sum_{j = 1}^{n_i} \bb{I} (e_j) = k_i, \; i = 1, \cd, \ell
\ri \},
\]
where the character string $\overline{e_1 \; e_2 \; \cd \; e_N}$ is a permutation of $U_1, U_2, \cd, U_N$, and $\bb{I} (.)$ denote the indicator
function such that for $Z \in \{ U_1, U_2, \cd, U_N\}$,  $\bb{I} (Z) = 1$  if $Z$ has the attribute and $\bb{I} (Z) = 0$ otherwise.  We need to
figure out the number of permutations in $\bs{\fra{V}}_\ell$. Note that a permutation in $\bs{\fra{V}}_\ell$ can be constructed by the following
procedure.

\bed

\item The first step is to choose $k_1$ units from $A$ and $n_1 - k_1$ units from $B$.  Since $A$ has $p N$ units and $B$ has $N - p N$ units, there
are $\bi{p N}{k_1} \bi{N - p N}{n_1 - k_1}$ possible choices, among which each choice allows for $n_1!$ permutations.  Hence, there are $n_1!
\bi{p N}{k_1} \bi{N - p N}{n_1 - k_1}$  ways to obtain a character string of length $n_1$.

\item For $i = 2, \cd, \ell$, the $i$-th step is to choose $k_i - k_{i - 1}$ units from the remainder of $A$ and $n_i - n_{i-1} - (k_i - k_{i -
1})$ units from the remainder of $B$.  Since there are $p N - k_{i - 1}$ units remaining in $A$ and $N - n_{i -1} - p N + k_{i - 1}$ units
remaining in $B$ after the $(i-1)$-th step, there are $\bi{ p N - k_{i - 1}} {k_i - k_{i - 1}} \bi{N - n_{i -1} - p N + k_{i - 1}}{n_i - n_{i -
1} - k_i + k_{i - 1}}$ possible choices, among which each choice allows for $(n_i - n_{i - 1})!$ permutations.  Hence, there are $(n_i - n_{i -
1})! \bi{ p N - k_{i - 1}} {k_i - k_{i - 1}} \bi{N - n_{i -1} - p N + k_{i - 1}}{n_i - n_{i - 1} - k_i + k_{i - 1}}$ ways to obtain a character
string of length $n_i - n_{i - 1}$.

\item After the $\ell$-th step, the total number of units remaining in $A$ and $B$ is $N - n_\ell$. These $N - n_\ell$ units allows for $(N -
n_\ell)!$ permutations, which correspond to $(N - n_\ell)!$ ways of obtaining a character string of length $N - n_\ell$.

\item Connect all partial character strings obtained at all steps to make a complete character string like $\overline{e_1 \; e_2 \; \cd \;
e_N}$.

\eed

From the above procedure, it can be seen that there are $(N - n_\ell)! \; C_\ell (k_1, \cd, k_{\ell - 1}, k_\ell)$ permutations in
$\bs{\fra{V}}_\ell$, where
\[
C_\ell (k_1, \cd, k_{\ell - 1}, k_\ell) \DEF \prod_{i = 1}^\ell (n_i - n_{i - 1})! \bi{ p N - k_{i - 1}} {k_i - k_{i - 1}} \bi{N - n_{i -1} - p
N + k_{i - 1}}{n_i - n_{i - 1} - k_i + k_{i - 1}}.
\]
Notice that we have used $n_0 = k_0 = 0$ for purpose of simplifying notations.  Based on the above analysis, we have that there are $(N -
n_{\ell+1} )! \; C_{\ell + 1} (k_1, \cd, k_{\ell}, k_{\ell + 1})$ permutations in $\bs{\fra{V}}_{\ell + 1}$, where
\[
C_{\ell + 1} (k_1, \cd, k_{\ell}, k_{\ell + 1}) =  (n_{\ell +1} - n_\ell )! \; \bi{ p N - k_\ell } { k_{\ell + 1} - k_\ell } \bi{ N - n_\ell - p
N + k_\ell } { n_{\ell +1} - n_\ell - k_{\ell +1} + k_\ell } \times C(k_1, \cd, k_{\ell - 1}, k_\ell) \] and
\[
\bs{\fra{V}}_{\ell + 1} = \li \{ \overline{e_1 \; e_2 \; \cd \; e_N} \in \bs{\fra{U}}: \sum_{j = 1}^{n_i} \bb{I} (e_j) = k_i, \; i = 1, \cd,
\ell + 1 \ri \}.
\]
Note that each permutation of $U_1, U_2, \cd, U_N$ corresponds to a sample point in the sample space. Invoking the established fact there are
$(N - n_\ell)! \; C_\ell (k_1, \cd, k_{\ell - 1}, k_\ell)$ permutations in $\bs{\fra{V}}_\ell$, we have that the number of sample points in
event $\{ ( K_1, \cd, K_{\ell - 1}) \in \fra{D}, \; K_\ell = k_\ell \}$ is equal to $( N - n_\ell )! \sum_{(k_1, \cd, k_{\ell - 1}) \in \fra{D}
} C_\ell (k_1, \cd, k_{\ell - 1}, k_\ell)$.  Since each sample point of the sample space has the same probability $\f{1}{N!}$, we have \be
\la{how3} \Pr \{ ( K_1, \cd, K_{\ell - 1}) \in \fra{D}, \; K_\ell = k_\ell \} =  \f{( N - n_\ell )!}{N!}  \times \sum_{(k_1, \cd, k_{\ell - 1})
\in \fra{D} } C_\ell (k_1, \cd, k_{\ell - 1}, k_\ell). \ee Recalling  that the number of permutations in $\bs{\fra{V}}_{\ell + 1}$ is equal to
$(N - n_{\ell+1} )! \; C_{\ell + 1} (k_1, \cd, k_{\ell}, k_{\ell + 1})$, we have
 \bel &  & \Pr \{ ( K_1, \cd, K_{\ell - 1}) \in \fra{D},
\; K_\ell = k_\ell, \; K_{\ell + 1} = k_{\ell + 1} \} \nonumber\\
& = & \f{1}{N!} \times (N - n_{\ell + 1} )! \times \sum_{(k_1, \cd, k_{\ell - 1}) \in \fra{D}
} C_{\ell + 1} (k_1, \cd, k_{\ell}, k_{\ell + 1}) \nonumber\\
& = & \f{1}{N!} \times (N - n_{\ell + 1} )! \times \sum_{(k_1, \cd, k_{\ell - 1}) \in \fra{D}
} C_\ell (k_1, \cd, k_{\ell - 1}, k_\ell) \nonumber\\
&   & \times \bi{ p N - k_\ell } { k_{\ell + 1} - k_\ell } \times \bi{ N - n_\ell - p N + k_\ell } { n_{\ell +1} - n_\ell - k_{\ell +1} + k_\ell
} (n_{\ell +1} - n_\ell )! \nonumber \\
 & = & \f{1}{N!} \times ( N - n_\ell )! \times \sum_{(k_1, \cd,
k_{\ell - 1}) \in \fra{D}} C_\ell (k_1, \cd, k_{\ell - 1}, k_\ell) \nonumber \\
&   & \times \bi{ p N - k_\ell } { k_{\ell + 1} - k_\ell } \times \bi{ N - n_\ell - p N + k_\ell } { n_{\ell +1} - n_\ell - k_{\ell +1} + k_\ell
} (n_{\ell +1} - n_\ell )! \times \f{ (N - n_{\ell + 1} )!}{ ( N - n_\ell )! }. \la{how899} \eel  Combining ( \ref{how3} ) and ( \ref{how899} )
leads to (\ref{want8}) and consequently proves the recursive relationship (\ref{recur18}).

\sect{Proof of Theorem \ref{One_sided_Exact} }
\la{One_sided_Exact_Ap}

For arbitrary parametric values $\se_0 < \se_1$ in $\Se$, by the
assumption that $\bs{\varphi}_n$ converges in probability to $\se$,
we have that  $\Pr \{   \bs{\varphi}_n \geq \f{\se_0 + \se_1}{2}
\mid \se_0 \} \leq \Pr \{ | \bs{\varphi}_n - \se_0 | \geq \f{\se_1 -
\se_0}{2} \mid \se_0 \} \to 0$ and $\Pr \{ \bs{\varphi}_n \leq
\f{\se_0 + \se_1}{2} \mid \se_1 \} \leq \Pr \{ | \bs{\varphi}_n -
\se_1 | \geq \f{\se_1 - \se_0}{2} \mid \se_1 \} \to 0$ as $n \to
\iy$.  This shows that $\ovl{n}$ exists and is finite. By the
definition of the testing plan, we have {\small \bee \Pr \{
\tx{Accept} \; \mscr{H}_0 \mid \se \} & = & \sum_{\ell = 1}^s \Pr \{
\tx{Accept} \; \mscr{H}_0, \; \bs{l} = \ell \mid \se \}
 \leq  \sum_{\ell = 1}^s \Pr \{ \bs{D}_\ell = 1 \mid \se \}\\
& = & \sum_{\ell = 1}^s \Pr \li \{ \wh{\bs{\se}}_\ell \leq \udl{F}
(n_\ell, \se_0, \se_1, \al_0, \ba_1) \mid \se \ri \} \leq \sum_{\ell
= 1}^s \Pr \li \{ \wh{\bs{\se}}_\ell \leq \wh{F} (n_\ell, \se_1,
\ba_1) \mid \se \ri \}. \eee} Since $F_{\wh{\bs{\se}}_\ell} (z,
\se)$ is non-decreasing with respect to $z \in
I_{\wh{\bs{\se}}_\ell}$ for any given $\se \in {\Se}$, we have $\Pr
\{ \wh{\bs{\se}}_\ell \leq \wh{F} (n_\ell, \se_1, \ba_1) \mid \se \}
\leq \Pr \{ F_{\wh{\bs{\se}}_\ell} (\wh{\bs{\se}}_\ell, \se_1) \leq
\ba_1 \mid \se \}$ for $\ell = 1, \cd, s$. Since
$\wh{\bs{\se}}_\ell$ is a ULE of $\se$, by Lemma \ref{ULE_Basic}, we
have that $F_{\wh{\bs{\se}}_\ell} (z, \se) = \Pr \{
\wh{\bs{\se}}_\ell \leq z \mid \se \}$ is non-increasing with
respect to $\se$ no less than $z$. This implies that $\Pr \{
F_{\wh{\bs{\se}}_\ell} (\wh{\bs{\se}}_\ell, \se_1) \leq \ba_1 \mid
\se \} \leq \Pr \{ F_{\wh{\bs{\se}}_\ell} (\wh{\bs{\se}}_\ell, \se)
\leq \ba_1 \mid \se \}, \; \ell = 1, \cd, s$ for $\se \in {\Se}$ no
less than $\se_1$.  Therefore, $\Pr \{ \tx{Accept} \; \mscr{H}_0
\mid \se \} \leq \sum_{\ell = 1}^s \Pr \{ \bs{D}_\ell = 1 \mid \se
\} \leq \sum_{\ell = 1}^s \Pr \{ F_{\wh{\bs{\se}}_\ell}
(\wh{\bs{\se}}_\ell, \se )  \leq \ba_1 \mid \se \}  \leq s \ba_1$
for $\se \in {\Se}$ no less than $\se_1$, where the last inequality
follows from Lemma \ref{ProbTrans}.  By a similar method, we can
show that $\Pr \{ \tx{Reject} \; \mscr{H}_0 \mid \se \} \leq
\sum_{\ell = 1}^s \Pr \{ \bs{D}_\ell = 2 \mid \se \} \leq s \al_0$
for $\se \in {\Se}$ no greater than $\se_0$. By the definition of
the testing plan and the assumption that the likelihood ratio is
monotonically increasing with respect to $\wh{\bs{\se}}_\ell$, we
have that the test procedure is a generalized SPRT.  Hence, the
monotonicity of $\Pr \{ \tx{Accept} \; \mscr{H}_0 \mid \se \}$ with
respect to $\se$ is established. This concludes the proof of the
theorem.

\sect{Proof of Theorem \ref{Normal_Mean_unknown_variance} } \la{App_Normal_Mean_unknown_variance}

We need some preliminary results.

\beL \la{lemch1} For any $\de \in (0, 1)$, $\f{ t_{n, \de} }{\sq{n}}
$ is monotonically decreasing to $0$ as $n$ increases from $2$ to
$\iy$. \eeL

\bpf

For simplicity of notations, let $\psi (n) = \f{ t_{n, \de}
}{\sq{n}}$. Then, $\de = \Pr \{ \f{|U|}{ \sq{Z \sh n} } > t_{n, \de}
\}  = \Pr \{ \f{|U|}{ \sq{Z} } > \psi (n) \}$, where $U$ and $Z$ are
independent random variables such that $U$ is a Gaussian variable
with zero mean and unit variance and that $Z$ is a chi-squared
variable of $n$ degrees of freedom. Since $\f{U}{ \sq{Z \sh n} }$
possesses a Student's $t$-distribution of $n$ degrees of freedom,
its mean and variance are, respectively, $0$ and $\f{n}{n - 2}$.
Accordingly, the mean and variance of $\f{U}{ \sq{Z} }$ are,
respectively, $0$ and $\f{1}{n - 2}$.  By Chebyshev's inequality,
{\small $\Pr \{ \f{|U|}{ \sq{Z} } > \psi \} \leq \f{1}{(n - 2) [
\psi (n) ]^2}$}, leading to $\de < \f{1}{(n - 2) [ \psi (n) ]^2}$,
i.e., $\psi (n) < \f{1}{\sq{(n - 2) \de} } \to 0$ as $n \to \iy$.
This proves $\lim_{n \to \iy} \f{ t_{n, \de} }{\sq{n}} = 0$.

To show the monotonicity, it suffices to show that, for any fixed $t
> 0$,  $\Pr \{ |U| \sh \sq{Z}  > t \}$ decreases monotonically with
respect to $n$.  Let $V_1, \cd, V_n, V_{n + 1}$ be i.i.d. Gaussian
random variables which have zero mean, unity variance  and are
independent with $U$. Then, {\small $\Pr \{ |U| \sh \sq{Z}
> t \} = \Pr \{ |U| \sh \sq{ \sum_{i = 1}^n V_i^2  } >
t \}$}.   In view of  {\small $\Pr \{ |U| \sh \sq{\sum_{i = 1}^n
V_i^2 } > t  \} > \Pr \{ |U| \sh \sq{ \sum_{i = 1}^{n + 1} V_i^2 }
> t  \}$}  and {\small $ \Pr  \{ |U| \sh \sq{ \sum_{i = 1}^n V_i^2 }
> \psi(n)  \} = \Pr \{ |U| \sh \sq{ \sum_{i = 1}^{n+ 1} V_i^2
}  > \psi(n + 1) \} = \de$}, we have {\small $\Pr \{ |U| \sh \sq{
\sum_{i = 1}^{n+ 1} V_i^2 }  > \psi(n + 1) \}
> \Pr \{ |U| \sh \sq{ \sum_{i = 1}^{n+ 1} V_i^2 }  >
\psi(n) \}$,} which implies $\psi(n + 1) < \psi (n)$.  This
completes the proof of the lemma.

\epf

\beL \la{lemch} {\small $\lim_{\de \to 0} \f{ \mcal{Z}_\de } { \sq{
2 \ln \f{1}{\de} } } = 1$.} \eeL

\bpf

For simplicity of notations, we abbreviate $\mcal{Z}_\de$ as $z$
when this can be done without introducing confusion. By virtue of
the well-known inequality {\small $1 - \Phi(z) < \f{1}{\sq{2 \pi}}
\exp \li ( - \f{z^2}{2} \ri ) \li ( \f{1}{z} \ri )$}, we have
{\small $\de < \f{1}{\sq{2 \pi}} \exp \li ( - \f{z^2}{2} \ri ) \li (
\f{1}{z} \ri )$}, or equivalently, {\small $\f{2 \ln
\f{1}{\de}}{z^2}
> \f{2 \ln (\sq{2 \pi} z)}{z^2}  + 1$}, which implies
{\small $\liminf_{z \to \iy} \f{2 \ln \f{1}{\de}}{z^2} \geq
1$} and, consequently, {\small $\limsup_{\de \to 0} \f{ \mcal{Z}_\de
} { \sq{ 2 \ln \f{1}{\de} } } \leq 1$}.  On the other hand, making
use of the well-known inequality {\small $\f{1}{\sq{2 \pi}} \exp \li
( - \f{z^2}{2} \ri ) \li ( \f{1}{z} - \f{1}{z^3} \ri ) < 1 -
\Phi(z)$}, we have {\small $\de
> \f{1}{\sq{2 \pi}} \exp \li ( - \f{z^2}{2} \ri ) \li ( \f{1}{z} \ri
) \li (1 - \f{1}{z^2} \ri )$}, which implies {\small $\f{2 \ln
\f{1}{\de}}{z^2} < \f{2}{z^2} \ln \li (\f{\sq{2 \pi} z^3}{z^2 - 1}
\ri )  + 1$} and thus {\small $\liminf_{\de \to 0} \f{ \mcal{Z}_\de
} { \sq{ 2 \ln \f{1}{\de} } } \geq 1$}.  This establishes {\small
$\lim_{\de \to 0} \f{ \mcal{Z}_\de } { \sq{ 2 \ln \f{1}{\de} } } =
1$.}

\epf

\beL \la{lemCH} Let $X$ be a chi-squared random variable with $n$
degrees of freedom. Then, $\Pr \{  X \geq n(1 + \ka) \} \leq [ (1 +
\ka)  e^{-\ka }  ]^{ \f{n}{2} }$ for any $\ka > 0$ and $\Pr \{ X
\leq n(1 - \ka) \} \leq  [ (1 - \ka) e^{\ka} ]^{ \f{n}{2} }$ for $0
< \ka < 1$.  \eeL

\bpf For simplicity of notations,  let $c = n(1 + \ka)$.  Then, \bee
\Pr \li \{ X \geq c \ri \} & \leq & \inf_{\ro  > 0} \bb{E} \li [
e^{\ro  (X - c)} \ri ] =  \inf_{\ro  > 0} \int_0^\iy \f{ 1 } { 2^{
\f{n}{2} }  \Ga \li( \f{n}{2} \ri) } x^{\f{n}{2} - 1}e^{-
\f{x}{2}} e^{\ro  (x - c)} dx\\
& = & \inf_{\ro  > 0} e^{- \ro  c} ( 1 - 2 \ro  )^{- \f{n}{2} }
\int_0^\iy \f{ 1 } { 2^{n}  \Ga \li( \f{n}{2} \ri) } y^{\f{n}{2} -
1}e^{- \f{y}{2}} dy = \inf_{\ro  > 0} e^{- \ro  c} ( 1 - 2 \ro  )^{-
\f{n}{2} }, \eee where we have introduced a change of variable $\li(
\f{1}{2} - \ro  \ri ) x = \f{y}{2}$ in the integration.  Note that
{\small $\f{ d  } {d \ro }  [  e^{- \ro  c} ( 1 - 2 \ro  )^{ -
\f{n}{2} }  ] = ( \f{n}{1 - 2 \ro } - c ) e^{- \ro  c} ( 1 - 2 \ro
)^{ - \f{n}{2} }$}, which equals $0$ for $\ro  = \f{c - n}{2c} > 0$.
Therefore, \bee \Pr \li \{ X \geq n(1 + \ka) \ri \} & \leq & \exp
\li( - \f{c - n}{2c} c \ri ) \li ( \f{1}{1 - 2 \f{c - n}{2c} } \ri
)^{ \f{n}{2} } = \li ( \f{ 1 + \ka } { e^\ka } \ri )^{ \f{n}{2} }
\eee for any $\ka > 0$. Similarly, $\Pr \li \{ X \leq n(1 - \ka) \ri
\} \leq \li ( \f{ 1 - \ka } { e^{-\ka} } \ri )^{ \f{n}{2} }$ for $0
< \ka < 1$.  This completes the proof of the lemma. \epf

The following result is due to Wallace \cite{Wallace}.

\beL \la{Wallace} Let $F(t)$ be Student's t-distribution of $n$
degrees of freedom. Let $x(t)$ be the root of equation $\Phi(x) =
F(t)$ with respect to $x$. Then, {\small $\sqrt{ \li ( n - \f{1}{2}
\ri ) \ln \li ( 1+ \f{t^2}{n} \ri ) } \leq x(t) \leq \sqrt{ n \ln
\li ( 1+ \f{t^2}{n} \ri ) }$} for any $t > 0$.  \eeL

\beL \la{keyy} For any $\ep > 0$, there exists a number $\ze^* > 0$
such that {\small $\li | \f{t_{n, \bs{\al}} - t_{n, \bs{\ba}} }{\sq{
n} } \ri | < \ep$} for any $\ze \in (0, \ze^*)$ and all {\small $n
\geq \ka(\ze, \vro) = \min \li \{ \li \lf  \vro \;
\mcal{Z}_{\sq{\bs{\al}}}^2 \ri \rf, \; \li \lf  \vro \;
\mcal{Z}_{\sq{\bs{\ba}}}^2  \ri \rf \ri \} > 1$}, where $\ze^*$ is
independent of $n$ and $\vro > 0$. \eeL

\bpf  Define {\small \[ h (\ze, n) = \li [  \ln \li ( 1 + \f{t_{n,
\bs{\al}}^2}{n} \ri ) \ri ] \li [ \ln \li ( 1 + \f{t_{n,
\bs{\ba}}^2}{n} \ri ) \ri ]^{-1}
\]}
for $n \geq \ka(\ze, \vro)$.  We shall first show that $h (\ze, n)$
tends to $1$ uniformly for $n \geq \ka(\ze, \vro)$ as $\ze \to 0$.
Applying Lemma \ref{Wallace}, we have  {\small \be \la{equi}
 \f{  \mcal{Z}_{\bs{\al}}^2
} { n } \leq \ln \li ( 1+ \f{t_{n, \bs{\al}}^2}{n} \ri ) \leq \f{
\mcal{Z}_{\bs{\al}}^2  } {  n - \f{1}{2}  }, \qqu \f{
\mcal{Z}_{\bs{\ba}}^2 } { n } \leq \ln \li ( 1+ \f{t_{n,
\bs{\ba}}^2}{n} \ri ) \leq \f{ \mcal{Z}_{\bs{\ba}}^2 } { n -
\f{1}{2} } \ee} and thus {\small
\[
\li (1 - \f{1}{2  \ka(\ze, \vro) } \ri )  \li ( \f{
\mcal{Z}_{\bs{\al}} }{ \mcal{Z}_{\bs{\ba}} } \ri )^2 < \f{ n -
\f{1}{2} } { n  } \li ( \f{ \mcal{Z}_{\bs{\al}} }{
\mcal{Z}_{\bs{\ba}} } \ri )^2 \leq h (\ze, n) \leq \f{ n  } { n -
\f{1}{2} } \li ( \f{ \mcal{Z}_{\bs{\al}} }{ \mcal{Z}_{\bs{\ba}} }
\ri )^2  < \li (1 + \f{1}{2  \ka(\ze, \vro) - 1 } \ri )  \li ( \f{
\mcal{Z}_{\bs{\al}} }{ \mcal{Z}_{\bs{\ba}} } \ri )^2  \]} for $n
\geq \ka(\ze, \vro)$. By Lemma \ref{lemch}, we have {\small
\[ \lim_{\ze \to 0} \f{ \mcal{Z}_{\bs{\al}} }{ \mcal{Z}_{\bs{\ba}}} =
\lim_{\ze \to 0} \li [ \f{ \mcal{Z}_{\bs{\al}} }{ \sq{ 2 \ln
\f{1}{\bs{\al}} } } \times \f{ \sq{ 2 \ln \f{1}{\bs{\al}} } } { \sq{
2 \ln \f{1}{\bs{\ba}} } } \li \sh \f{ \mcal{Z}_{\bs{\ba}} }{ \sq{ 2
\ln \f{1}{\bs{\ba}} } } \ri. \ri ] = 1.
\]}
It follows that $h (\ze, n)$ tends to $1$ uniformly for $n \geq
\ka(\ze, \vro)$ as $\ze \to 0$.  By virtue of (\ref{equi}), we have
\[
\ln \li ( 1+ \f{t_{n, \bs{\al}}^2}{n} \ri ) \leq \f{
\mcal{Z}_{\bs{\al}}^2  } {  n - \f{1}{2}  } \leq \f{
\mcal{Z}_{\bs{\al}}^2 } { \ka(\ze, \vro) - \f{1}{2}  }  \to
\f{2}{\vro}
\]
and
\[
\ln \li ( 1+ \f{t_{n, \bs{\ba}}^2}{n} \ri ) \leq \f{
\mcal{Z}_{\bs{\ba}}^2  } {  n - \f{1}{2}  } \leq \f{
\mcal{Z}_{\bs{\ba}}^2 } { \ka(\ze, \vro) - \f{1}{2}  }  \to
\f{2}{\vro}
\]
uniformly for $n \geq \ka(\ze, \vro)$ as $\ze \to 0$.  Therefore,
both {\small $\f{t_{n, \bs{\al}}^2}{n}$} and {\small $\f{t_{n,
\bs{\ba}}^2}{n}$} are bounded uniformly for all $n \geq \ka(\ze,
\vro)$ and any $\ze \in (0, 1)$.  By virtue of this result and
recalling that $h (\ze, n)$ tends to $1$ uniformly for $n \geq
\ka(\ze, \vro)$ as $\ze \to 0$,  we have that {\small $\ln ( 1+
\f{t_{n, \bs{\al}}^2} {n} ) - \ln ( 1+ \f{t_{n, \bs{\ba}}^2}{n} )$}
tends to $0$ and thus {\small $ \f{t_{n, \bs{\al}} - t_{n, \bs{\ba}}
}{\sq{ n} }$} tends to $0$ uniformly for $n \geq \ka(\ze, \vro)$ as
$\ze \to 0$.  This completes the proof of the lemma.

\epf

\beL \la{lemgood9}

For any $\vDe > 0$, $\sum_{n = \ka(\ze, \vro) + 1}^\iy \Pr \{ | \f{
\ovl{X}_n } { \wt{\si}_n } - \se | \geq \vDe \mid \se \} \to 0$ as
$\ze \to 0$, where {\small $\ka(\ze, \vro) = \min \li \{ \li \lf
\vro \; \mcal{Z}_{\sq{\bs{\al}}}^2 \ri \rf, \; \li \lf  \vro \;
\mcal{Z}_{\sq{\bs{\ba}}}^2  \ri \rf \ri \}$}.  \eeL

\bpf

We shall first show that {\small $\sum_{n = \ka(\ze, \vro) + 1}^\iy
\Pr \li \{ \f{ \ovl{X}_n } { \wt{\si}_n } \leq  \se - \vDe \mid \se
\ri \} \to 0$} as $\ze \to 0$ by considering two cases: (i) $\se
\geq \vDe$; (ii) $\se < \vDe$.

In the case of $\se \geq \vDe$, let $\eta$ be a positive number such
that $(1 + \eta) (\se - \vDe) < \se$. Then, \bel \Pr \li \{ \f{
\ovl{X}_n } { \wt{\si}_n } \leq \se - \vDe \mid \se \ri \} & \leq &
\Pr \li \{ \f{ \ovl{X}_n } { \wt{\si}_n } \leq \se - \vDe, \;
\wt{\si}_n \leq (1 + \eta) \si \mid \se \ri \} + \Pr \{ \wt{\si}_n
> (1 + \eta) \si \mid \se \} \nonumber\\
& \leq &  \Pr \{ \ovl{X}_n \leq (1 + \eta) \si (\se - \vDe) \mid \se
\} + \Pr \{ \wt{\si}_n > (1 + \eta) \si \mid \se \} \nonumber\\
& = & \Pr \{  U \geq \sq{n} [ (1 + \eta) \vDe - \eta \se ] \} + \Pr
\{ \chi_{n-1}^2 > n (1 + \eta)^2 \} \nonumber\\
& < & \Pr \{  U \geq \sq{n} [ (1 + \eta) \vDe - \eta \se ] \} + \Pr
\{ \chi_{n-1}^2 > (n - 1) (1 + \eta) \}, \la{com1} \eel where $U$ is
a Gaussian random variable with zero mean and unit variance and
$\chi_{n-1}^2$ is a chi-square variable of $n - 1$ degrees of
freedom.  By the choice of $\eta$, we have $ (1 + \eta) \vDe - \eta
\se > 0$ as a consequence of $(1 + \eta) (\se - \vDe) < \se$. Hence,
\be \la{com2} \Pr \li \{  U \geq \sq{n} [ (1 + \eta) \vDe - \eta \se
] \ri \} < \exp \li ( - \f{n}{2}  [ (1 + \eta) \vDe - \eta \se ]^2
\ri ). \ee On the other hand, by Lemma \ref{lemCH}, we have \be
\la{com3} \Pr \{ \chi_{n-1}^2 > (n - 1) (1 + \eta) \} \leq [ (1 +
\eta) e^{- \eta} ]^{(n - 1) \sh 2}. \ee Combining (\ref{com1}),
(\ref{com2}) and (\ref{com3}) yields
\[
\sum_{n = \ka(\ze, \vro) + 1}^\iy \Pr \li \{ \f{ \ovl{X}_n } {
\wt{\si}_n } \leq  \se - \vDe \mid \se \ri \} < \sum_{n = \ka(\ze,
\vro) + 1}^\iy \li [  \exp \li ( - \f{n}{2}  [ (1 + \eta) \vDe -
\eta \se ]^2 \ri ) +   [ (1 + \eta) e^{- \eta} ]^{(n - 1) \sh 2} \ri
],
\]
where the right side tends to $0$ as $\ze \to 0$ because $\ka(\ze,
\vro) \to \iy$ as $\ze \to 0$.

In the case of $\se < \vDe$, let $\eta \in (0, 1)$ be a number such
that $(1 - \eta) (\se - \vDe) < \se$. Then, \bel \Pr \li \{
\f{\ovl{X}_n } { \wt{\si}_n } \leq \se - \vDe \mid \se \ri \} & \leq
& \Pr \li \{ \f{ \ovl{X}_n } { \wt{\si}_n } \leq \se - \vDe, \;
\wt{\si}_n \geq (1 - \eta) \si \mid \se \ri \} + \Pr \{ \wt{\si}_n <
(1 - \eta) \si \mid \se \} \nonumber\\
& \leq &  \Pr \{ \ovl{X}_n \leq (1 - \eta) \si (\se - \vDe) \mid \se
\} + \Pr \{ \wt{\si}_n < (1 - \eta) \si \mid \se \} \nonumber\\
& = & \Pr \{  U \geq \sq{n} [\eta \se + (1 - \eta) \vDe] \} + \Pr \{
\chi_{n-1}^2 < n (1 - \eta)^2 \}. \la{com11} \eel  By the choice of
$\eta$, we have $\eta \se + (1 - \eta) \vDe > 0$ as a consequence of
$(1 - \eta) (\se - \vDe) < \se$. Hence, \be \la{com22} \Pr \{  U
\geq \sq{n} [ \eta \se + (1 - \eta) \vDe  ] \} < \exp \li ( -
\f{n}{2}  [ \eta \se + (1 - \eta) \vDe ]^2 \ri ). \ee For small
enough $\ze > 0$, we have $n > \ka(\ze, \vro) > \f{1}{\eta}$ and
thus \be \la{com33} \Pr \{ \chi_{n-1}^2 < n (1 - \eta)^2 \} < \Pr \{
\chi_{n-1}^2 < (n-1) (1 - \eta) \}  \leq [ (1 - \eta) e^{ \eta}
]^{(n - 1) \sh 2}, \ee where the last inequality follows from Lemma
\ref{lemCH}.  Combining (\ref{com11}), (\ref{com22}) and
(\ref{com33}) yields
\[
\sum_{n = \ka(\ze, \vro) + 1}^\iy \Pr \li \{ \f{ \ovl{X}_n } {
\wt{\si}_n } \leq  \se - \vDe \mid \se \ri \} < \sum_{n = \ka(\ze,
\vro) + 1}^\iy \li [  \exp \li ( - \f{n}{2}  [ \eta \se + (1 - \eta)
\vDe  ]^2 \ri ) +   [ (1 - \eta) e^{\eta} ]^{(n - 1) \sh 2} \ri ],
\]
where the right side tends to $0$ as $\ze \to 0$ because $\ka(\ze,
\vro) \to \iy$ as $\ze \to 0$.    This proves that $\sum_{n =
\ka(\ze, \vro) + 1}^\iy \Pr \{ \f{ \ovl{X}_n } { \wt{\si}_n } \leq
\se - \vDe \mid \se \} \to 0$ as $\ze \to 0$. In a similar manner,
we can show that $\sum_{n = \ka(\ze, \vro) + 1}^\iy \Pr \{ \f{
\ovl{X}_n } { \wt{\si}_n } \geq \se + \vDe \mid \se \} \to 0$ as
$\ze \to 0$. This concludes the proof of the lemma. \epf

\beL \la{make8}

Let $\bs{\de} = O(\ze) \in (0, 1)$. If $\ze > 0$ is sufficiently
small, then
 \[ \f{1}{|\se|} \li ( \f{
t_{n, \bs{\de}} } { \sq{n} } - \f{ t_{n, \sq{\bs{\de}} } }{ \sq{n} }
\ri ) >  \exp \li ( \f{ \ln \f{1}{\bs{\de} } }{ 4 n } \ri ) > 1 \]
for $2 \leq n < \li \lf \mcal{Z}_{\sq{\bs{\de}}}^2 \ri \rf$, where
$0 < \vro < \f{1}{4 ( 1 + |\se|)^2}$.

\eeL

\bpf  From Wallace's inequality restated in Lemma \ref{Wallace},  we
have
\[
\sq{ \exp \li (  \f{ \mcal{Z}_\de^2 }{ n }  \ri )  - 1 } \leq \f{
t_{n, \de} } { \sq{n} } \leq \sq{ \exp \li (  \f{ \mcal{Z}_\de^2 }{
n - \f{1}{2} }  \ri )  - 1 }, \qqu \fa \de \in (0, 1)
\]
and thus \[ \f{1}{|\se|} \li ( \f{ t_{n, \bs{\de}} } { \sq{n} } -
\f{ t_{n, \sq{\bs{\de}} } }{ \sq{n} } \ri ) > \f{1}{|\se|} \li [
\sq{ \exp \li ( \f{ \mcal{Z}_{\bs{\de}}^2 }{ n }  \ri )  - 1 } -
\sq{ \exp \li ( \f{ \mcal{Z}_{\sq{\bs{\de}}}^2 }{ n - \f{1}{2} } \ri
) - 1 } \ri ].
\]
Therefore, to show the lemma, it suffices to show that \be
\la{make2} \f{1}{|\se|} \li [ \sq{ \exp \li ( \f{
\mcal{Z}_{\bs{\de}}^2 }{ n }  \ri )  - 1 } - \sq{ \exp \li ( \f{
\mcal{Z}_{\sq{\bs{\de}}}^2 }{ n - \f{1}{2} } \ri ) - 1 } \ri ] >
\exp \li ( \f{ \ln \f{1}{\bs{\de} } }{ 4 n } \ri ) > 1 \ee for $2
\leq n < \li \lf \mcal{Z}_{\sq{\bs{\de}}}^2 \ri \rf$ if $\ze > 0$ is
small enough.  By Lemma \ref{lemch}, for small enough $\ze
> 0$, we have $\ln \f{1}{\sq{\bs{\de}} } < \f{2}{3}
\mcal{Z}_{\sq{\bs{\de}}}^2$ and thus \bee \f{ \exp \li ( \f{
\mcal{Z}_{\sq{\bs{\de}}}^2 }{ n } \ri ) }{\exp \li ( \f{ \ln
\f{1}{\sq{\bs{\de}} } }{  n } \ri ) } - \f{1}{ \exp \li ( \f{ \ln
\f{1}{\sq{\bs{\de}} } }{  n  } \ri )} &
> & \exp \li ( \f{ \mcal{Z}_{\sq{\bs{\de}}}^2 }{ 3 n } \ri )  -
1  >  \exp \li ( \f{ \mcal{Z}_{\sq{\bs{\de}}}^2 }{ 3
(n + 1) } \ri )  - 1 \\
& \geq & \exp \li ( \f{ 1 }{ 3 \vro } \ri )  - 1 > \f{ 1 }{ 3 \vro }
> \f{ 4 ( 1 + |\se|)^2 }{ 3 } > 1 \eee
for $2 \leq n < \li \lf \mcal{Z}_{\sq{\bs{\de}}}^2 \ri \rf$.  Hence,
\[
\sq{ \exp \li ( \f{ \mcal{Z}_{\sq{\bs{\de}}}^2 }{ n - \f{1}{2} } \ri
) - 1 } > \sq{ \exp \li ( \f{ \mcal{Z}_{\sq{\bs{\de}}}^2 }{ n } \ri
) - 1 } > \exp \li ( \f{ \ln \f{1}{\bs{\de} } }{ 4 n } \ri ) > 1
\]
for $2 \leq n < \li \lf \mcal{Z}_{\sq{\bs{\de}}}^2 \ri \rf$ if $\ze$
is small enough. Therefore, to guarantee (\ref{make2}), it suffices
to make $\ze$ small enough and ensure that
\[
\sq{ \exp \li ( \f{ \mcal{Z}_{\bs{\de}}^2 }{ n }  \ri )  - 1 }
> (1 + |\se|)  \sq{ \exp \li ( \f{ \mcal{Z}_{\sq{\bs{\de}}}^2 }{ n - \f{1}{2} } \ri ) - 1 }.
\]
By Lemma \ref{lemch}, we have $\lim_{\ze \to 0} \f{
\mcal{Z}_{\bs{\de}}^2 } { \mcal{Z}_{\sq{\bs{\de}}}^2 } = 2$.  This
implies that,  if $\ze
> 0$ is sufficiently small, then $\f{ \mcal{Z}_{\bs{\de}}^2  }
 { \mcal{Z}_{\sq{\bs{\de}}}^2 } > \f{5}{3}$,  and consequently,
\[
\f{ \mcal{Z}_{\bs{\de}}^2 }{ n }  - \f{ \mcal{Z}_{\sq{\bs{\de}}}^2
}{ n - \f{1}{2} } = \f{ \mcal{Z}_{\sq{\bs{\de}}}^2 }{ n - \f{1}{2} }
\li ( \f{ n - \f{1}{2} } { n } \f{ \mcal{Z}_{\bs{\de}}^2  } {
\mcal{Z}_{\sq{\bs{\de}}}^2 }   - 1 \ri )
> \f{ 1 } { \vro } \li ( \f{ 2 - \f{1}{2} } { 2 } \times
\f{5}{3} - 1 \ri ) = \f{1}{4 \vro}
\]
for $2 \leq n < \li \lf \mcal{Z}_{\sq{\bs{\de}}}^2 \ri \rf$.  Hence,
\[
\f{  \exp \li ( \f{ \mcal{Z}_{\bs{\de}}^2 }{ n }  \ri )  - 1 }{ \exp
\li ( \f{ \mcal{Z}_{\sq{\bs{\de}}}^2 }{ n - \f{1}{2} } \ri ) - 1 } >
\f{ \exp \li ( \f{ \mcal{Z}_{\bs{\de}}^2 }{ n }  \ri ) - 1 }{ \exp
\li ( \f{ \mcal{Z}_{\sq{\bs{\de}}}^2 }{ n - \f{1}{2} } \ri ) }
> \f{  \exp \li ( \f{ \mcal{Z}_{\bs{\de}}^2 }{ n } \ri ) }{ \exp
\li ( \f{ \mcal{Z}_{\sq{\bs{\de}}}^2 }{ n - \f{1}{2} } \ri ) } - 1
>  \exp \li (  \f{1}{4 \vro} \ri ) - 1 > \f{1}{4 \vro} > ( 1 +
|\se|)^2
\]
for $2 \leq n < \li \lf \mcal{Z}_{\sq{\bs{\de}}}^2 \ri \rf$, and
consequently (\ref{make2}) is ensured if $\ze > 0$ is small enough.
This completes the proof of the lemma.

\epf

\beL \la{lem11m}  Let $\se^\prime < \se^{\prime \prime}$ and
$\ka(\ze, \vro) = \min \li \{ \li \lf  \vro \;
\mcal{Z}_{\sq{\bs{\al}}}^2  \ri \rf, \; \li \lf  \vro \;
\mcal{Z}_{\sq{\bs{\ba}}}^2 \ri \rf \ri \}$.  Then, {\small \be
\la{goal1} \lim_{\ze \to 0} \li [ \sum_{n = 2}^{\ka(\ze, \vro)} \Pr
\li \{ \f{ \ovl{X}_n } { \wt{\si}_n } \leq \se^{\prime \prime} - \f{
t_{n-1, \bs{\ba}} } { \sq{n - 1} } \mid \se \ri \} + \sum_{n =
\ka(\ze, \vro) + 1}^\iy \Pr \li \{ \f{ \ovl{X}_n } { \wt{\si}_n }
\leq \f{ \se^\prime + \se^{\prime \prime} }{2} + \f{ t_{n-1,
\bs{\al}} - t_{n-1, \bs{\ba}} } { 2 \sq{n - 1} } \mid \se \ri \} \ri
] = 0 \ee} for $\se \geq \se^{\prime \prime}$ provided that $0 <
\vro < \f{1}{6 (1 + |\se|)^2}$.  Similarly, {\small \be \la{goal2}
\lim_{\ze \to 0} \li [ \sum_{n = 2}^{\ka(\ze, \vro)} \Pr \li \{ \f{
\ovl{X}_n } { \wt{\si}_n } \geq \se^{\prime} + \f{ t_{n-1, \bs{\al}}
} { \sq{n - 1} } \mid \se \ri \} + \sum_{n = \ka(\ze, \vro) + 1}^\iy
\Pr \li \{ \f{ \ovl{X}_n } { \wt{\si}_n } \geq \f{ \se^\prime +
\se^{\prime \prime} }{2} + \f{ t_{n-1, \bs{\al}} - t_{n-1, \bs{\ba}}
} { 2 \sq{n - 1} } \mid \se \ri \} \ri ] = 0 \ee} for $\se \leq
\se^{\prime}$ provided that $0 < \vro < \f{1}{6 (1 + |\se|)^2}$.

\eeL

\bpf  Without loss of generality, assume that $\ze$ is sufficiently
small so that $\ka(\ze, \vro)$ is greater than $2$.  We shall first
show that \be \la{shw1}
 \lim_{\ze \to 0}
\sum_{n = 2}^{\ka(\ze, \vro)} \Pr \li \{ \f{ \ovl{X}_n } {
\wt{\si}_n } \leq \se^{\prime \prime} - \f{ t_{n-1, \bs{\ba}} } {
\sq{n - 1} } \mid \se \ri \} = 0 \ee for $\se \geq \se^{\prime
\prime}$.  Obviously, $\lim_{\ze \to 0} \Pr \li \{ \f{ \ovl{X}_n } {
\wt{\si}_n } \leq \se^{\prime \prime} - \f{ t_{n-1, \bs{\ba}} } {
\sq{n - 1} } \mid \se \ri \} = 0$ for $n = 2$ and $\se \geq
\se^{\prime \prime}$.  Hence, to show (\ref{shw1}), it remains to
show \be \la{shw2}
 \lim_{\ze \to 0} \sum_{n = 3}^{\ka(\ze, \vro)}
\Pr \li \{ \f{ \ovl{X}_n } { \wt{\si}_n } \leq \se^{\prime \prime} -
\f{ t_{n-1, \bs{\ba}} } { \sq{n - 1} } \mid \se \ri \} = 0 \ee for
$\se \geq \se^{\prime \prime}$.  We shall show (\ref{shw2}) by
considering three cases: (i) $\se = 0$; (ii) $\se < 0$;  (iii) $\se
> 0$.

In the case of $\se = 0 \geq \se^{\prime \prime}$, we have
\[ \sum_{n = 3}^{\ka(\ze, \vro)} \Pr \li \{ \f{ \ovl{X}_n } {
\wt{\si}_n } \leq \se^{\prime \prime} - \f{ t_{n-1, \bs{\ba}} } {
\sq{n - 1} } \mid \se \ri \} \leq \sum_{n = 3}^{\ka(\ze, \vro)} \Pr
\li \{  \f{ \sq{n} \; \ovl{X}_n } { \wh{\si}_n } \leq - t_{n-1,
\bs{\ba}} \mid \se \ri \} < \ka(\ze, \vro) \; \bs{\ba}.
\]
Noting that \[ \ka(\ze, \vro) \; \bs{\ba} \leq \vro \;
\mcal{Z}_{\sq{\bs{\ba}}}^2  \times \bs{\ba}  = \vro \times \f{
\mcal{Z}_{\sq{\bs{\ba}}}^2 }{ 2 \ln \f{1}{ \sq{\bs{\ba}} } } \times
\bs{\ba} \times  { 2 \ln \f{1}{ \sq{\bs{\ba}} } } \to 0
\]
as $\ze \to 0$, we have that (\ref{shw2}) is true for the case of
$\se = 0 \geq \se^{\prime \prime}$.  Hence, it remains to show that
(\ref{shw2}) is true for the cases of $\se < 0 $ and $\se > 0$.  Let
\[
\vDe_n = \se \sq{n - 1} \li ( 1 - \f{\si}{\wt{\si}_n} \ri ) +
t_{n-1, \sq{\bs{\ba}} }-
 t_{n-1, \bs{\ba}}, \qqu n = 3, 4, \cd.
\] Note that {\small \bee  \sum_{n = 3}^{\ka(\ze,
\vro)} \Pr \li \{ \f{ \ovl{X}_n } { \wt{\si}_n } \leq \se^{\prime
\prime} - \f{ t_{n-1, \bs{\ba}} } { \sq{n - 1} } \mid \se \ri \} & =
& \sum_{n = 3}^{\ka(\ze, \vro)} \Pr \li \{  \f{ \sq{n} \; \ovl{X}_n
} { \wh{\si}_n } \leq - t_{n-1, \bs{\ba}} +
\se^{\prime \prime} \sq{n - 1} \mid \se \ri \}\\
& \leq  &  \sum_{n = 3}^{\ka(\ze, \vro)} \Pr \li \{  \f{ \sq{n} (
\ovl{X}_n - \si \se) } { \wh{\si}_n } \leq - t_{n-1, \bs{\ba}} - \f{
\sq{n - 1} \si \se } { \wt{\si}_n } + \se \sq{n - 1}
\mid \se \ri \}\\
& = &  \sum_{n = 3}^{\ka(\ze, \vro)} \Pr \li \{ \f{ \sq{n} (
\ovl{X}_n - \si \se) } { \wh{\si}_n } \leq - t_{n-1, \sq{\bs{\ba}} } + \vDe_n \mid \se \ri \}\\
& \leq & \sum_{n = 3}^{\ka(\ze, \vro)} \Pr \li \{ \f{ \sq{n} (
\ovl{X}_n - \si \se) } { \wh{\si}_n } \leq - t_{n-1, \sq{\bs{\ba}}
} \mid \se \ri \} + \sum_{n = 3}^{\ka(\ze, \vro)}  \Pr \{ \vDe_n  \geq 0 \mid \se  \}\\
& \leq & \ka(\ze, \vro) \; \sq{\bs{\ba}} + \sum_{n = 3}^{\ka(\ze,
\vro)}  \Pr \{ \vDe_n \geq 0 \mid \se  \}. \eee} Clearly,
\[
\ka(\ze, \vro) \; \sq{\bs{\ba}} \leq \vro \;
\mcal{Z}_{\sq{\bs{\ba}}}^2  \times  \sq{\bs{\ba}}  = \vro \times \f{
\mcal{Z}_{\sq{\bs{\ba}}}^2 }{ 2 \ln \f{1}{ \sq{\bs{\ba}} } } \times
\sq{\bs{\ba}}  \times  { 2 \ln \f{1}{ \sq{\bs{\ba}} } } \to 0
\]
as $\ze \to 0$.  Hence, to show (\ref{shw2}), it suffices to show
$\lim_{\ze \to 0 }\sum_{n = 3}^{\ka(\ze, \vro)}  \Pr \{ \vDe_n \geq
0 \mid \se  \} = 0$ for $\se < 0$ and $\se
> 0$.

In the case of  $\se < 0$,  we have {\small \[ \Pr \{ \vDe_n \geq 0
\mid \se \} = \Pr \li \{ \f{\si}{\wt{\si}_n} - 1 \geq \f{1}{|\se|}
\li ( \f{ t_{n-1, \bs{\ba}} } { \sq{n - 1} } - \f{ t_{n-1,
\sq{\bs{\ba}} } }{ \sq{n - 1} } \ri )  \ri \} \leq  \Pr \li \{
\f{\si}{\wt{\si}_n} \geq \f{1}{|\se|} \li ( \f{ t_{n-1, \bs{\ba}} }
{\sq{n - 1} } - \f{ t_{n-1, \sq{\bs{\ba}} } }{\sq{n - 1} } \ri ) \ri
\}. \]} By Lemma \ref{make8}, for small enough $\ze
> 0$, we have
\[
\Pr \{ \vDe_n \geq 0  \mid \se \} \leq \Pr \li \{
\f{\si}{\wt{\si}_n} > \exp \li ( \f{ \ln \f{1}{\bs{\ba} } }{ 4 (n -
1) } \ri ) \ri \} = \Pr \li \{ \f{\si}{\wh{\si}_n} > \sq{ \f{n -
1}{n} } \exp \li ( \f{ \ln \f{1}{\bs{\ba} } }{ 4 (n - 1) } \ri ) \ri
\}
\]
for $3 \leq n \leq \ka(\ze, \vro)$.  By Lemma \ref{lemch}, we have
that $\ln \f{1}{\sq{\bs{\ba}}} > \f{1}{3}
\mcal{Z}_{\sq{\bs{\ba}}}^2$ if $\ze$ is small enough.  This implies
that
\[ \exp \li ( \f{ \ln \f{1}{\sq{ \bs{\ba}} } }{ 6 (n - 1) } \ri )
> \exp \li ( \f{ \mcal{Z}_{\sq{\bs{\ba}}}^2 }{ 18 (n - 1) } \ri )
> \exp \li ( \f{ \mcal{Z}_{\sq{\bs{\ba}}}^2 }{ 18 n } \ri ) \geq \exp
\li ( \f{1}{18 \vro}  \ri )
\]
and thus \bee \sq{ \f{n - 1}{n} } \exp \li ( \f{ \ln \f{1}{\bs{\ba}
} }{ 4 (n - 1) } \ri ) & > & \sq{ \f{2}{3} } \exp \li ( \f{ \ln
\f{1}{\bs{\ba} } }{ 4 (n - 1) } \ri ) = \sq{ \f{2}{3} } \exp \li (
\f{ \ln \f{1}{\sq{ \bs{\ba}} } }{ 6 (n - 1) } \ri ) \exp \li (
\f{ \ln \f{1}{\bs{\ba} } }{ 6 (n - 1) } \ri )\\
& > & \sq{ \f{2}{3} } \exp \li ( \f{1}{18 \vro}  \ri ) \exp \li (
\f{ \ln \f{1}{\bs{\ba} } }{ 6 (n - 1) } \ri )\\
& > & \sq{ \f{2}{3} } \exp \li ( \f{1}{18 \times \f{1}{6} }  \ri )
\exp \li ( \f{ \ln \f{1}{\bs{\ba} } }{ 6 (n - 1) } \ri ) > \exp \li
( \f{ \ln \f{1}{\bs{\ba} } }{ 6 (n - 1) } \ri )
 \eee for $3 \leq n \leq
\ka(\ze, \vro)$ if $\ze$ is small enough, where we have used the
assumption that $\vro < \f{1}{6 (1 + |\se|)^2} < \f{1}{6}$.
Therefore, for small enough $\ze
> 0$, we have
\[
\Pr \{ \vDe_n \geq 0  \mid \se \} < \Pr \li \{ \f{\si}{\wh{\si}_n}
> \exp \li ( \f{ \ln \f{1}{\bs{\ba} } }{ 6 (n - 1) } \ri ) \ri \}
\]
for $3 \leq n \leq \ka(\ze, \vro)$ and it follows that \bel \Pr \{
\vDe_n \geq 0 \mid \se \} & < & \Pr \li \{ \wh{\si}_n < \si
\bs{\ba}^{\f{1}{6 (n-1)}} \mid \se \ri \} = \Pr \li \{ \chi_{n-1}^2
< (n - 1) \; \bs{\ba}^{\f{1}{3(n-1)}} \ri \} \nonumber\\
& \leq & \li [ \bs{\ba}^{\f{1}{3(n-1)}} \; \exp \li ( 1 -
\bs{\ba}^{\f{1}{3(n-1)}} \ri ) \ri ]^{(n - 1) \sh 2} <
\bs{\ba}^{\f{1}{6}} \; e^{(n - 1) \sh 2} \la{use3} \eel for $3 \leq
n \leq \ka(\ze, \vro)$.  Noting that $\f{1}{2} \ka(\ze, \vro) < 2
\vro \ln \f{1}{ \sq{\bs{\ba}} }$ for small enough $\ze$ and invoking
the assumption that $0 < \vro < \f{1}{6 (1 + |\se|)^2} < \f{1}{6}$,
we have \be \la{use4}
 \bs{\ba}^{\f{1}{6}} \exp \li ( \f{\ka(\ze, \vro)}{2} \ri ) < \bs{\ba}^{\f{1}{6}}
\exp \li ( 2 \vro \ln \f{1}{ \sq{\bs{\ba}} } \ri ) =
\bs{\ba}^{\f{1}{6} - \vro} \to 0 \ee as  $\ze \to 0$. It follows
from (\ref{use3}) and (\ref{use4})  that, in the case of $\se < 0$,
\be \la{exm}
 \sum_{n = 3}^{\ka(\ze, \vro)} \Pr \{ \vDe_n \geq 0 \mid
\se \}
 < \bs{\ba}^{\f{1}{6}} \sum_{n = 3}^{\ka(\ze, \vro)} e^{(n - 1) \sh 2} = \bs{\ba}^{\f{1}{6}} \times
\f{ \exp \li ( \f{\ka(\ze, \vro)}{2} \ri ) - e } {\sq{e} - 1} < 2 \;
\bs{\ba}^{\f{1}{6}}  \exp \li ( \f{\ka(\ze, \vro)}{2} \ri )
 \to 0
\ee as $\ze \to 0$.

In the case of $\se > 0$, by virtue of Lemma \ref{make8}, we have
{\small \bee \Pr \{ \vDe_n \geq 0 \mid \se \} = \Pr \li \{ 1 -
\f{\si}{\wt{\si}_n}  \geq \f{1}{\se} \li ( \f{ t_{n-1, \bs{\ba}} } {
\sq{n - 1} } - \f{ t_{n-1, \sq{\bs{\ba}} } }{ \sq{n - 1} } \ri )
\mid \se  \ri \} \leq \Pr \li \{ 1 - \f{\si}{\wt{\si}_n} >  1 \ri \}
= \Pr \{ \wt{\si}_n < 0 \} = 0 \eee} for $3 \leq n \leq \ka(\ze,
\vro)$ provided that $\ze$ is small enough.  It follows that
$\sum_{n = 3}^{\ka(\ze, \vro)} \Pr \{ \vDe_n \geq 0 \mid \se \} = 0$
for $\se > 0$ if $\ze > 0$ is sufficiently small.  Therefore, we
have shown that (\ref{shw1}) holds for $\se \geq \se^{\prime
\prime}$.

Next, we shall show that \[ \lim_{\ze \to 0} \sum_{n = \ka(\ze,
\vro) + 1}^\iy \Pr \li \{ \f{ \ovl{X}_n } { \wt{\si}_n } \leq \f{
\se^\prime + \se^{\prime \prime} }{2} + \f{ t_{n-1, \bs{\al}} -
t_{n-1, \bs{\ba}} } { 2 \sq{n - 1} } \mid \se \ri \}  = 0
\]
for $\se \geq \se^{\prime \prime}$.  By Lemma \ref{keyy}, there
exist a number $\vDe > 0$ and $\ze^* \in (0, 1)$ such that
\[
\f{ \se^\prime + \se^{\prime \prime} }{2} + \f{ t_{n-1, \bs{\al}} -
t_{n-1, \bs{\ba}} } { 2 \sq{n - 1} } < \se - \vDe, \qqu \fa \se \geq
\se^{\prime \prime}
\]
for any $\ze \in (0, \ze^*)$.   It follows from Lemma \ref{lemgood9}
that \bel &  & \sum_{n = \ka(\ze, \vro) + 1}^\iy \Pr \li \{ \f{
\ovl{X}_n } { \wt{\si}_n } \leq \f{ \se^\prime + \se^{\prime \prime}
}{2} + \f{ t_{n-1, \bs{\al}} - t_{n-1, \bs{\ba}} } { 2 \sq{n - 1} }
\mid \se \ri \} \nonumber\\
&  \leq  & \sum_{n = \ka(\ze, \vro) + 1}^\iy \Pr \li \{ \f{
\ovl{X}_n } { \wt{\si}_n } \leq \se - \vDe \mid \se \ri \} < \sum_{n
= \ka(\ze, \vro) + 1}^\iy \Pr \li \{ \li | \f{ \ovl{X}_n } {
\wt{\si}_n } - \se \ri | \geq \vDe \mid \se \ri \} \to 0 \la{combq}
\eel as $\ze \to 0$.  Combining (\ref{shw1}) and (\ref{combq}) leads
to (\ref{goal1}).

Now we want to show that (\ref{goal2}) is true.  It suffices to show
that \be \la{gog1} \lim_{\ze \to 0} \sum_{n = 2}^{\ka(\ze, \vro)}
\Pr \li \{ \f{ \ovl{X}_n } { \wt{\si}_n } \geq \se^{\prime} + \f{
t_{n-1, \bs{\al}} } { \sq{n - 1} } \mid \se \ri \} = 0 \ee and \be
\la{gog2}
 \lim_{\ze \to 0} \sum_{n = \ka(\ze, \vro) + 1}^\iy \Pr \li
\{ \f{ \ovl{X}_n } { \wt{\si}_n } \geq \f{ \se^\prime + \se^{\prime
\prime} }{2} + \f{ t_{n-1, \bs{\al}} - t_{n-1, \bs{\ba}} } { 2 \sq{n
- 1} } \mid \se \ri \}  = 0 \ee for $\se \leq \se^{\prime \prime}$
under the assumption that $0 < \vro < \f{1}{6 (1 + |\se|)^2}$.
Clearly, for $n = 2$ and $\se \leq \se^{\prime \prime}$, $\Pr \li \{
\f{ \ovl{X}_n } { \wt{\si}_n } \geq \se^{\prime} + \f{ t_{n-1,
\bs{\al}} } { \sq{n - 1} } \mid \se \ri \} \to 0$ as $\ze \to 0$.
Hence, to show (\ref{gog1}), it suffices to show that \be \la{gog3}
\lim_{\ze \to 0} \sum_{n = 3}^{\ka(\ze, \vro)} \Pr \li \{ \f{
\ovl{X}_n } { \wt{\si}_n } \geq \se^{\prime} + \f{ t_{n-1, \bs{\al}}
} { \sq{n - 1} } \mid \se \ri \} = 0 \ee for $\se \leq \se^\prime$.
 We can show (\ref{gog3}) by considering three cases: (i) $\se < 0$;
(ii) $\se > 0$; (iii) $\se = 0$.

Note that, for $\se \leq \se^\prime$, {\small \bee  \sum_{n =
3}^{\ka(\ze, \vro)} \Pr \li \{ \f{ \ovl{X}_n } { \wt{\si}_n } \geq
\se^{\prime} + \f{ t_{n-1, \bs{\al}} } { \sq{n - 1} } \mid \se \ri
\} & =  & \sum_{n = 3}^{\ka(\ze, \vro)} \Pr \li \{  \f{ \sq{n} \;
\ovl{X}_n } { \wh{\si}_n } \geq  t_{n-1, \bs{\al}} +
\se^{\prime} \sq{n - 1} \mid \se \ri \}\\
& \leq  &  \sum_{n = 3}^{\ka(\ze, \vro)} \Pr \li \{  \f{ \sq{n} (
\ovl{X}_n  - \si \se ) } { \wh{\si}_n } \geq  t_{n-1, \bs{\al}} -
\f{ \sq{n - 1} \si \se } { \wt{\si}_n } +  \se \sq{n - 1}
\mid \se \ri \}\\
& = &  \sum_{n = 3}^{\ka(\ze, \vro)}
 \Pr \li \{ \f{ \sq{n} (\ovl{X}_n - \si \se  ) } { \wh{\si}_n }
 \geq  t_{n-1, \sq{\bs{\al}} } + \vDe_n \mid \se \ri \}\\
& \leq & \sum_{n = 3}^{\ka(\ze, \vro)} \Pr \li \{ \f{ \sq{n} (
\ovl{X}_n  - \si \se ) } { \wh{\si}_n } \geq  t_{n-1, \sq{\bs{\al}}
} \mid \se \ri \} + \sum_{n = 3}^{\ka(\ze, \vro)}  \Pr \{ \vDe_n  \leq 0 \mid \se  \}\\
& < & \ka(\ze, \vro) \; \sq{\bs{\al}} + \sum_{n = 3}^{\ka(\ze,
\vro)} \Pr \{ \vDe_n \leq 0 \mid \se  \}, \eee} where
\[
\vDe_n = \se \sq{n - 1} \li ( 1 - \f{\si}{\wt{\si}_n} \ri ) -
t_{n-1, \sq{\bs{\al}} } +
 t_{n-1, \bs{\al}}, \qqu n = 3, 4, \cd
\]
and
\[
\ka(\ze, \vro) \; \sq{\bs{\al}} \leq \vro \;
\mcal{Z}_{\sq{\bs{\al}}}^2  \times  \sq{\bs{\al}}  = \vro \times \f{
\mcal{Z}_{\sq{\bs{\al}}}^2 }{ 2 \ln \f{1}{ \sq{\bs{\al}} } } \times
\sq{\bs{\al}}  \times  { 2 \ln \f{1}{ \sq{\bs{\al}} } } \to 0
\]
as $\ze \to 0$.

In the case of $\se > 0$,  by Lemma \ref{make8}, we have \bee \Pr \{
\vDe_n \leq 0 \mid \se \} & = & \Pr \li \{ 1 - \f{\si}{\wt{\si}_n}
\leq - \f{1}{\se} \li ( \f{ t_{n-1, \bs{\al}} } {\sq{n - 1} } - \f{
t_{n-1, \sq{\bs{\al}} } }{ \sq{n -
1} } \ri ) \ri \}\\
& \leq &   \Pr \li \{ \f{\si}{\wt{\si}_n} \geq  \f{1}{\se} \li ( \f{
t_{n-1, \bs{\al}} } {\sq{n - 1} } - \f{ t_{n-1, \sq{\bs{\al}} }
}{\sq{n - 1} }  \ri ) \ri \} \leq  \Pr \li \{ \f{\si}{\wt{\si}_n}
> \exp \li ( \f{ \ln \f{1}{\bs{\al}} }{4 ( n - 1 )} \ri ) \ri \}
\eee for $3 \leq n \leq \ka(\ze, \vro)$ if  $\ze$ is small enough.
Hence, by a similar method as that for proving (\ref{exm}), we have
$\lim_{\ze \to 0} \sum_{n = 3}^{\ka(\ze, \vro)} \Pr \{ \vDe_n \leq 0
\mid \se \} \to 0$ as $\ze \to 0$.

 In the case of $\se < 0$, by Lemma \ref{make8}, we have \bee \Pr \{ \vDe_n  \leq 0 \mid
\se \} = \Pr \li \{ 1 - \f{\si}{\wt{\si}_n} \geq \f{1}{|\se|} \li (
\f{ t_{n-1, \bs{\al}} } { \sq{n - 1} } - \f{ t_{n-1, \sq{\bs{\al}} }
}{ \sq{n - 1} } \ri ) \ri \} \leq \Pr \li \{ 1 - \f{\si}{\wt{\si}_n}
> 1 \ri \} = 0 \eee for $3 \leq n \leq \ka(\ze, \vro)$ if $\ze$ is small
enough.  Hence, $\lim_{\ze \to 0} \sum_{n = 3}^{\ka(\ze, \vro)} \Pr
\{ \vDe_n \leq 0 \mid \se \} = 0$ for $\se < 0$ if $\ze > 0$ is
small enough.

In the case of $\se = 0 \leq \se^\prime$, we have {\small \bee
\sum_{n = 3}^{\ka(\ze, \vro)} \Pr \li \{ \f{ \ovl{X}_n } {
\wt{\si}_n } \geq \se^{\prime} + \f{ t_{n-1, \bs{\al}} } { \sq{n -
1} } \mid \se \ri \}  \leq   \ka(\ze, \vro) \bs{\al} \to 0 \eee} as
$\ze \to 0$.  Therefore, (\ref{gog3}) is true for all three cases.
As a result, (\ref{gog1}) is true for $\se \leq \se^\prime$.

By a similar method as that for (\ref{combq}), we can show that
(\ref{gog2}) is true.  Finally,  combining (\ref{gog1}) and
(\ref{gog2}) leads to (\ref{goal2}).  This completes the proof of
the lemma.

\epf

\bsk

Now we are in a position to prove the theorem.  Note that \be
\la{impa0} \Pr \{ \tx{Reject} \; \mscr{H}_j \mid \se \} \leq \sum_{i
= 1}^{j} \Pr \{ \tx{Accept} \; \mscr{H}_{i-1} \mid \se \} + \sum_{i
= j + 1}^{m - 1} \Pr \{ \tx{Accept} \; \mscr{H}_i \mid \se \}. \ee
By Lemma \ref{lemch1}, we have
\[
f_{\ell, i} \leq \se_{i}^{\prime \prime}, \qqu g_{\ell, i} \geq
\se_{i}^{\prime}, \qqu f_{\ell, i} \leq \f{ \se_{i}^{\prime} +
\se_{i}^{\prime \prime} }{2} + \f{ t_{n_\ell - 1, \al_i} - t_{n_\ell
- 1, \ba_i} } {2 \sq{n_\ell - 1} } \leq g_{\ell, i}
\]
for $i = 1, \cd, m - 1$ and $\ell = 1, \cd, s$.  Hence, by the
definition of the testing plan, we have \bel \Pr \{ \tx{Accept} \;
\mscr{H}_{i-1} \mid \se \}  & < &  \sum_{n = 2}^{\ka} \Pr \li \{ \f{
\ovl{X}_n } { \wt{\si}_n } \leq \se_{i}^{\prime
\prime} - \f{  t_{n - 1, \ba_i} } { \sq{n - 1} } \mid \se \ri \} \nonumber\\
&    &  + \sum_{n = \ka + 1}^\iy \Pr \li \{ \f{ \ovl{X}_n } {
\wt{\si}_n } \leq \f{ \se_i^{\prime} + \se_i^{\prime \prime} }{2} +
\f{ t_{n - 1, \al_i} - t_{n - 1, \ba_i} } {2 \sq{n - 1} } \mid \se
\ri \} \la{goodg} \eel for $i = 1, \cd, m$, where $\ka$ can be any
integer greater than $2$.  Making use of (\ref{goodg}) and applying
Lemma \ref{lem11m} with  $\ka = \ka(\ze, \vro)$, we have that \be
\la{impa} \lim_{\ze \to 0} \Pr \{ \tx{Accept} \; \mscr{H}_{i-1} \mid
\se \} = 0, \qqu \fa \se \geq \se_i^{\prime \prime},  \qqu i = 1,
\cd, m. \ee

Similarly, by the definition of the testing plan, we have \bel \Pr
\{ \tx{Accept} \; \mscr{H}_i \mid \se \}  & < &  \sum_{n = 2}^{\ka}
\Pr \li \{ \f{ \ovl{X}_n } { \wt{\si}_n } \geq
\se_{i}^\prime + \f{  t_{n - 1, \al_i} } { \sq{n - 1} } \mid \se \ri \} \nonumber\\
&    &  + \sum_{n = \ka + 1}^\iy \Pr \li \{ \f{ \ovl{X}_n } {
\wt{\si}_n } \geq \f{ \se_i^{\prime} + \se_i^{\prime \prime} }{2} +
\f{ t_{n - 1, \al_i} - t_{n - 1, \ba_i} } {2 \sq{n - 1} } \mid \se
\ri \} \la{goodg2} \eel for $i = 1, \cd, m - 1$, where $\ka$ can be
any integer greater than $2$.  Making use of (\ref{goodg2}) and
applying Lemma \ref{lem11m} with $\ka = \ka(\ze, \vro)$, we have
that \be \la{impa2} \lim_{\ze \to 0} \Pr \{ \tx{Accept} \;
\mscr{H}_i \mid \se \} = 0, \qqu \fa \se \leq \se_i^\prime,  \qqu i
= 1, \cd, m -1. \ee Therefore, Theorem
\ref{Normal_Mean_unknown_variance} follows from (\ref{impa0}),
(\ref{impa}) and (\ref{impa2}).

\sect{Proofs of Theorems \ref{Ratio_known_mean} and
\ref{Ratio_unknown_mean} } \la{proof ratio}

As a consequence of the definitions of the sampling schemes,
Theorems \ref{Ratio_known_mean} and \ref{Ratio_unknown_mean} can be
proved by the same argument, which relies on a preliminary result as
stated by the following lemma.

\beL \la{FChernoff}
 Let $Z$ be a random variable possessing an
 $F$-distribution of $m$ and $n$ degrees of freedom.
Then, for $r$ greater than $1$, both $\Pr \li \{ Z > r \ri \}$ and
$\Pr \li \{ Z < \f{1}{r} \ri \}$ are less than $2 [g(r)]^d$, where
$d = \min (m, n) / 2$ and $g(x) = \f{1}{\sq{x}} \exp ( 1 -
\f{1}{\sq{x}}  )$.  \eeL

\bpf Clearly,  $Z$ can be expressed as $\f{U}{V}$, where $U$ and $V$
are independent random variables possessing
 $\chi^2$-distributions of $m$ and $n$ degrees of freedom
respectively.  Note that $\Pr \li \{ Z > r \ri \} = \Pr \li \{
\f{U}{V} > r \ri \} \leq \Pr \{ U > \sq{r} \} + \Pr  \{ V <
\f{1}{\sq{r}} \} < [ g ( \f{1}{\sq{r}} )  ]^{m \sh 2} + [ g(r) ]^{n
\sh 2}$,  where the last inequality follows from Chernoff bound.
Observing that $g ( \f{1}{\sq{r}} ) = g(r) = 1$ for $r = 1$ and that
the derivative of {\small $\f{ g ( \f{1}{\sq{r}} ) } { g(r) }$} with
respective to $r$ is negative for $r$ greater than $1$, we have $g (
\f{1}{\sq{r}} ) < g(r)$ for $r > 1$. It follows that $\Pr \li \{ Z >
r \ri \} < 2 [g(r)]^d$ for $r > 1$. Since $\f{1}{Z}$ is a random
variable possessing an $F$-distribution of $n$ and $m$ degrees of
freedom, it follows from the established result that $\Pr \{
\f{1}{Z} > r \} < 2 [g(r)]^d$ for $r > 1$. This completes the proof
of the lemma.

\epf

We are now in a position to prove the theorems. Let \[ l^* = 1 + 2
\max_{ i \in \{1, \cd, m - 1 \}} \max \li \{ \f{ \ln \f{\al_i}{2} }
{ \ln g (k_i) }, \f{ \ln \f{\ba_i }{2} } { \ln g (k_i) } \ri \} = O
\li ( \ln \f{1}{\ze} \ri ), \] where $k_i = \sq{ \f{ \se_{i}^{\prime
\prime} } { \se_{i}^{\prime} }}$ for $i = 1, \cd, m - 1$. Then, \be
\la{vip} \f{1}{2} \li [ \min \{ n_\ell^X, n_\ell^Y \} - 1 \ri ] \geq
\f{\ell - 1}{2} \geq \max_{ i \in \{1, \cd, m - 1 \}} \max \li \{
\f{ \ln \f{\al_i}{2} } { \ln g (k_i)  }, \f{ \ln \f{\ba_i }{2} } {
\ln g (k_i) } \ri \} \ee for $\ell \geq l^*$. Making use of
(\ref{vip}) and Lemma \ref{FChernoff}, we have $\Pr \{
\wh{\bs{\se}}_\ell > k_i \se \} < \al_i, \; \Pr \{
\wh{\bs{\se}}_\ell < \f{\se}{k_i} \} < \ba_i$ and consequently $\Up
(n_\ell^X - 1, n_\ell^Y - 1, 1 - \al_i) < k_i, \; \Up (n_\ell^X - 1,
n_\ell^Y - 1, \ba_i) > \f{1}{k_i}$ for $i = 1, \cd, m-1$.  This
implies that $\f{1}{k_i} < v_{\ell, i} \leq u_{\ell, i} < k_i$  and
thus $\se_{i}^{\prime \prime} v_{\ell, i} \geq \se_{i}^{\prime}
u_{\ell, i}$ for
 $i = 1, \cd, m-1$.  It follows that $\{ \bs{l} \leq l^* \}$ is a sure
 event and consequently, for any $\se \in \varTheta_i$ and $i = 0,
 1, \cd, m - 1$,
 \[
\Pr \{ \tx{Reject} \; \mscr{H}_i \mid \se \} \leq \sum_{\ell =
1}^{l^*} \sum_{j = 1}^{m - 1} (\al_j + \ba_j) = O \li ( \ln
\f{1}{\ze} \ri ) O( \ze ) \to 0
 \]
as $\ze \to 0$.  This completes the proof of the theorems.

\sect{Proof of Theorem \ref{Exact_Computation_SPRT}} \la{Exact_Computation_SPRT_app}

Under the assumption that $u(m) > \se$ and $\li [ \f{ \exp ( \eta (\se) u  - \psi (\se) ) } { \exp ( \eta (u) u - \psi (u) ) } \ri ]^m < \ep$,
making use of the definition of the stopping rule of SPRT and the likelihood ratio bound established in \cite{ChenR}, we have
\[
\Pr \{ \mbf{n} > m \mid \se \} \leq \Pr \li \{    \f{\sum_{i = 1}^m X_i}{m} > u(m)  \mid \se \ri \} \leq \li [ \f{ \exp ( \eta (\se) u  - \psi
(\se) ) } { \exp ( \eta (u) u - \psi (u) ) } \ri ]^m < \ep. \]  It follows that
\[
\Pr \{  \tx{Accept $\mscr{H}_0$}, \; \mbf{n} \leq m   \mid \se \} \leq \Pr \{  \tx{Accept $\mscr{H}_0$} \mid \se \} \leq \Pr \{  \tx{Accept
$\mscr{H}_0$}, \; \mbf{n} \leq m \mid \se \} + \ep,
\]
\[
\Pr \{  \tx{Accept $\mscr{H}_1$}, \; \mbf{n} \leq m   \mid \se \} \leq \Pr \{  \tx{Accept $\mscr{H}_1$} \mid \se \} \leq \Pr \{  \tx{Accept
$\mscr{H}_1$}, \; \mbf{n} \leq m \mid \se \} + \ep.
\]
Similarly, under the assumption that $v(m) < \se$ and $\li [ \f{ \exp ( \eta (\se) v  - \psi (\se) ) } { \exp ( \eta (v) v - \psi (v) ) } \ri
]^m < \ep$, we have
\[
\Pr \{ \mbf{n} > m \} \leq \Pr \li \{    \f{\sum_{i = 1}^m X_i}{m} < v(m)  \ri \}  \leq  \li [ \f{ \exp ( \eta (\se) v  - \psi (\se) ) } { \exp
( \eta (v) v - \psi (v) ) } \ri ]^m < \ep \] and thus
\[
\Pr \{  \tx{Accept $\mscr{H}_0$}, \; \mbf{n} \leq m   \mid \se \} \leq \Pr \{  \tx{Accept $\mscr{H}_0$} \mid \se \} \leq \Pr \{  \tx{Accept
$\mscr{H}_0$}, \; \mbf{n} \leq m \mid \se \} + \ep,
\]
\[
\Pr \{  \tx{Accept $\mscr{H}_1$}, \; \mbf{n} \leq m   \mid \se \} \leq \Pr \{  \tx{Accept $\mscr{H}_1$} \mid \se \} \leq \Pr \{  \tx{Accept
$\mscr{H}_1$}, \; \mbf{n} \leq m \mid \se \} + \ep.
\]

Now we consider the computation of $\bb{E} [ \mbf{n} ]$.  Note that $\bb{E} [ \mbf{n} ] = \sum_{n = 1}^\iy \Pr \{  \mbf{n} > n  \}$.  It
suffices to bound $\sum_{n = m}^\iy \Pr \{  \mbf{n} > n  \}$.  By the assumption that $u(m)
> \se$, we have $u(n) > \se$ for all $n \geq m$, since $u(n)$ is increasing with respect to $n$.  It follows that
\bee \Pr \{ \mbf{n} > n \} & \leq & \Pr \li \{    \f{\sum_{i = 1}^n X_i}{n} > u(n)  \mid \se \ri \} \\
& \leq & \li [ \f{ \exp ( \eta (\se) u (n)  - \psi (\se) ) } { \exp ( \eta (u(n)) u(n) - \psi (u(n)) ) } \ri ]^n\\
& \leq & \li [ \f{ \exp ( \eta (\se) u (m)  - \psi (\se) ) } { \exp ( \eta (u(m)) u(m) - \psi (u(m)) ) } \ri ]^n  \eee for $n \geq m$ and
consequently,
\[
\sum_{n = m}^\iy \Pr \{  \mbf{n} > n  \} \leq \sum_{n = m}^\iy \li [ \f{ \exp ( \eta (\se) u  - \psi (\se) ) } { \exp ( \eta (u) u - \psi (u) )
} \ri ]^n  = \f{ \li [ \f{ \exp ( \eta (\se) u - \psi (\se) ) } { \exp ( \eta (u) u - \psi (u) ) } \ri ]^m } { 1 - \f{ \exp ( \eta (\se) u -
\psi (\se) ) } { \exp ( \eta (u) u - \psi (u) ) } } < \ep
\]
where $u = u(m)$, provided that $u(m)
> \se$ and $\li [ \f{ \exp ( \eta (\se) u - \psi (\se) ) } { \exp ( \eta (u) u - \psi (u) ) } \ri ]^m < \ep \li [ 1 - \f{ \exp ( \eta (\se) u -
\psi (\se) ) } { \exp ( \eta (u) u - \psi (u) ) } \ri ]$.

Similarly, by the assumption that $v(m) < \se$, we have $v(n) < \se$ for all $n \geq m$, since $v(n)$ is decreasing with respect to $n$.  It
follows that
\bee \Pr \{ \mbf{n} > n \} & \leq & \Pr \li \{    \f{\sum_{i = 1}^n X_i}{n} < v(n)  \mid \se \ri \} \\
& \leq & \li [ \f{ \exp ( \eta (\se) v (n)  - \psi (\se) ) } { \exp ( \eta (v(n)) v(n) - \psi (v(n)) ) } \ri ]^n\\
& \leq & \li [ \f{ \exp ( \eta (\se) v (m)  - \psi (\se) ) } { \exp ( \eta (v(m)) v(m) - \psi (v(m)) ) } \ri ]^n  \eee for $n \geq m$ and
consequently,
\[
\sum_{n = m}^\iy \Pr \{  \mbf{n} > n  \} \leq \sum_{n = m}^\iy \li [ \f{ \exp ( \eta (\se) v  - \psi (\se) ) } { \exp ( \eta (v) v - \psi (v) )
} \ri ]^n  = \f{ \li [ \f{ \exp ( \eta (\se) v - \psi (\se) ) } { \exp ( \eta (v) v - \psi (v) ) } \ri ]^m } { 1 - \f{ \exp ( \eta (\se) v -
\psi (\se) ) } { \exp ( \eta (v) v - \psi (v) ) } } < \ep
\]
where $v = v(m)$, provided that $v(m) < \se$ and $\li [ \f{ \exp ( \eta (\se) v - \psi (\se) ) } { \exp ( \eta (v) v - \psi (v) ) } \ri ]^m <
\ep \li [ 1 - \f{ \exp ( \eta (\se) v - \psi (\se) ) } { \exp ( \eta (v) v - \psi (v) ) } \ri ]$.  This completes the proof of the theorem.

\sect{Proof of Theorem \ref{recursive}} \la{recursive_app}

 Note that inequality (\ref{assign}) can be written as \be
\la{ineq189}
 \f{ 1 + \vep_\ell } { 1 + \vep_{\ell + 1}  } (1 +
\eta_\ell) (1 + \ga_\ell) \leq \f{ 1 - \vep_\ell }{ 1 - \vep_{\ell + 1} } (1 - \eta_\ell). \ee By  virtue of (\ref{con4}) and (\ref{ineq189}),
we have \bel \f{1}{2} \f{ 1 + \vep_\ell }{1 + \vep_{\ell + 1} } (1 + \eta_\ell ) \li ( \ovl{L}_\ell + \ovl{U}_\ell \ri ) & < & \f{1}{2} \f{ 1 +
\vep_\ell }{1 + \vep_{\ell + 1} }  (1 + \eta_\ell
) (1 + \ga_\ell) \li ( \udl{L}_\ell + \udl{U}_\ell \ri ) \nonumber\\
& \leq &  \f{1}{2} \f{ 1 - \vep_\ell }{1 - \vep_{\ell + 1} }  (1 - \eta_\ell ) \li ( \udl{L}_\ell + \udl{U}_\ell  \ri ). \la{cta} \eel Making
use of (\ref{con1}) and the assumption (\ref{imp88}),  we have \be \la{ctb}
 g_{\ell+1} (y) \leq (1 + \vep_\ell)
\sum_{i=1}^{m_\ell} h_{\ell, i} \; \ovl{I} (A_{\ell,i}, B_{\ell, i}, C, D) < \f{1}{2} (1 + \vep_\ell) (1 + \eta_\ell) \li ( \ovl{L}_\ell
 + \ovl{U}_\ell \ri ). \ee  Similarly, making
use of (\ref{con2}) and the assumption (\ref{imp88}),  we have \be \la{ctc}
 g_{\ell+1} (y)  \geq  (1 - \vep_\ell) \sum_{i=1}^{m_\ell}
h_{\ell, i} \; \udl{I} (A_{\ell,i}, B_{\ell, i}, C, D) >  \f{1}{2} (1 - \vep_\ell) (1 - \eta_\ell) \li ( \udl{L}_\ell + \udl{U}_\ell \ri ). \ee
Combining (\ref{cta}), (\ref{ctb}) and (\ref{ctc}) leads to  $(1 - \vep_{\ell+1}) h_{\ell+1} < g_{\ell+1} (y) < (1 + \vep_{\ell+1}) h_{\ell+1}$
for any $y \in [C, \; D]$.

From (\ref{con1}) and (\ref{con2}),  we have {\small \bee &  & (1 + \ga_\ell) (\udl{L}_\ell + \udl{U}_\ell) >  \f{ 2 (1 + \ga_\ell) }{ 1 +
\eta_\ell } \sum_{i=1}^{m_\ell}  h_{\ell, i} \; \udl{I}
(A_{\ell,i}, B_{\ell, i}, C, D),\\
&  & \ovl{L}_\ell + \ovl{U}_\ell < \f{2}{1 - \eta_\ell} \sum_{i=1}^{m_\ell}  h_{\ell, i} \; \ovl{I} (A_{\ell,i}, B_{\ell, i}, C, D).  \eee}   By
the assumptions on $\udl{I}$ and $\ovl{I}$, we have that

{\small $| \sum_{i=1}^{m_\ell} h_{\ell, i} \; \ovl{I} (A_{\ell,i}, B_{\ell, i}, C, D) - \sum_{i=1}^{m_\ell} h_{\ell, i} \; \udl{I} (A_{\ell,i},
B_{\ell, i}, C, D) | \to 0$} as $D - C \to 0$. It follows that (\ref{con4}) is satisfied if $(1 + \ga_\ell) (1 - \eta_\ell)
> 1 + \eta_\ell$ and $D - C$ is sufficiently small.  This completes
the proof of the theorem.

\end{document}